\documentclass[a4paper,12pt]{article}

\usepackage{amsmath,amsthm,amscd}
\usepackage{amssymb}

\newtheorem{thm}{Theorem}[section]
\newtheorem{cor}[thm]{Corollary}
\newtheorem{lem}[thm]{Lemma}
\newtheorem{prop}[thm]{Proposition}
\newtheorem{defn}[thm]{Definition}
\newtheorem{ntt}[thm]{Notation}

\theoremstyle{definition}
\newtheorem{exmp}[thm]{Example}
\newtheorem{rem}[thm]{Remark}

\newcommand{\cC}{\mathcal{C}}

\newcommand{\cH}{\mathcal{H}}

\newcommand{\cJ}{\mathcal{J}}

\newcommand{\cL}{\mathcal{L}}

\newcommand{\cO}{\mathcal{O}}

\newcommand{\fa}{\mathfrak{a}}
\newcommand{\fb}{\mathfrak{b}}
\newcommand{\fd}{\mathfrak{d}}

\newcommand{\fm}{\mathfrak{m}}

\newcommand{\fp}{\mathfrak{p}}

\newcommand{\bB}{\mathbf{B}}

\newcommand{\ga}{\alpha}
\newcommand{\gb}{\beta}

\newcommand{\gd}{\delta}

\DeclareMathOperator{\ord}{ord}

\DeclareMathOperator{\Spec}{Spec}

\makeatletter
\def\overbracket#1{\mathop{\vbox{\ialign{##\crcr\noalign{\kern3\p@}
	\downbracketfill\crcr\noalign{\kern3\p@\nointerlineskip}
	$\hfil\displaystyle{#1}\hfil$\crcr}}}\limits}
\def\underbracket#1{\mathop{\vtop{\ialign{##\crcr
	$\hfil\displaystyle{#1}\hfil$\crcr\noalign{\kern3\p@\nointerlineskip}
	\upbracketfill\crcr\noalign{\kern3\p@}}}}\limits}
\def\overparenthesis#1{\mathop{\vbox{\ialign{##\crcr\noalign{\kern3\p@}
	\downparenthfill\crcr\noalign{\kern3\p@\nointerlineskip}
	$\hfil\displaystyle{#1}\hfil$\crcr}}}\limits}
\def\underparenthesis#1{\mathop{\vtop{\ialign{##\crcr
	$\hfil\displaystyle{#1}\hfil$\crcr\noalign{\kern3\p@\nointerlineskip}
	\upparenthfill\crcr\noalign{\kern3\p@}}}}\limits}
\def\downparenthfill{$\m@th\braceld\leaders\vrule\hfill\bracerd$}
\def\upparenthfill{$\m@th\bracelu\leaders\vrule\hfill\braceru$}
\def\upbracketfill{$\m@th\makesm@sh{\llap{\vrule\@height3\p@\@width.7\p@}}
  \leaders\vrule\@height.7\p@\hfill
  \makesm@sh{\rlap{\vrule\@height3\p@\@width.7\p@}}$}
\def\downbracketfill{$\m@th
  \makesm@sh{\llap{\vrule\@height.7\p@\@depth2.3\p@\@width.7\p@}}
  \leaders\vrule\@height.7\p@\hfill
  \makesm@sh{\rlap{\vrule\@height.7\p@\@depth2.3\p@\@width.7\p@}}$}
\makeatother
\begin{document}
\title
{JUMPING NUMBERS OF A SIMPLE COMPLETE IDEAL IN A TWO-DIMENSIONAL REGULAR 
LOCAL RING}
\author{\textsc{Tarmo J\"arvilehto}}
\maketitle
\vspace{30pt}
\thispagestyle{empty}
\tableofcontents
\newpage 

\section{Introduction}

In recent years multiplier ideals have emerged as an important tool in algebraic geometry and commutative algebra.
For the history of multiplier ideals and for more details concerning their general theory we refer to \cite{Lz}. 
Given an ideal $I$ in a regular local ring, its multiplier ideals  $\cJ(cI)$ form a family parametrized by non-negative 
rational numbers $c$. This family is totally ordered by inclusion and the parametrization is order reversing, i.e., 
$c<c'$ implies $\cJ(c'I) \subset \cJ(cI)$. The jumping numbers associated to the ideal $I$ are those rational numbers 
$\xi$ satisfying $\cJ(\xi I) \subsetneq \cJ((\xi-\varepsilon)I)$ for all $\varepsilon>0$. 
Jumping numbers of an ideal encode information about the singularities of the corresponding subscheme.
The first one of these numbers, the log-canonical threshold, has been much studied in birational geometry.
A good source for information about jumping numbers is the fundamental article \cite{ELSV} by 
Ein, Lazarsfeld, Smith and Varolin.

The purpose of the present manuscript is to determine the jumping numbers of a simple complete ideal $\fa$ in a two 
dimensional regular local ring $\ga$, and investigate their connection to other singularity invariants associated to 
the ideal. Recall that an ideal is called simple if it is not a product of two proper ideals. Simple ideals play a 
fundamental role in Zariski's theory of complete ideals, his famous theorem about the unique factorization of complete 
ideals saying that every complete ideal can be expressed uniquely as a product of simple ideals. Zariski thought of 
complete ideals as linear systems of curves satisfying "infinitely near base conditions". His theorem about the unique 
factorization of a complete ideal then corresponds to the factorization of a curve into irreducible branches. This leads 
us to consider the jumping numbers of an analytically irreducible plane curve. In particular, we are interested to 
compare the jumping numbers of a simple complete ideal to those of the analytically irreducible plane curve defined by 
a general element of the ideal.

Based on the proximity relations between infinitely near points, our methods are very much arithmetical in nature. 
The advantage of this approach is that tecniques needed are elementary. Moreover, our results hold in every 
characteristic.

Let us explain our results in more detail. We first define the notion of a log-canonical threshold of an ideal with 
respect to another ideal. It turns out in Proposition \ref{Ha} that jumping numbers are log-canonical thresholds with 
respect to suitable ideals. We then utilize the proximity relations to calculate these numbers 
(see Propositions \ref{ci} and \ref{c=RH}). In our main Theorem \ref{main} we give the promised formula for the jumping 
numbers. Consider the composition of point blow-ups 
\begin{equation}\label{resI}
\mathcal X=\mathcal X_{n+1}\xrightarrow{\pi_n}\cdots\xrightarrow{\pi_2}\mathcal X_2\xrightarrow{\pi_1}
\mathcal X_1=\Spec\ga,
\end{equation}         
obtained by blowing up the base points of $\fa$. It turns out that the set of jumping numbers is 
$\cH_\fa=H_0 \cup\cdots\cup H_{g^*}$, where the sets $H_0,\dots,H_{g^*}$ correspond to the stars of the dual graph 
arising from the configuration of the irreducible exceptional divisors on $\mathcal X$. 
A star is a vertex associated to an irreducible 
exceptional divisor which intersects more than two other irreducible exceptional divisors (see page \pageref{sdg}). 
More precisely, assume that $(a_1,\dots,a_n)$ is the point basis of $\fa$ (see page \pageref{pb}) and let 
$\{\gamma_1,\dots,\gamma_{g^*}\}$ be the set of indices corresponding to the star vertices. 
Set $\gamma_0=1$ and $\gamma_{g^*+1}=n$, and for every $\nu=0,\dots,g^*$ write 
\[
b_\nu:=\frac{a_1^2+\cdots+a_{\gamma_{\nu+1}}^2}{a_{\gamma_\nu}}.
\] 
Then 
\[
H_\nu:=\left\lbrace \frac{s+1}{a_{\gamma_\nu}}+\frac{t+1}{b_\nu}+\frac{m}{a_{\gamma_{\nu+1}}}
\hspace{2pt}\middle |\hspace{2pt} s,t,m\in\mathbb N, 
\frac{s+1}{a_{\gamma_\nu}}+\frac{t+1}{b_\nu}\le\frac{1}{a_{\gamma_{\nu+1}}}\right\rbrace
\]
for $\nu=0,\dots,g^*-1$, while
\[
H_{g^*}:=\left\lbrace\frac{s+1}{a_{\gamma_{g^*}}}+\frac{t+1}{b_{g^*}}
\hspace{2pt}\middle|\hspace{2pt} s,t\in\mathbb N\right\rbrace.
\]
We point out in Remark \ref{Za} that the numbers $a_{\gamma_0}, b_0,b_{1},\dots$ above are in fact 
the "Zariski exponents" of the ideal. They are conventionally denoted by $\bar\gb_\nu$ 
while the integers $a_{\gamma_\nu}$ are often denoted by $e_\nu$ for every $\nu$.

As a consequence of the main result, it turns out in Theorem \ref{1maincor} that the set of jumping numbers gives 
equivalent data to the information obtained from the point basis. The proof of this result is based on our 
Theorem \ref{ord}, which shows that the three smallest jumping numbers determine the order of the ideal.

In fact, to obtain the point basis from the jumping numbers we only need to know the elements 
$\xi'_\nu=1/a_{\gamma_\nu}+1/b_\nu=\min H_\nu$ in every set $H_0,\dots,H_{g^*}$. This is proved in Corollary 
\ref{maincor}. 
We shall observe in Proposition 
\ref{höl} that $\xi'_\nu=\min\{\xi\in\cH_\fa\mid\xi\ge1/a_{\gamma_\nu}\}$. Note that $\xi'_0$ is the log-canonical 
threshold of the ideal, which is already enough to determine the point basis and thereby all the jumping 
numbers, in the case the ideal is monomial. In a way, we may regard the sequence $(\xi'_0,\dots,\xi'_{g^*})$ as a 
generalization of the log-canonical threshold. In Corollary \ref{HS} we give also a formula $e(\fa)=(\xi'-1)^{-1}$, 
where $\xi':=\min\{\xi\in\cH_\fa\mid\xi>1\}$, for the Hilbert-Samuel multiplicity of the ideal $\fa$.

Utilizing the equivalence between simple complete ideals and analytically irreducible plane curves, we can now 
determine the jumping numbers of an analytically irreducible plane 
curve as well. 
It follows from \cite[Proposition 9.2.28]{Lz} (see also Proposition \ref{curveJ}) that 
the jumping numbers of our simple complete ideal $\fa$
coincide in the interval $[0,1[$ to those of the analytically irreducible plane curve corresponding to a "general" 
element of $\fa$. In fact, these numbers determine the jumping numbers of $\fa$ as soon 
as the integer $n$ appearing in the resolution (\ref{resI}) is known (see Theorem \ref{paxi} and Remark \ref{pox}). 
It is also worth to note that the word "general" can be interpreted here in the sense of Spivakovsky: Consider the 
resolution (\ref{resI}) above. Following \cite[Definition 7.1]{S1} and \cite[Definition 1]{CPR}, an element $f\in\fa$ 
is defined to be general, if the corresponding curve $\cC_f$ is analytically irreducible 
and the strict transform of $\cC_f$ intersects the strict transform of any
exceptional divisor passing through the center of $\pi_n$ transversely at this point.

Our formula for the jumping numbers of an arbitrary analytically irreducible plane curve in Theorem \ref{maincurve} shows 
that the jumping numbers depend only on the equisingularity class of the curve. Remarkably, it turns out in Theorem 
\ref{curveJh} that the jumping numbers actually determine the equisingularity class. Thus one can say that information 
about the topological type of the curve singularity is encoded in the set of the jumping numbers.

One should note that in the case of an analytically irreducible plane curve singularity over the complex numbers our 
formula can be obtained 
more directly using the general theory of jumping numbers. Indeed, in \cite[p. 1191]{V1} and \cite[p. 390]{V2} 
Vaqui\'e observed that the jumping numbers can be read off from the Hodge-theoretic spectrum defined by Steenbrink and 
Varchenko. This spectrum has in turn been calculated by several authors, starting from an unpublished preprint of M. 
Saito \cite[Theorem 1.5]{Sa}. One should also note that Igusa found a formula for the log-canonical threshold of an 
analytically irreducible plane curve with an isolated singularity (see \cite{I}), and that this result was generalized 
to reducible curves by Kuwata in \cite[Theorem 1.2]{K}. 

After finishing this paper I received a manuscript \cite{ST} by 
Smith and Thompson, which treats a similar problem from a different perspective.

\section{Preliminaries on complete ideals}

To begin with, we will briefly review some basic facts from the Zariski-Lipman theory of complete ideals. For more 
details, we refer to 
\cite{C}, \cite{L1}, \cite{L2}, \cite{L3} and \cite{Z}. Let $\mathcal K$ be a field. A two-dimensional regular local 
ring with the fraction field $\mathcal K$ is called a \textit{point}. 
The maximal ideal of a point $\ga$ is denoted by $\fm_\ga$. Write $\ord_\ga$ for the unique valuation of $\mathcal K$ 
such that for every $x \in \ga\smallsetminus\{0\}$
\[
\ord_\ga(x)=\max\{\nu\mid x\in\fm_\ga^\nu\}.
\]
A point $\gb$ is \textit{infinitely near} to a point $\ga$, if 
\[
\gb\supset\ga.
\]
Then the residue field extension $\ga/\fm_\ga\subset\gb/\fm_\gb$ is finite. In the following we will always consider 
points infinitely near to a fixed point, which has an algebraically closed residue field $\Bbbk$.

Take an element $x \in \fm_\ga\smallsetminus \fm_\ga^2$. A \textit{quadratic transform} of the point $\ga$ is a 
localization of the ring $\ga[\fm_\ga/x]$ at a maximal ideal of $\ga[\fm_\ga/x]$. Any two points $\ga\subset\gb$ can be 
connected by a unique sequence of quadratic transforms $\ga=\ga_1\subset\cdots\subset\ga_n=\gb$.

Assume that a point $\gb$ is a quadratic transform of the point $\ga$. Then the ideal $\fm_\ga\gb$ is generated by a 
single element, say $b\in\gb$. Let $\fa$ be an ideal in $\ga$. The \textit{transform} of an $\fa$ at $\gb$ is 
$\fa^{\gb}:=b^{-\ord_\ga(\fa)}\fa\gb$. For any point $\gb\supset\ga$, the transform $\fa^{\gb}$ is then defined 
inductively using the sequence of quadratic transforms connecting $\ga$ and $\gb$. It follows that if 
$\gamma \supset\gb\supset\ga$ are any points, then $\fa^\gamma=(\fa^\gb)^\gamma$. If $\fa$ has a finite colength, then 
so does $\fa^{\gb}$. Moreover, if $\fa$ is complete, i.e., integrally closed, then so is $\fa^\gb$.

The non-negative integer $\ord_\gb(\fa^\gb)$ is called the \textit{multiplicity} of $\fa$ \textit{at} $\gb$. In the 
case this is strictly positive we say that the point $\gb\supset\ga$ is a \textit{base point} of the ideal $\fa$.
The \textit{support} of $\fa$ is the set of the base points, which is known to be a finite set. The 
\textit{point basis}\label{pb} of the ideal $\fa$ is the family of multiplicities 
\[
\bB(\fa):=\{\ord_\gb(\fa^\gb)\mid\gb\supset\ga\}.
\]

Since the transform preserves products, i.e., $(\fa\fb)^{\gb}=\fa^{\gb}\fb^{\gb}$ for any finite colength ideals 
$\fa,\fb\subset\ga$, we have $\bB(\fa\fb)=\bB(\fa)+\bB(\fb)$. If in addition $\fa$ and $\fb$ are complete, then the 
condition $\ord_\gb(\fa^\gb)\ge \ord_\gb(\fb^\gb)$ for every point $\gb\supset\ga$ implies $\fa\subset\fb$. Moreover, 
$\fa=\fb$ exactly when $\ord_\gb(\fa^\gb)=\ord_\gb(\fb^\gb)$ for every point $\gb\supset\ga$.

By the famous result of Zariski \cite[p. 385]{Z} a complete ideal in $\ga$ factorizes uniquely into a product of simple 
complete ideals. Recall that an ideal is simple, if it is not a product of two proper ideals. There is a one to one 
correspondence between the simple complete ideals of finite colength in $\ga$ and the points containing $\ga$ (see e.g. 
\cite[p. 226]{L2}). The simple ideal corresponding to a point $\gb\supset\ga$ is then the unique ideal $\fa$ in $\ga$, 
whose transform at $\gb$ is the maximal ideal of $\gb$. The base points of $\fa$ are totally ordered by inclusion, and 
if $\ga=\ga_1\subset\dots\subset\ga_n$ is the sequence of the base points of $\fa$, then $\ga_n=\gb$. 
The point basis of $\fa$ is the vector 
$I:=(a_1,\dots,a_n)\in\mathbb N^n$, where $a_i:=\ord_{\ga_i}(\fa^{\ga_i})$.
Moreover, there is a one to one correspondence between points containing $\ga$ and divisorial valuations of $\ga$ given 
by $\gb\mapsto v:=\ord_\gb$ (\cite[pp. 389--391]{Z}). Note also that the ideals 
$\fm_\ga=\fp_1\supsetneq\cdots\supsetneq\fp_n=\fa$ corresponding to the points $\ga_1\subset\dots\subset\ga_n$ are the 
simple $v$-ideals containing $\fa$ (\cite[pp. 392]{Z}).

A point $\gb$ is said to be \textit{proximate} to a point $\ga$, if $\gb\supsetneq\ga$ and the valuation ring of 
$\ord_\ga$ contains $\gb$, in which case we write $\gb\succ\ga$. The notion of proximity can be interpreted 
geometrically as follows. The sequence of quadratic transforms $\ga=\ga_1\subset\dots\subset\ga_n$ corresponds to a 
sequence of regular surfaces
\begin{equation}\label{res}
\pi:\mathcal X=\mathcal X_{n+1}\xrightarrow{\pi_n}\cdots\xrightarrow{\pi_2}\mathcal X_2\xrightarrow{\pi_1}
\mathcal X_1=\Spec\ga,
\end{equation}
where $\pi_i:\mathcal X_{i+1}\rightarrow\mathcal X_i$ is the blow up of $\mathcal X_i$ at a closed point 
$\varsigma_i\in\mathcal X_i$ with $\cO_{\mathcal X_i,\varsigma_i}=\ga_i$ for every $i=1,\dots,n$, and 
$\pi=\pi_n\circ\cdots\circ\pi_1$. Set 
$E_i'=\pi_i^{-1}\{\varsigma_i\}\subset\mathcal X_{i+1}$. 
Write $E_i^*=(\pi_n\circ\cdots\circ \pi_{i+1})^*E_i'\subset\mathcal X$ for the total transform of $E_i'$ and let 
$E_i^{(j)}$ denote the strict transform of $E_i'$ on $\mathcal X_{j}$. Especially, set $E_i:=E_i^{(n+1)}$. 
If $j>i$, then $\varsigma_j\in E_i^{(j)}$ if and only if $\ga_j\succ\ga_i$, which is occasionally abbreviated by the 
notation $j\succ i$.

Recall that $\ga_{i}$ is always proximate to $\ga_{i-1}$, and there is at most one point $\ga_j$ such that 
$\ga_i\succ\ga_j$ and $j\neq i-1$. In the case such a point $\ga_j$ exists, $\ga_i$ is said to be a \textit{satellite} 
to $\ga_j$ (or shorter $i$ is satellite to $j$), and then also $\ga_\nu\succ\ga_j$ for every $j<\nu\le i$. In the 
opposite case $\ga_i$ (or just $i$) is said to be \textit{free}. Note that $\ga_2$ is always free. Also $\ga_1$ is 
regarded as a free point.

The proximity relations between the base points $\ga=\ga_1\subset\dots\subset\ga_n$ of a simple complete ideal $\fa$ 
of finite colength can be represented in the \textit{proximity matrix} of $\fa$. Following 
\cite[Definition-Lemma 1.5]{C}, it is 
\begin{equation}\label{prox}
P:=(p_{i,j})_{n\times n},\text{ where }p_{i,j}=
\left\lbrace
\begin{array}{rl}
   1,&\text{if }i=j;\\
  -1,&\text{if }i\succ j;\\
   0,&\text{otherwise.}
\end{array}
\right.
\end{equation}
Note that this is the transpose of the proximity matrix Lipman gives in \cite[p. 6]{L3}.  
Let $P^{-1}=(x_{i,j})_{n\times n}$ denote the inverse of the proximity matrix $P$, and write $X_i$ for the $i$:th row 
of $P^{-1}$. Observe that by e.g. \cite[Corollary 3.1]{L3} these rows correspond to the point basis vectors 
of the simple complete $\fm_\ga$-primary $v$-ideals containing $\fa$, and that they are in the descending order, so that 
$X_1$ represents $\fm_\ga$ and $X_n$ is the point basis of $\fa$. The next well known proposition shows that the 
proximity matrix of a simple complete ideal is totally determined by the point basis of the ideal. Because it will be 
crucial for our arguments in the sequel, we shall give a proof for the convenience of the reader.

\begin{prop}\label{1}
Let $\fa$ be a simple complete ideal of finite colength in a 
two-dimensional regular 
local ring $\ga$, and let $I=(a_1,\dots,a_n)$ be the point basis of $\fa$. For every index $i<k\le n$, we have 
\[
a_i=a_{i+1}+\cdots+a_k,
\] 
if and only if $\{\nu\mid\ga_\nu\succ\ga_i\}=\{i+1,\dots,k\}$. Then necessarily 
\[
a_i\ge a_{i+1}=\cdots=a_{k-1}\ge a_k,
\] 
and if $k>i+1$, then $a_{k-1}=a_k$ exactly when $k+1$ is free or $k=n$. Moreover, $a_{n-1}=a_n=1$. 
\end{prop}

\begin{proof}
As noted above, the point basis vector of $\fa$ is the bottom row of the inverse of the proximity matrix 
$P=(p_{i,j})_{n\times n}$ of $\fa$. Since $P^{-1}P=1$, we have 
\begin{equation}\label{ap} 
\sum_{\nu = 1}^n a_\nu p_{\nu,i}=\delta_{n,i}. 
\end{equation}  
for every $i=1,\dots,n$. 
Suppose that $i<n$. If $\{\nu\mid\ga_\nu\succ\ga_i\}=\{i+1,\dots,k\}$, then $p_{i+1,i}=\cdots=p_{k,i}=-1$ and 
$p_{j,i}=0$ for $j\notin\{i,\dots,k\}$, while $p_{i,i}=1$. It follows from 
Equation (\ref{ap}) that $a_i=a_{i+1}+\dots+a_k$.

Conversely, if $a_i=a_{i+1}+\dots+a_k$ and $k':=\max\{\nu\mid\ga_\nu\succ\ga_i\}$, then  
$p_{\nu,i}=-1$ if and only if $\nu\in\{i+1,\dots k'\}$, but then by Equation (\ref{ap}) we obtain 
\[
a_{i+1}+\cdots+a_k=a_i=a_{i+1}+\cdots+a_k'.
\] 
This forces the equality $k=k'$. Especially, choosing $i=n$ we recover the fact that $a_n=1$, and choosing $i=n-1$ 
yields $a_{n-1}=a_n$.

Let $i<j<k-1$. If $\ga_\nu\succ\ga_j$ for some $\nu>j+1$, then $\ga_{j+2}\succ\ga_j$. This is impossible, because 
$\ga_{j+2}$ is already proximate to $\ga_{j+1}$ and $\ga_i$. Therefore $\ga_{j+1}$ is the only point proximate to 
$\ga_j$, which implies that $a_{j}=a_{j+1}$. Now the second assertion is clear.

Suppose that $k-1>i$. We already observed that $a_{n-1}=a_n$, so we may assume $k<n$. If $\ga_{k+1}\succ \ga_j$ for 
some $j<k$, then also $\ga_k\succ \ga_j$. Because $\ga_k$ is already proximate to both $\ga_{k-1}$ and $\ga_i$, we must 
have $j\in\{i,k-1\}$. Since $k=\max\{\nu\mid\ga_\nu\succ\ga_i\}$ and $\ga_{k+1}\succ\ga_j$, the only possibility is 
$j=k-1$. Then both $\ga_k\succ\ga_{k-1}$ and $\ga_{k+1}\succ\ga_{k-1}$, which implies $a_{k-1}\ge a_k+a_{k+1}$. 
Especially, $a_{k-1}>a_k$. On the other hand, as $a_{k-1}=a_k+a_{k+1}+\cdots+a_{k'}$ for some $k'\ge k$, we see that 
$a_{k-1}>a_k$ implies $\ga_{k+1}\succ\ga_{k-1}$.
\end{proof}

The lattice $\Lambda:=\mathbb ZE_1+\cdots+\mathbb ZE_n$ of the exceptional divisors on $\mathcal X$ has two other 
convenient bases besides $\{E_i\mid i=1,\dots,n\}$, namely  $\{E_i^*\mid i=1,\dots,n\}$ and 
$\{\hat E_i\mid i=1,\dots,n\}$, where the $\hat E_i$:s are such that the intersection product with $E_j$ is the 
negative Kronecker delta, i.e., 
\[
E_i\cdot\hat E_j=-\delta_{i,j}
\] 
for $i,j=1,\dots,n$. Write 
\begin{equation}\label{E}
E:=(E_1,\dots,E_n)^\textsc{t}, E^*:=(E^*_1,\dots,E^*_n)^\textsc{t}\text{ and } 
\hat E:=(\hat E_1,\dots,\hat E_n)^\textsc{t},
\end{equation}
where \textsc{t} stands for the transpose. The proximity matrix $P$ is the base change matrix between $E$ and $E^*$, more 
precisely $E=P^\textsc{t}E^*$. 
Furthermore, $E^*_i\cdot E^*_j=-\delta_{i,j}$ (cf. \cite[p. 174]{C}). Then 
\begin{equation}\label{PP}
(E_i\cdot E_j)_{n\times n}=EE^\textsc{t}=P^\textsc{t}E^*(P^\textsc{t}E^*)^\textsc{t}=P^\textsc{t}(E_i^*\cdot E_j^*)_{n\times n}P
=-P^\textsc{t}P. 
\end{equation}
Assume that $E_i=\lambda_1\hat E_1+\cdots+\lambda_n \hat E_n$. Then $E_i\cdot E_j=-\lambda_j$, since 
$E_i\cdot\hat E_j=-\delta_{i,j}$. Therefore $E=P^\textsc{t}P\hat E$. This further implies that $\hat E=P^{-1}E^*$. 
Let $D\in\Lambda$ and let $d,d^*$ and $\hat d$ denote the row vectors in $\mathbb Z^n$ satisfying 
$D=dE=d^*E^*=\hat d\hat E$. Then \begin{equation}\label{d*}
d^*=dP^\textsc{t}\text{ and }\hat d=dP^\textsc{t}P=d^*P. 
\end{equation}
Note that the intersection product of 
$D$ and $F=fE\in\Lambda$ is  
\begin{equation}\label{d}
D\cdot F=dEE^\textsc{t}f^\textsc{t}=-dP^\textsc{t}Pf^\textsc{t}=-d^*(f^*)^\textsc{t}=-(d_1^*f_1^*+\cdots+d_n^*f_n^*).
\end{equation}

Recall that a divisor $D\in\Lambda$ is \textit{antinef}, if the intersection product $E_i\cdot D$ is non-positive for 
every $i$. According to Equation (\ref{d}) this means that $\hat d=dP^\textsc{t}P$ is non-negative at every entry, 
which can be abbreviated $\hat d=d^*P\ge0$. This is to say that the row vector $d^*=(d^*_1,\dots,d^*_n)$ satisfies the 
\textit{proximity inequalities} 
\begin{equation}\label{pxe}
d^*_i\ge\sum_{j\succ i}d^*_j.
\end{equation}
By \cite[Theorem 2.1]{L3} the proximity inequalities guarantee that there exist a unique complete 
ideal $\fd$ of finite colength in $\ga$ having 
the point basis $\bB(\fd)=d^*$. Then $\fd\cO_{\mathcal X}=\cO_{\mathcal X}(-D)$ so that 
$\fd=\Gamma(\mathcal X,\cO_{\mathcal X}(-D))$, as $\fd$ is complete. Thus we recover the following proposition 
(cf. \cite[\S18, p. 238-239]{L}):

\begin{prop}\label{nef}
There is a one to one correspondence between the antinef\\ divisors in $\Lambda$ and the complete ideals of finite 
colength in $\ga$ generating invertible $\cO_{\mathcal X}$-sheaves, given by 
$D\leftrightarrow\Gamma(\mathcal X,\cO_{\mathcal X}(-D))$. 
\end{prop}

\begin{rem}\label{VF}
Let $\fd$ be a complete ideal of finite colength in $\ga$ such that $\fd\cO_{\mathcal X}$ is invertible, and let 
$D=dE\in\Lambda$ be the antinef divisor corresponding to $\fd$ so that $\fd\cO_{\mathcal X}=\cO_{\mathcal X}(-D)$. The 
vector $\mathbf V(\fd):=d$ is called the \textit{valuation vector} of $\fd$. Observe that $d_i=\ord_{\ga_i}(\fd)$. 
Recall also that if $\fp_1,\dots,\fp_n$ are the simple $v$-ideals containing $\fa$, then 
$\fp_i\cO_{\mathcal X}=\cO_{\mathcal X}(-\hat E_i)$ so that $\fp_i=\Gamma(\mathcal X,\cO_{\mathcal X}(-\hat E_i)$. 
Because $D=\hat d\hat E$, we see that $\fd=\fp_1^{\hat d_1}\cdots\fp_n^{\hat d_n}$. We say that 
$\mathbf F(\fd):=\hat d$ is the \textit{factorization vector} of $\fd$. 
\end{rem}

For any divisor $D\in\Lambda$, there exists by e.g. \cite[Lemma 1.2]{LW} a minimal one among the antinef divisors 
$D^\sim$ satisfying $D^\sim\ge D$, which is to say that $D^\sim-D$ is effective. This is called the 
\textit{antinef closure} of $D$. 
According to \cite[\S18, p. 238]{L} an antinef divisor is effective. We can construct $D^\sim$ by the so called 
\textit{Laufer-algorithm} \label{LA} described in \cite[Proposition 1]{R}: Set $D_1=D$. For $i\ge1$, let $D_i=D^\sim$ 
when $D_i$ is antinef. Otherwise there exsist $\nu_i\in\{1,\dots,n\}$ such that $D_i\cdot E_{\nu_i}>0$. In this case 
set $D_{i+1}=D_i+E_{\nu_i}$. We have 
\begin{equation}\label{GL}
\Gamma(\mathcal X,\cO_{\mathcal X}(-D))=\Gamma(\mathcal X,\cO_{\mathcal X}(-D^\sim))
\end{equation} 
for any divisor $D\in\Lambda$ by \cite[Lemma 1.2]{LW}.

\section{Arithmetic of the point basis}

In this section we want to concentrate on the structure of the point basis of a simple ideal. From now on, let 
$\ga=\ga_1\subset\cdots\subset\ga_n$ be the quadratic sequence of the base points of a simple complete ideal $\fa$ of 
finite colength in a two-dimensional regular local ring $\ga$, and let $I=(a_1,\dots,a_n)$ denote the point basis of 
$\fa$.

\begin{defn}\label{rm11}
If $\ga_\gamma$ is a 
satellite point and $\ga_{\gamma+1}$ is not, then $\ga_\gamma$ is a \textit{terminal satellite} point. If $\ga_\tau$ is 
free and $\ga_{\tau+1}$ is not, then we say $\ga_\tau$ is a \textit{terminal free} point, and if $\ga_{\upsilon}$ is a 
free point but $\ga_{\upsilon-1}$ is not, then $\ga_{\upsilon}$ is a \textit{leading free} point. 
\end{defn}

\begin{rem}
Note that  $\ga_n$ is either a terminal satellite or a terminal free point. The first point $\ga_1$ is considered as a 
leading free point, and if $\ga_\gamma$ is a satellite point, then $2<\gamma\le n$. We say that the quadratic sequence 
$\ga=\ga_1\subset\cdots\subset\ga_n$ (or $\fa$) has a 
free point (a leading free point, a terminal free point, a satellite, a terminal satellite) at $i$, if $\ga_i$ is a 
free point (a leading free point, a terminal free point, a satellite, a terminal satellite, respectively). We may also 
say that $i$ is free (a satellite, a terminal satellite, resp.), if there is no confusion. 
\end{rem}

Clearly, every free point except $\ga_1$ belongs to a sequence of free points preceded by either $\ga_1$ or a terminal 
satellite. Moreover, except $\ga_1$, the leading free points are exactly the points immediately following a terminal 
satellite point. We use the following notation for terminal satellites and terminal free points.

\begin{ntt}\label{gamma}
Let us write $\Gamma_\fa$, or just $\Gamma$ if there is no confusion, for the set
\[
\{\gamma_1,\dots,\gamma_g\}:=\{\gamma\mid\ga_\gamma\text{ is a terminal satellite point of }\fa\}, 
\] 
where $\gamma_1<\cdots<\gamma_g$. Write also
\[
\Gamma^*:=\{\gamma\in\Gamma\mid\gamma<n\}=\{\gamma_1,\dots,\gamma_{g^*}\}.
\]
Moreover, let $\ga_{\tau_1},\dots,\ga_{\tau_g}$ denote the terminal free points of $\fa$ anterior to 
$\ga_{\gamma_g}$, so that 
\begin{equation}\label{tgt}
\tau_1<\gamma_1<\tau_2<\gamma_2<\cdots<\tau_g<\gamma_g.
\end{equation}
Set $\gamma_0:=1=:\tau_0$ and $\gamma_{g+1}:=n=:\tau_{g+1}$, and then define
\[
\bar\Gamma:=\Gamma\cup\{\gamma_{g+1}\}. 
\] 
\end{ntt}

Note that the number $g^*$ of elements in $\Gamma^*$ is either $g-1$ or $g$ depending on whether $\ga_n$ is a satellite 
or not. Thus $\gamma_{g^*+1}=n$ and $\bar\Gamma=\Gamma^*\cup\{n\}$. If necessary, we may write 
$\gamma_\nu=\gamma_\nu^{\fa}$ or $\tau_\nu=\tau_\nu^{\fa}$ to specify the ideal in question.

\begin{exmp}
If $\Gamma^*=\emptyset$, then the ideal is monomial, i.e., there exist a regular system of parameters for $\ga$ such 
that $\fa$ is monomial in these parameters.
\end{exmp}

\begin{prop}\label{taudit}
The multiplicity at a non-terminal free point is equal to the multiplicity at the terminal satellite point preceding 
it or to the multiplicity at the first point, if there is no preceding satellite point. The multiplicity at the 
subsequent terminal free point is strictly less, or the terminal free point is the last point, in which case they both 
are equal to one. Moreover, the multiplicities at a terminal satellite point and at the point immediately preceding it 
are equal.
\end{prop}

\begin{proof}
It follows from Proposition \ref{1} that $a_i\ge a_j$, if $j>i$. Moreover, $a_{i-1}>a_i$ if and only if $i+1$ is a 
satellite to $i-1$. By Definition \ref{rm11} we see that $i+1$ is free for $\gamma_{\nu-1}<i<\tau_\nu$. If 
$1\le\nu\le g$, then $\tau_\nu+1$ must be a satellite to $\tau_\nu-1$ as $\tau_\nu$ is free. Note also that 
$\gamma_\nu+1$ is not a satellite. Therefore 
\[
a_{\gamma_{\nu-1}}=\cdots=a_{\tau_\nu-1}>a_{\tau_\nu}\ge a_{\gamma_\nu-1}=a_{\gamma_\nu}
\] 
for $\nu=1,\dots,g$ .
Furthermore, because $\gamma_g$ is the last satellite, we must have $a_{\gamma_g}=\cdots=a_n$ and $a_n=1$, as observed 
in Proposition \ref{1}.
\end{proof}

For $\nu=1,\dots,g+1$, the sequence $a_{\gamma_{\nu-1}},\dots,a_{\gamma_\nu}$ of multiplicities may be written as 
\begin{equation}\label{sv}
\overbrace{s_{\nu,1},\dots,s_{\nu,1}}^{r_{\nu,1}\text{ times }},
\overbrace{s_{\nu,2},\dots,s_{\nu,2}}^{r_{\nu,2}\text{ times }},\dots,
\overbrace{s_{\nu,m_\nu},\dots,s_{\nu,m_\nu}}^{r_{\nu,m_\nu}+1\text{ times }}, 
\end{equation}
where $s_{\nu,1}>\dots>s_{\nu,m_\nu}$. Proposition \ref{taudit} implies that $m_\nu>1$ 
and $r_{\nu,m_\nu}\ge1$ for $1\le\nu\le g$, while $m_{g+1}=1$ and $r_{g+1,1}:=\gamma_{g+1}-\gamma_g\ge0$.
Since there are no terminal satellites 
between $\gamma_{\nu-1}$ and $\gamma_\nu$, Proposition \ref{1} yields
\begin{equation}\label{sr}
s_{\nu,j-1}=r_{\nu,j}s_{\nu,j}+s_{\nu,j+1} 
\end{equation}
for $1<j\le m_\nu$, where we set $s_{\nu,m_\nu+1}:=a_{\gamma_\nu}$ for $1\le\nu\le g+1$. Because of 
Proposition \ref{taudit}, the same holds also for $j=1$, if we set $s_{\nu,0}:=a_{\gamma_{\nu-1}}+\cdots+a_{\tau_\nu}$ 
for every $\nu=1,\dots,g+1$. When $j=m_\nu$, Equation (\ref{sr}) yields $s_{\nu,m_\nu}\mid s_{\nu,m_\nu-1}$. This shows that we 
may obtain the sequence $a_{\gamma_{\nu-1}},\dots,a_{\gamma_\nu}$ from 
$a_{\gamma_{\nu-1}}$ and $a_{\gamma_{\nu-1}}+\cdots+a_{\tau_\nu}$ by the Euclidian division algorithm.

\begin{ntt}\label{dpuis}
Let $\nu \in \{1,\dots,g+1\}$. We write $\gb'_\nu(\fa)$, or just $\gb'_\nu$ if 
the ideal is clear from the context, for the positive rational number
\[
\gb'_\nu:=\frac{a_{\gamma_{\nu-1}}+\cdots+a_{\tau_\nu}}{a_{\gamma_{\nu-1}}}. 
\]
We also set $\gb'_\nu:=0$ for $\nu>g+1$. 
\end{ntt}

\begin{prop}\label{aiaj1}
The point basis of $\fa$ is totally determined by the rational numbers 
$\gb'_1,\dots,\gb'_{g^*+1}$. In particular, the numbers $\gb'_1,\dots,\gb'_g$ yield the multiplicities 
$a_1,\dots,a_{\gamma_g}$. Moreover, for $\nu\in\{1,\dots,g+1\}$
\[
 \gcd\{a_1,\dots,a_{\gamma_\nu}\}
=\gcd\{a_{\gamma_{\nu-1}},a_{\tau_\nu}\}=a_{\gamma_\nu}.
\] 
\end{prop}

\begin{proof}
Using the Euclidian division algorithm described above we get  
\begin{equation}\label{euc}
\gcd\{a_{\gamma_{\nu-1}},a_{\tau_\nu}\}=\gcd\{s_{\nu,j},s_{\nu,j+1}\}=a_{\gamma_\nu}
\end{equation} 
for every $0\le j\le m_\nu$. Especially, $a_{\gamma_\nu}$ divides $a_{\gamma_{\nu-1}}$, and so $a_{\gamma_\nu}$ divides 
$a_i$ for every $i\le\gamma_\nu$. Hence $\gcd\{a_1,\dots,a_{\gamma_\nu}\}=a_{\gamma_\nu}$. This shows the last claim.

To prove the first two claims, we observe that for every $\nu=1,\dots,g+1$
\[
\gb'_\nu=\frac{s_{\nu,0}}{s_{\nu,1}}=\frac{\textsc n_\nu}{\textsc d_\nu}
\]
where $\textsc n_\nu$ and $\textsc d_\nu$ are integers with $\gcd\{\textsc n_\nu,\textsc d_\nu\}=1$. Suppose that we
know the pair $(\gb'_\nu,a_{\gamma_\nu})$. By Equation (\ref{euc}) we have  
$s_{\nu,0}=a_{\gamma_\nu}\textsc n_\nu$ and $s_{\nu,1}=a_{\gamma_{\nu-1}}=a_{\gamma_\nu}\textsc d_\nu$. Then we obtain 
the multiplicities $a_{\gamma_{\nu-1}},\dots,a_{\gamma_\nu}$ by the Euclidian division algorithm. 
Recall that $a_{\gamma_g}=a_n=1$ by Proposition \ref{taudit}. 
Starting from the pair $(\gb'_{g^*+1},a_n)$, or $(\gb'_g,a_{\gamma_g})$, 
we then get all the multiplicities
\[
a_1,\dots,a_n, \text{ or }a_1,\dots,a_{\gamma_g},
\]
respectively. 
\end{proof}

\begin{rem}\label{puis}
It follows from the Euclidian division algorithm (see Formula (\ref{sv}), Equation (\ref{sr}) and Proposition 
\ref{aiaj1}), that each $\gb'_\nu$ can be obtained from the integers $r_{\nu,1},\dots,r_{\nu,m_\nu}$ as a continued 
fraction
\[
\gb'_\nu=r_{\nu,1}+\frac{1}{r_{\nu,2}+\cdots+\frac{1}{r_{\nu,m_\nu}+1}}.
\] 
Note that these numbers are the \textit{Puiseux exponents} Spivakovsky defines in \cite[Definition 6.4]{S1}. 
\end{rem}

We now want to investigate the relationship between the point bases of the ideals $\fp_i$ to that of the ideal $\fa$. 
Let $P$ be the proximity matrix of $\fa$, and write $X_i$ for the $i$:th row of $(x_{i,j})_{n\times n}=P^{-1}$. 
Recall that $(x_{i,1},\dots,x_{i,i})$ is the point basis of $\fp_i$ and that $X_n=I$.

We first observe the following.

\begin{prop}\label{tg}
Let $i,j\in\{1,\dots,n\}$, with $i\ge j$. Then $\fp_i$ has a satellite or a free point at $j$, 
if and only if $\fa$ has a satellite or a free point at $j$, respectively. 
Moreover, this point is terminal for $\fp_i$ if and only if it is terminal for $\fa$ or $j=i$. 
If $i>j$, then $x_{i,j}=x_{i,j+1}$ if and only if 
$a_j=a_{j+1}$ or $i=j+1$.
\end{prop}

\begin{proof} 
The first two claims are obvious, since the base points of $\fp_i$ are $\ga_1\supset\cdots\supset\ga_i$, which are base 
points of $\fa$. Suppose then that $i>j$. By Proposition \ref{1} we have $x_{i,j}=x_{i,j+1}$ if and only if $\ga_{j+1}$ 
is the only point among the base points of $\fp_i$ proximate to $\ga_j$. An equivalent condition to this is that either 
$\ga_{j+2}$ is not proximate to $\ga_j$ or $i=j+1$, in other words, either $a_j=a_{j+1}$ or $i=j+1$.
\end{proof}

As a corollary to Proposition \ref{taudit} we now get

\begin{prop}\label{tau}
Let $\nu\in\{0,\dots,g\}$. Then $x_{i,j}=1$ for $\gamma_\nu\le j\le i\le\tau_{\nu+1}$.    
\end{prop}

\begin{proof}
Consider the ideal $\fp_i$. Because $\fa$ has a satellite point at $\gamma_\nu\le i$ and a free point at every $j$ 
satisfying $\gamma_\nu<j\le i$, Proposition \ref{tg} shows that $\gamma_\nu^{\fp_i}=\gamma_\nu$ and 
$\tau_{\nu+1}^{\fp_i}=i$. In particular, $\fp_i$ has the last satellite point at $\gamma_\nu$ and 
$\gamma_{\nu+1}^{\fp_i}=i$. It follows from Proposition \ref{taudit} that 
$x_{i,\gamma_\nu}=\cdots=x_{i,i}=1$.
\end{proof}

Note that the base points of the transform $\fp_m^{\ga_k}$ are $\ga_k\subset\cdots\subset\ga_m$ and the transform of 
$\fp_m^{\ga_k}$ at ${\ga_m}$ is the maximal ideal, i.e., $\fp_m^{\ga_k}$ is the simple complete ideal in $\ga_k$ 
corresponding to the point $\ga_m$. It follows that $(x_{m,k},\dots,x_{m,m})$ is the point basis of $\fp_m^{\ga_k}$. 
Observe also that the sequence of multiplicities $(a_k,\dots,a_m)$, where $1\le k\le m\le n$, is a point basis of a 
complete ideal, though not necessarily simple, because this sequence satisfies the proximity inequalities (\ref{pxe}). 
Indeed, we have  
\[
a_i-\sum_{m\ge j\succ i}a_{j}\ge a_i-\sum_{j\succ i}a_{j}\ge\gd_{i,n}\ge0
\]
for every $k\le i\le m$. This motivates the following notation.

\begin{ntt}
For $X=(x_1,\dots,x_n) \in \mathbb N^n$ and for $i,j\in\{1,\dots,n\}$ we write 
\begin{equation*}
X^{[i,j]}:=(\stackrel{(1)}{0},\dots,\stackrel{(i-1)}{0},x_i,\dots,x_j,\stackrel{(j+1)}{0},\dots,\stackrel{(n)}{0})
\end{equation*}
Moreover, we set 
\[
X^{[i,j)}:=X^{[i,j-1]}\text{ and }X^{(i,j]}:=X^{[i+1,j]}.
\]
In addition to that, we write 
\[
X^{\le i}:=X^{[1,i]}, X^{<i}:=X^{[1,i)}, X^{\ge i}:=X^{[i,n]}\text{ and }X^{>i}:=X^{(i,n]}.
\] 
\end{ntt}

For the truncated rows of $P^{-1}$ we obtain the following result.

\begin{prop}\label{tr}
Let $i,k\in\{1,\dots,n\}$. 
\[
X_i^{\le k}=
\left\lbrace
\begin{array}{ll}
X_i,&\text{ if $k\ge i$;}\\
x_{i,k}X_k,&\text{ if }k<j\text{ and }k+1\text{ is free;}\\ 
x_{i,k}X_{k}+\varrho_{i,k}X_h,&\text{ if }k<j\text{ and }k+1\text{ is a satellite to }h,
\end{array}
\right.
\]
where $\varrho_{i,k}:=x_{i,h}-(x_{i,h+1}+\cdots+x_{i,k})$. In particular,
\[
X_i^{\le\gamma_\nu}=x_{i,\gamma_\nu}X_{\gamma_\nu}, 
\]
where $\nu$ is such that $i\ge\gamma_\nu$.
\end{prop}

\begin{proof}
The case $k\ge i$ being trivial, we may restrict ourselves to the case $k<i$.
Proposition \ref{tg} now implies that we may replace $\fa$ by $\fp_i$, and so we may assume that $i=n$.
Consider the factorization vector $F:=X_n^{\le k}P$ so that 
\[
X_n^{\le k}=FP^{-1}=\sum_{j=1}^nF_jX_j.
\]
Now $F_j=0$ for $k<j\le n$, and if $1\le j\le k$, then 
\[
F_j=a_j-\sum_{k\ge\nu\succ j}a_\nu. 
\]
Because $a_j=\sum_{\nu\succ j}a_\nu+\gd_{j,n}$ by Proposition \ref{1}, we see that 
$F_j>0$, if and only if $k+1\succ j$ or $j=k=n$. Clearly, $F_k=a_k$.

It now follows that if $k+1$ is free or $k=n$, then $F_j=0$ for every $j\neq k$, implying $X_n^{\le k}=a_kX_k$. Assume 
then that $k+1$ is a satellite to some (unique) $h<k$. Now $F_h=a_h-(a_{h+1}+\cdots+a_k)$, while $F_j=0$ for every 
$j\notin\{h,k\}$. Then  
\[
X_n^{\le k}=a_kX_k+\left(a_h-\cdots-a_k\right)X_h, 
\]
as wanted.
\end{proof}

For certain computations, it is convenient to introduce the following notation.

\begin{ntt} 
Let $P^{-1}=(x_{i,j})_{n\times n}$ be the inverse of the proximity matrix of $\fa$.
For any $\nu\in\{0,\dots,g\}$ and $\gamma_\nu<i\le n$, write 
\[
\rho_{i,\gamma_\nu}=\rho_{i,\gamma_\nu}^{\fa}:=x_{i,\gamma_\nu+1}+\cdots+x_{i,\tau_{\nu+1}}\text{ and }\rho_\nu:=\rho_{\gamma_{\nu+1},\gamma_\nu}.
\]
In the case $i\le\gamma_\nu$ we set $\rho_{i,\gamma_\nu}=0$. 
\end{ntt}

\begin{rem}\label{rhob}
Observe that $\gb'_\nu=1+\rho_{n,\gamma_{\nu-1}}/a_{\gamma_{\nu-1}}$ for $\nu=1,\dots,g+1$. 
\end{rem}

\begin{cor}\label{ax}
Let $\gamma\in\{\gamma_0,\dots,\gamma_{g+1}\}$, and take an element $x_{i,j}$ of $P^{-1}$ with $j\le\gamma\le i$. 
Then $x_{i,j}=x_{i,\gamma}x_{\gamma,j}$. It follows that 
\begin{equation*}
a_{\gamma_\nu}=a_{\gamma_{\nu+1}}x_{\gamma_{\nu+1},\gamma_\nu}\text{ and }
a_{\gamma_\nu}=x_{\gamma_{\nu+1},\gamma_\nu}\cdots x_{\gamma_{g+1},\gamma_g}
\end{equation*} 
for $\nu=0,\dots,g$. Moreover, if $\gamma_{\nu+1}\le i$, then we have $\rho_{i,\gamma_\nu}=x_{i,\gamma_{\nu+1}}\rho_\nu$.
\end{cor}

\begin{proof}
For the first claim we observe that either $\gamma+1$ is free or $\gamma=n$. The latter case is trivial, while the former 
implies by Proposition \ref{tr} that $X_i^{\le\gamma}=x_{i,\gamma}X_{\gamma}$. Especially, this gives 
$x_{i,j}=x_{i,\gamma}x_{\gamma,j}$. Thus the first claim holds. Choosing 
$j=\gamma_\nu$, $\gamma=\gamma_{\nu+1}$ and $i=n$ we see that 
$a_{\gamma_\nu}=x_{\gamma_{\nu+1},\gamma_\nu}a_{\gamma_{\nu+1}}$, and then the second equality follows by induction. 
As $\tau_{\nu+1}\le\gamma_{\nu+1}$ we see that $x_{i,j}=x_{i,\gamma_\nu+1}x_{\gamma_\nu+1,j}$  
for every $j=\gamma_\nu+1,\dots,\tau_{\nu+1}$ and $i\ge\gamma_{\nu+1}$, which proves the last claim.
\end{proof}

\begin{prop}\label{rhow}
Let $\fp_i$ be the simple complete $v$-ideal corresponding to the row $X_i$ and 
let $\mu$ be the integer satisfying $\gamma_\mu<i\le\gamma_{\mu+1}$ unless $i=1$, in which case we set $\mu:=0$.
For every $\nu\in\{0,\dots,\mu+1\}$, we have 
\[
\rho_{i,\gamma_\nu^{\fp_i}}^{\fp_i}=\rho_{i,\gamma_\nu}.
\] 
\end{prop}

\begin{proof}
Suppose first that $0\le\nu<\mu$. Then $\gamma_\nu<i$. It follows from Proposition \ref{tg} that 
$\gamma_\nu^{\fp_i}=\gamma_\nu$. Moreover, $\tau_{\nu+1}^{\fp_i}=\min\{i,\tau_{\nu+1}\}$. 
As $x_{i,j}=0$ for $j>i$ we then get
\begin{equation}\label{tpi}
\rho_{i,\gamma_\nu^{\fp_i}}^{\fp_i}=x_{i,\gamma_\nu^{\fp_i}+1}+\cdots+x_{i,\tau_{\nu+1}^{\fp_i}}
=x_{i,\gamma_\nu+1}+\cdots+x_{i,\tau_{\nu+1}}=\rho_{i,\gamma_\nu}.
\end{equation}
If $\nu=\mu$, then $\gamma_\nu<i$, but it may happen that $\tau_\nu\ge i$. If this is the case, then $\fa$ has a free 
point at $i$. By Proposition \ref{tg} this means that $\tau_{\nu+1}^{\fp_i}=i$, and since $x_{i,j}=0$ for $j>i$, 
we observe that Equation (\ref{tpi}) holds. 
If $\nu=\mu+1$, then $\gamma_\nu\ge\gamma_\nu^{\fp_i}=i$ and by definition  
$\rho_{i,\gamma_\nu^{\fp_i}}^{\fp_i}=\rho_{i,\gamma_\nu}=0$. 
\end{proof}

\begin{cor}
if $\gamma_\nu\le i$, then $\gb'_\nu(\fp_i)$ doesn't depend on $i$.
\end{cor}

\begin{proof}
Assuming $\gamma_\nu\le i$ we get by using Corollary \ref{ax} 
\[
\frac{\rho_{i,\gamma_{\nu-1}}}{x_{i,\gamma_{\nu-1}}}
=\frac{x_{i,\gamma_\nu}\rho_{\nu-1}}{x_{i,\gamma_\nu}x_{\gamma_\nu,\gamma_{\nu-1}}}
=\frac{a_{\gamma_\nu}\rho_{\nu-1}}{a_{\gamma_\nu}x_{\gamma_\nu,\gamma_{\nu-1}}}
=\frac{\rho_{n,\gamma_{\nu-1}}}{a_{\gamma_{\nu-1}}}.
\]
Since $\gamma_{\nu-1}\le i$, it follows from Proposition \ref{rhow} that 
$\rho^{\fp_i}_{i,\gamma_{\nu-1}^{\fp_i}}=\rho_{i,\gamma_{\nu-1}}$. Because 
$\gb'_\nu(\fp_i)=1+\rho^{\fp_i}_{i,\gamma_{\nu-1}^{\fp_i}}/x_{i,\gamma_{\nu-1}}$ and
$\gb'_\nu=1+\rho_{n,\gamma_{\nu-1}}/a_{\gamma_{\nu-1}}$, we then observe that
$\gb'_\nu(\fp_i)=\gb'_\nu$. Thus we get the claim.
\end{proof}

\section{The Dual graph}

Let $\fa$ be a simple complete ideal of finite colength in a two-dimensional regular local ring $\ga$. Consider the 
resolution (\ref{res}) of $\fa$.
The configuration of the exceptional divisors $E_1,\dots,E_n$ on $\mathcal X$ arising from (\ref{res}) can 
be represented by a weighted graph (see, e.g., \cite[p. 111, §5. pp. 124 -- 129 ]{S}). This graph is called the 
\textit{dual graph} of $\fa$.  Its vertices are $\epsilon_1,\dots,\epsilon_n$, where $\epsilon_i$ corresponds to 
$E_i$ for every $i$. The \textit{weight} of a vertex $\epsilon_i$ is $w_i:=-E_i\cdot E_i$. Two vertices are 
\textit{adjacent}, if they are joined by an edge. This takes place, if and only if the intersection of the 
corresponding divisors on $\mathcal X$ is nonempty. A vertex is called an \textit{end}, if it is 
adjacent to at most one vertex. If it has more than two adjacent vertices, it is said to be a \textit{star}\label{sdg}.

\begin{prop}\label{adj}
Let $\epsilon_i$ and $\epsilon_j$ be vertices with $i<j$ in the dual graph of $\fa$. They are adjacent, if and only if 
$j=\max\{\nu\mid \nu\succ i\}$. In particular, if this is the case, then $j$ is uniquely determined.
\end{prop}

\begin{proof}
Vertices $\epsilon_i$ and $\epsilon_j$ are adjacent, if and only if $E_i\cdot E_j \neq 0$. By  
Equation (\ref{PP}) and by the definition of the proximity matrix we have 
\[
-E_i\cdot E_j=\sum_{\nu=j}^\mu p_{\nu,i}p_{\nu,j}, 
\]
where $\mu:=\max\{\nu\mid \nu\succ i\}$. Clearly, $j>\mu$ implies $E_i\cdot E_j=0$. If $j<\mu$, then 
\[
-E_i\cdot E_j=p_{j+2,i}p_{j+2,j}+\cdots+p_{\mu,i}p_{\mu,j}.
\] 
Observe that this is nonzero, if and only if 
$p_{j+2,i}p_{j+2,j}\neq 0$, but this is impossibile, because then $j+2$ would be proximate to $i,j$ and $j+1$. 
Then we see that $E_i\cdot E_j\neq0$ with $i<j$ if and only if $j=\mu$, in which case $E_i\cdot E_j=1$. 
\end{proof}

\begin{rem}\label{adj1}
Observe that $w_i=-E_i\cdot E_i=1+\mu-i$, where we write $\mu:=\max\{\nu\mid \nu\succ i\}$ so that $\mu-i$ is the number of points 
proximate to $i$. Thus vertices $\epsilon_i$ and $\epsilon_j$ with $i<j$ are adjacent, if and only if 
\[
j=i+w_i-1.
\]
Furthermore, if $i\neq j$, then $E_i\cdot E_j=\delta_{i,i+w_i-1}$. Thereby we observe that the matrix $P^{\textsc t}P$, 
which represents the dual graph of $\fa$ by Equation (\ref{PP}), is totally determined by the sequence of the weights 
$(w_1,\dots,w_n)$. Moreover, we obtain the proximity matrix $P=(p_{i,j})_{n\times n}$ from the $w_i$:s, since 
$p_{i,j}=-1$, if and only if $i<j<i+w_i$. Thus the sequence $(w_1\dots,w_n)$ gives equivalent data to the point basis 
of $\fa$.
\end{rem}

\begin{prop}\label{star}
The stars of the dual graph of $\fa$ are precisely the vertices $\epsilon_\gamma$ such that $\gamma<n$ and $\fa$ has a 
terminal satellite at $\gamma$. The ends of the dual graph are exactly the vertices $\epsilon_\tau$ such that $\tau=1$ 
or $\fa$ has a terminal free point at $\tau$.
\end{prop}

\begin{proof}
Suppose that $\gamma<n$ is a terminal satellite. By Proposition \ref{adj} we know that $\epsilon_\mu$ is adjacent to 
$\epsilon_\gamma$ for $\mu=\max\{i\mid i\succ\gamma\}$. Because $\gamma$ is a terminal satellite, we get for some 
$\nu<\gamma-1$
\[
\gamma=\max\{i\mid i\succ\gamma-1\}=\max\{i\mid i\succ \nu\}. 
\] 
By using Proposition \ref{adj} again, we see that $\epsilon_\nu$, $\epsilon_{\gamma-1}$ and $\epsilon_\mu$ are adjacent 
to $\epsilon_\gamma$. Thus $\epsilon_\gamma$ is a star.

Conversely, suppose that $\epsilon_\gamma$ is a star. Then there are three vertices adjacent to $\epsilon_\gamma$, say 
$\epsilon_i$, $\epsilon_j$ and $\epsilon_k$. Noting that a point cannot be proximate to three different points, it 
follows from Proposition \ref{adj} that exactly one of the indices is greater that $\gamma$. Therefore we may suppose 
that $j<k<\gamma<i$. In particular, $\gamma<n$. If $\epsilon_{\gamma+1}$ is a satellite to some $m$, then also 
$\gamma\succ m$, which implies that $m\in \{j,k\}$. On the other hand, as $\epsilon_j$ and $\epsilon_k$ are adjacent to 
$\epsilon_\gamma$ we see by Proposition \ref{adj} that $\gamma+1$ is proximate to neither $j$ nor $k$, which is a 
contradiction. Therefore $\gamma+1$ is free, and so $\gamma$ is a terminal satellite.

Suppose that $\fa$ has a terminal free point at $\tau>1$. Suppose also that $\epsilon_i$ is adjacent to $\epsilon_\tau$. 
If $i<\tau$, then $\tau\succ i$ by Proposition \ref{adj}, and because $\tau$ is not a satellite to any $i<\tau-1$, it 
follows that $i=\tau-1$, but since $\epsilon_i$ is adjacent to $\epsilon_\tau$,
Proposition \ref{adj} implies that $\tau+1\nsucc\tau-1$. On the other hand, we know that $\tau+1$ is not free as 
$\tau$ is terminal. This means that $\tau=n$, and $\epsilon_{\tau-1}$ is the only vertex adjacent to $\epsilon_\tau$, 
i.e., $\tau$ is an end. Furthermore, we see that $\tau<n$ implies $i>\tau$, but then $i$ is uniquely determined by 
Proposition \ref{adj}. Hence $\epsilon_\tau$ is an end. For the same reason $\epsilon_1$ is an end.

Let us then prove the converse.  As we just saw, $\epsilon_1$ is always an end. Proposition \ref{adj} yields that 
$\epsilon_j$ adjacent to $\epsilon_n$ whenever $n\succ j$. Thus $\epsilon_n$ is an end, if and only if $n$ is free. 
Suppose then that $\epsilon_i$ is an end with $1<i<n$. Then Proposition \ref{adj} implies that the only vertex adjacent 
to $\epsilon_i$ is $\epsilon_\mu$, where $\mu:=\max\{\nu\mid\nu\succ i\}$ is greater than $i$. Moreover, if $i\succ j$ 
for some $j<i$, then $i+1\succ j$, too. Especially, this shows that $i+1\succ i-1$, and since $i+1\succ i$, it follows 
that $i\nsucc j$ for any $j<i-1$. Therefore $i$ is free, while $i+1$ is a satellite. So, if $\epsilon_i$ in an end for 
some $i\in\{1,\dots,n\}$, then $i=1$ or $\fa$ has a terminal free point at $i$. 
\end{proof}

\begin{rem}
By Proposition \ref{star} we observe that the vertex $\epsilon_i$ is a star exactly, when $i\in\Gamma^*$, and $g^*$ is 
the number of the star vertices of the dual graph.
\end{rem}

We will now recall how the dual graph can be constructed from the point basis $(a_1,\dots,a_n)$ of $\fa$. Let 
$\{\gamma_0,\dots,\gamma_{g+1}\}$ be as given in Notation \ref{gamma}. Let also the integers $s_{\nu,\mu}$, 
$r_{\nu,\mu}$ and $m_\nu$ be as in Formula (\ref{sv}), where $\nu=1,\dots,g+1$ and $\mu=1,\dots,m_\nu$. Set  
\begin{equation}\label{kappa}
\kappa_{\nu,\mu}:=\sum_{i=1}^{\nu-1}\sum_{j=1}^{m_i}r_{i,j}+\sum_{j=1}^{\mu}r_{\nu,j}
\end{equation}
for $\nu=1,\dots,g+1$ and $\mu=0,\dots,m_\nu$. It follows that for $\mu>0$
\begin{equation*}
\kappa_{\nu,\mu}=\kappa_{\nu,\mu-1}+r_{\nu,\mu},
\text{ and }
\kappa_{\nu,\mu}=\kappa_{\nu-1,m_{\nu-1}}+\sum_{j=1}^{\mu}r_{\nu,j} 
\end{equation*}
when $\nu>1$. Moreover, for every $\nu=1,\dots,g+1$ we have
\[
\kappa_{\nu,0}=\gamma_{\nu-1}-1 \text{ and } \kappa_{\nu,m_\nu}=\gamma_\nu-1.
\]
Note that $a_i=s_{\nu,\mu}$ for $\kappa_{\nu,\mu-1}<i\le\kappa_{\nu,\mu}$, and $a_{\kappa_{\nu,m_\nu}}=a_{\gamma_\nu}$.

\begin{lem}\label{kah}
Let $\epsilon_1,\dots,\epsilon_n$ be the vertices of the dual graph of $\fa$. Assume that 
$\kappa_{\nu,\mu-1}<i<\kappa_{\nu,\mu}$ or $i=\kappa_{\nu,m_\nu}$ for $\nu\in\{1,\dots,g+1\}$ and 
$\mu\in\{1,\dots,m_\nu\}$. If $j>i$, then $\epsilon_j$ is adjacent to $\epsilon_i$, if and only if $j=i+1$. Moreover, 
$\epsilon_{\kappa_{\nu,\mu}}$ is adjacent to $\epsilon_{\kappa_{\nu,\mu+1}+1}$ for $\mu\in\{1,\dots,m_\nu-1\}$
\end{lem}

\begin{proof} 
If $\kappa_{\nu,\mu-1}<i<\kappa_{\nu,\mu}$ or $i=\kappa_{\nu,m_\nu}$, then we have $a_i=a_{i+1}$. It follows from Proposition 
\ref{1} that $j\succ i$ if and only if $j=i+1$, which is to say that $i$ has no satellites. 
By Proposition \ref{adj} this takes place exactly when $\epsilon_i$ is adjacent $\epsilon_{i+1}$.

Assume that $\mu\in\{1,\dots,m_\nu-1\}$. An application of Equation (\ref{sr}) then shows that and 
$i=\kappa_{\nu,\mu}$. As we observed above,
$a_{\kappa_{\nu,\mu}}=a_{\kappa_{\nu,\mu}+1}+\cdots+a_{\kappa_{\nu,\mu+1}+1}$. So 
$\kappa_{\nu,\mu+1}+1=\max\{j\mid j\succ\kappa_{\nu,\mu}\}$ by Proposition \ref{1}, and the claim follows from 
Proposition \ref{adj}.
\end{proof}

Using Lemma \ref{kah} together with Proposition \ref{star} and Remark \ref{adj1}, we are able to construct the 
dual graph of $\fa$ from the point basis. 
The figure below describes a fragment of the dual graph of $\fa$. It illustrates the organization of the vertices 
$\epsilon_{\gamma_{\nu-1}},\dots,\epsilon_{\gamma_\nu}$ and the corresponding multiplicities 
$a_{\gamma_{\nu-1}},\dots,a_{\gamma_\nu}$ in the point basis of $\fa$.

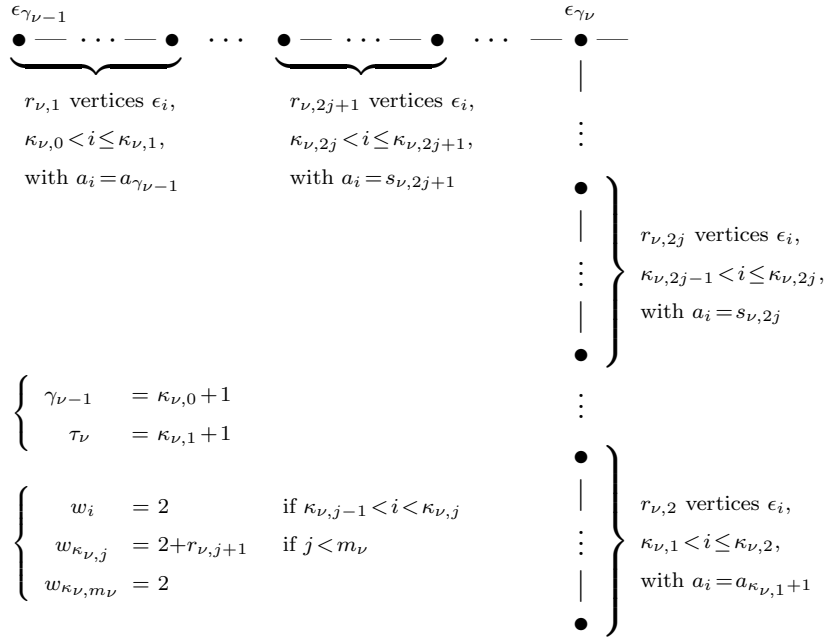
\[
\phantom{...}
\begin{array}{cclclccclc}
  &
  &\scriptstyle \epsilon_{\gamma_{\nu-1}}
  &&&
  &\scriptstyle \epsilon_{\gamma_\nu}
  &&&\vspace{-2pt}
\\
  &
  &\underbrace{\bullet\hspace{3pt}\text{|}\hspace{-1pt}\phantom{_\chi}{\cdots}\hspace{3pt}\text{|}\hspace{3pt}\bullet}
  &\cdots
  &\underbrace{\bullet\hspace{3pt}\text{|}\hspace{-1pt}\phantom{_\chi}{\cdots}\hspace{3pt}\text{|}\hspace{3pt}\bullet}
  &\cdots\hspace{6pt}\text{|}\hspace{-11pt}
  &\bullet&\hspace{-12,5pt}\text{|}&\vspace{-7pt}
\\
  &
  &&&&&\mid
  &&\vspace{-6pt}
\\
  &&\makebox[0\height][l]{$\begin{array}{l}
             \scriptstyle r_{\nu,1}\text{ vertices }\epsilon_i,\\ 
             \scriptstyle \kappa_{\nu,0}\hspace{1pt}<\hspace{1pt}i\hspace{1pt}\le\hspace{1pt} \kappa_{\nu,1},\\
             \scriptstyle \text{with }a_i \hspace{1pt}=\hspace{1pt} a_{\gamma_{\nu-1}}
             \end{array}$}
  &&\makebox[0\height][l]{$\begin{array}{l}
             \scriptstyle r_{\nu,2j+1}\text{ vertices }\epsilon_i,\\ 
             \scriptstyle \kappa_{\nu,2j}\hspace{1pt}<\hspace{1pt}i\hspace{1pt}\le\hspace{1pt}\kappa_{\nu,2j+1},\\
             \scriptstyle \text{with }a_i\hspace{1pt}=\hspace{1pt}s_{\nu,2j+1}
             \end{array}$}
  &&\vdots 
  &&\vspace{-10pt}
\\
  &&
  &&&&\begin{array}{c}\bullet\\\mid\\\vdots\\\mid\\\bullet\end{array}
  &\hspace{-25pt}\left.\begin{array}{c}\phantom.\\\\\\\\\phantom.\end{array}\right\rbrace
  &\hspace{-11pt}
    \begin{array}{l}
      \scriptstyle r_{\nu,2j}\text{ vertices }\epsilon_i,\\ 
      \scriptstyle \kappa_{\nu,2j-1}\hspace{1pt}<\hspace{1pt}i\hspace{1pt}\le\hspace{1pt}\kappa_{\nu,2j},\\
      \scriptstyle \text{with }a_i \hspace{1pt}=\hspace{1pt}s_{\nu,2j} 
    \end{array}
\\
 &&\vspace{0pt}\makebox[0\height][l]{$\scriptstyle\left\lbrace
   \begin{array}{cl}
     \scriptstyle\gamma_{\nu-1}\hspace{4pt}&\scriptstyle=\hspace{3pt}\kappa_{\nu,0}\hspace{1pt}+\hspace{1pt}1\\
     \scriptstyle\tau_\nu\hspace{-4pt}&\scriptstyle=\hspace{3pt}\kappa_{\nu,1}\hspace{1pt}+\hspace{1pt}1
   \end{array}
   \right.$}
 &&&
 &\vdots &&\vspace{-5pt}
\\
 &&\makebox[0pt][l]{\vspace{100pt}$\scriptstyle 
   \left\lbrace
   \begin{array}{cll}
     \scriptstyle w_i\hspace{-6pt}&\scriptstyle=\hspace{3pt}\scriptstyle2
                                 &\scriptstyle\text{ if }\kappa_{\nu,j-1}\hspace{1pt}<
                                  \hspace{1pt}i\hspace{1pt}<\hspace{1pt}\kappa_{\nu,j}\\
     \scriptstyle w_{\kappa_{\nu,j}}\hspace{-6pt}&\scriptstyle=\hspace{3pt}\scriptstyle2+r_{\nu,j+1}
                                 &\scriptstyle\text{ if }j\hspace{1pt}<\hspace{1pt}m_\nu\\
     \scriptstyle w_{\kappa_{\nu,m_\nu}}\hspace{-6pt}&\scriptstyle=\hspace{3pt}\scriptstyle2
   \end{array}
   \right.$}
 &&&&\begin{array}{c}\bullet\\\mid\\\vdots\\\mid\\\bullet\end{array}
 &\hspace{-25pt}\left.\begin{array}{c}\phantom.\\\\\\\\\phantom.\end{array}\right\rbrace
 &\hspace{-11pt}   
   \begin{array}{l}
     \scriptstyle r_{\nu,2}\text{ vertices }\epsilon_i, \\ 
     \scriptstyle \kappa_{\nu,1}\hspace{1pt}<\hspace{1pt}i\hspace{1pt}\le\hspace{1pt}\kappa_{\nu,2},\\
     \scriptstyle \text{with }a_i \hspace{1pt}=\hspace{1pt} a_{\kappa_{\nu,1}+1}
    \end{array}
\\
\\
\end{array}
\]

Observe that the vertices $\epsilon_i$ for $\kappa_{\nu,m_\nu-1}<i\le\kappa_{\nu,m_\nu}$ with the multiplicity 
$a_i=a_{\gamma_\nu}$ lie at the horizontal or vertical branch, depending on whether $m_\nu$ is odd or even, respectively. 
Note also that for $\nu=1,\dots,g$ the vertex $\epsilon_{\gamma_\nu}$ belongs to the next segment of the dual graph.

\begin{exmp}
Let $\fa$ be a simple complete ideal in a two-dimensional regular local ring $\ga$ having the resolution (\ref{res}). 
and let $n=8$ in this resolution. Suppose that $i\succ j$ for $i,j\in\{1,\dots,8\}$ with $i>j$, 
if and only if $i=j+1$ or $(i,j)\in\{(3,1),(6,4),(7,4)\}$. The proximity 
matrix of $\fa$ is then 
\[
\phantom{^{-1}P^{\textsc t}}P=
\left[
\begin{array}{cccccccc}
  1&\cdot&\cdot&\cdot&\cdot&\cdot&\cdot&\phantom{-}0\phantom{-}\\  
 -1\phantom{-}&1&\cdot&\cdot&\cdot&\cdot&\cdot&\phantom{-}\cdot\phantom{-}\\  
 -1\phantom{-}&-1\phantom{-}&1&\cdot&\cdot&\cdot&\cdot&\phantom{-}\cdot\phantom{-}\\  
  \cdot&\cdot&-1\phantom{-}&1&\cdot&\cdot&\cdot&\phantom{-}\cdot\phantom{-}\\ 
  \cdot&\cdot&\cdot&-1\phantom{-}&1&\cdot&\cdot&\phantom{-}\cdot\phantom{-}\\ 
  \cdot&\cdot&\cdot&-1\phantom{-}&-1\phantom{-}&1&\cdot&\phantom{-}\cdot\phantom{-}\\  
  \cdot&\cdot&\cdot&-1\phantom{-}&\cdot&-1\phantom{-}&1&\phantom{-}\cdot\phantom{-}\\  
  0&\cdot&\cdot&\cdot&\cdot&\cdot&-1\phantom{-}&\phantom{-}1\phantom{-}\\  
\end{array}
\right],
\]
and the inverse of $P$ is 
\[
\phantom{P^{\textsc t}}P^{-1}=
\left[
\begin{array}{cccccccccc}
 \phantom{-}1\phantom{-}&\cdot&\cdot&\cdot&\cdot&\cdot&\cdot&\phantom{-}0\phantom{-}\\  
            1&\phantom{-}1\phantom{-}&\cdot&\cdot&\cdot&\cdot&\cdot&\phantom{-}\cdot\phantom{-}\\  
            2&1&\phantom{-}1\phantom{-}&\cdot&\cdot&\cdot&\cdot&\phantom{-}\cdot\phantom{-}\\  
            2&1&1&\phantom{-}1\phantom{-}&\cdot&\cdot&\cdot&\phantom{-}\cdot\phantom{-}\\ 
            2&1&1&1&\phantom{-}1\phantom{-}&\cdot&\cdot&\phantom{-}\cdot\phantom{-}\\ 
            4&2&2&2&1&\phantom{-}1\phantom{-}&\cdot&\phantom{-}\cdot\phantom{-}\\  
            6&3&3&3&1&1&\phantom{-}1\phantom{-}&\phantom{-}\cdot\phantom{-}\\  
            6&3&3&3&1&1&1&\phantom{-}1\phantom{-}
\end{array}
\right].
\]
Thus the point basis of $\fa$ is $I=(a_1,\dots,a_8)=(6,3,3,3,1,1,1,1)$. The dual graph of $\fa$ is presented in the 
matrix
\[
\phantom{^{-1}}P^{\textsc t}P=
\left[
\begin{array}{cccccccc}
  3&\cdot&-1\phantom{-}&\cdot&\cdot&\cdot&\cdot&\phantom{-}0\phantom{-}\\  
  \cdot&2&-1\phantom{-}&\cdot&\cdot&\cdot&\cdot&\phantom{-}\cdot\phantom{-}\\  
 -1\phantom{-}&-1\phantom{-}&2&-1\phantom{-}&\cdot&\cdot&\cdot&\phantom{-}\cdot\phantom{-}\\  
  \cdot&\cdot&-1\phantom{-}&4&\cdot&\cdot&-1\phantom{-}&\phantom{-}\cdot\phantom{-}\\ 
  \cdot&\cdot&\cdot&\cdot&2&-1\phantom{-}&\cdot&\phantom{-}\cdot\phantom{-}\\ 
  \cdot&\cdot&\cdot&\cdot&-1\phantom{-}&2&-1\phantom{-}&\phantom{-}\cdot\phantom{-}\\  
  \cdot&\cdot&\cdot&-1\phantom{-}&\cdot&-1\phantom{-}&2&-1\phantom{-}\\  
  0&\cdot&\cdot&\cdot&\cdot&\cdot&-1\phantom{-}&\phantom{-}1\phantom{-}\\  
\end{array}
\right],
\]
We may draw the dual graph as follows: 
\[
\begin{array}{ccccccccc}
\stackrel{(\epsilon_1)}3&\text{|}&\stackrel{(\epsilon_3)}2&\text{|}&\stackrel{(\epsilon_4)}4 
&\text{|}&\stackrel{(\epsilon_7)}2&\text{|}&\stackrel{(\epsilon_8)}1\\
&&|&&&&|&&\\
&&\stackrel{(\epsilon_2)}2&&&&\stackrel{(\epsilon_6)}2&&\\
&&&&&&|&&\\
&&&&&&\stackrel{(\epsilon_5)}2&&
\end{array}
\]
The stars of the dual graph are $\epsilon_3$ and $\epsilon_7$, while the ends are $\epsilon_1,\epsilon_2,\epsilon_5$ 
and $\epsilon_8$. Thus $g=2$, $\Gamma=\{3,7\}$ and $\{\tau_0,\dots,\tau_3\}=\{1,2,5,8\}$. Now $\gb'_1=9/6$, 
$\gb'_2=7/3$ and $\gb'_3=2$. On the other hand, starting from $(\gb'_1,\dots,\gb'_{g+1}) = (3/2,7/3,2)$ we get first  
$(a_{\gamma_0},a_{\gamma_1},a_{\gamma_2})=(6,3,1)$ and then $(a_{\gamma_0}\gb'_1,\dots,a_{\gamma_2}\gb'_3)=(9,7,2)$. 
By using the Euclidian algorithm as above we get $(a_1,\dots,a_{\gamma_1})=(6,3,3)$, 
$(a_{\gamma_1},\dots,a_{\gamma_2})=(3,3,1,1,1)$ and $(a_{\gamma_2},\dots,a_n)=(1,1)$, and so we may reconstruct the 
point basis, which yields the proximity matrix by Proposition \ref{1}. Note that the dual graph can be obtained 
directly from the point basis by the following manner: From every entry $a_1,\dots,a_{n-1}$ in the point basis, draw an 
arc extending to $a_k$ such that $a_i=a_{i+1}+\cdots+a_{k}$, i.e.,
\[
 6\overparenthesis{\phantom{tällai}3\underparenthesis{\phantom{tällai}}_{(2)}}^{(3)}
 3\overparenthesis{\phantom{tällai}}^{(2)}3\overparenthesis{\phantom{tällai}
 1\underparenthesis{\phantom{tällai}}_{(2)}1\underparenthesis{\phantom{tällai}}_{(2)}}^{(4)}
 1\overparenthesis{\phantom{tällai}}^{(2)}1
\]
The entries and the arcs correspond the vertices and the edges of the dual graph. Moreover, the lengt $k-i+1$ of the 
arc starting from $a_i$, indicated in parenthesis, is exactly the weight $w_i$, while $w_n = 1$. 
\end{exmp}

\section{On certain intersection products}

Let $\fa$ be a simple complete ideal of finite colength in a regular local ring $\ga$ with the base points 
$\ga=\ga_1\subset\cdots\subset\ga_n$. Let $X_i$ denote the $i$:th row of the inverse $(x_{i,j})_{n\times n}:=P^{-1}$ of 
the proximity matrix. For any rows $X_i$ and $X_j$ of $P^{-1}=(x_{i,j})_{n\times n}$, write $X_i\cdot X_j$ for the dot 
product, i.e., 
\[
X_i\cdot X_j:=x_{i,1}x_{j,1}+\cdots+x_{i,n}x_{j,n}, 
\]
and for any $\nu\in\{0,\dots,g\}$, write  
\[
[X_i\cdot X_j]_\nu:=X_i^{(\gamma_\nu,\gamma_{\nu+1}]}\cdot X_j^{(\gamma_\nu,\gamma_{\nu+1}]}, 
\]
so that 
\[
X_i\cdot X_j=x_{i,1}x_{j,1}+[X_i\cdot X_j]_0+\cdots+[X_i\cdot X_j]_g. 
\] 
Assume that the integers $\kappa_{\nu,\mu}$ attached to 
the ideal $\fa$ are as in Equation (\ref{kappa}). Then define 
\begin{equation}\label{U}
U := \{i\mid\kappa_{\nu,\mu-1}<i\le\kappa_{\nu,\mu}\text{ for any }\nu\text{ and }\mu\notin2\mathbb N,\text{ or }i=n\}.
\end{equation}

\begin{rem}\label{keyU}
Note that since $\gamma_\nu=\kappa_{\nu+1,0}+1$ for $\nu\in\{0,\dots,g+1\}$, we have 
\[
\gamma_\nu\in U\text{ for }\nu\in\{0,\dots,g+1\}.
\] 
Obviously, this gives $\tau_0,\tau_{g+1}\in U$. If $\nu\in\{1,\dots,g\}$, then 
$\tau_\nu=\kappa_{\nu,1}+1$ by Proposition \ref{taudit}, which implies that 
\[
\tau_\nu\notin U\text{ for }\nu\in\{1,\dots,g\}.
\]
\end{rem}

\begin{prop}\label{key1}
Let $1\le i\le j\le n$ and let $\nu\in\{0,\dots,g\}$. Then
\[
[X_i\cdot X_j]_\nu=\min\{x_{i,\gamma_\nu}\rho_{j,\gamma_\nu},x_{j,\gamma_\nu}\rho_{i,\gamma_\nu}\}=
\left\lbrace
\begin{array}{ll}
x_{j,\gamma_\nu}\rho_{i,\gamma_\nu},&\text{if } i\in U; \\
x_{i,\gamma_\nu}\rho_{j,\gamma_\nu},&\text{if } i\notin U. 
\end{array}
\right.
\] 
Moreover, if $\gamma_{\nu+1}\le i,j$, then 
$[X_i\cdot X_j]_\nu=x_{i,\gamma_\nu}\rho_{j,\gamma_\nu}=x_{j,\gamma_\nu}\rho_{i,\gamma_\nu}$.
\end{prop}

\begin{proof}
Let us show that it is enough to consider the case $\gamma_\nu\le i\le j\le\gamma_{\nu+1}$.
Indeed, suppose first that $i<\gamma_\nu$. Then $x_{i,k}=0$ for 
$\gamma_\nu\le k\le n$. Subsequently, $x_{i,\gamma_\nu}=0=\rho_{i,\gamma_\nu}$, and
\[
[X_i\cdot X_j]_\nu=X_i^{(\gamma_\nu,\gamma_{\nu+1}]}\cdot X_j^{(\gamma_\nu,\gamma_{\nu+1}]}=0
=x_{i,\gamma_\nu}\rho_{j,\gamma_\nu}
=x_{j,\gamma_\nu}\rho_{i,\gamma_\nu}.
\] 
So the claim is clear in this case, and we may assume $\gamma_\nu\le i$.

If $\gamma_{\nu+1}\le j$, then we obtain by using Proposition \ref{tr} 
\[
[X_i\cdot X_j]_\nu=[X_i\cdot X_j^{\le\gamma_{\nu+1}}]_\nu=x_{j,\gamma_{\nu+1}}[X_i\cdot X_{\gamma_{\nu+1}}]_\nu.
\]
It follows from Corollary \ref{ax} that  
\[
\begin{array}{lll}
x_{i,\gamma_\nu}\rho_{j,\gamma_\nu}&=&x_{j,\gamma_{\nu+1}}x_{i,\gamma_\nu}\rho_{\gamma_{\nu+1},\gamma_\nu}\\
x_{j,\gamma_\nu}\rho_{i,\gamma_\nu}&=&x_{j,\gamma_{\nu+1}}x_{\gamma_{\nu+1},\gamma_\nu}\rho_{i,\gamma_\nu},
\end{array}
\]
and so we are reduced to the case $j=\gamma_{\nu+1}$.
Similarily, if $\gamma_{\nu+1}<i$ (so that $\gamma_{\nu+1}<j$), then
\[
[X_i\cdot X_j]_\nu=[X_i^{\le\gamma_{\nu+1}}\cdot X_j^{\le\gamma_{\nu+1}}]_\nu=
x_{i,\gamma_{\nu+1}}x_{j,\gamma_{\nu+1}}[X_{\gamma_{\nu+1}}^2]_\nu.
\]
and
\[
\begin{array}{lll}
x_{j,\gamma_\nu}\rho_{i,\gamma_\nu}&=
&x_{i,\gamma_{\nu+1}}x_{j,\gamma_{\nu+1}}(x_{\gamma_{\nu+1},\gamma_\nu}\rho_{\gamma_{\nu+1},\gamma_\nu})
\\
x_{i,\gamma_\nu}\rho_{j,\gamma_\nu}&=
&x_{i,\gamma_{\nu+1}}x_{j,\gamma_{\nu+1}}(x_{\gamma_{\nu+1},\gamma_\nu}\rho_{\gamma_{\nu+1},\gamma_\nu}).
\end{array}
\]
Then we are reduced to the case $i=j=\gamma_{\nu+1}$.

Let $(a_1,\dots,a_n)$ denote the point basis of $\fa$ and let $\gamma_\nu\le i\le j\le\gamma_{\nu+1}$. It follows from 
Proposition \ref{tg} that if $k<u$, then $a_{k-1}>a_k$ exactly when $x_{u,k-1}>x_{u,k}$. We may rewrite the sequence 
$x_{u,\gamma_\nu},\dots,x_{u,i}$ of multiplicities for any $u\ge i$ as
\begin{equation}\label{tele11}
\overbrace{s^{u}_1,\dots,s^{u}_1}^{r_1\text{ times }},
\overbrace{s^{u}_2,\dots,s^{u}_2}^{r_2\text{ times }},\dots,
\overbrace{s^{u}_\lambda,\dots,s^{u}_\lambda}^{r_{\lambda}\text{ times }}, 
\overbrace{s^{u}_{\lambda+1},\dots,s^{u}_{\lambda+1}}^{r'\text{ times }}, 
\end{equation}
where $s^u_1>\cdots>s^u_\lambda\ge s^u_{\lambda+1}$ and $s^u_\lambda>s^u_{\lambda+1}$ whenever $u>i$. 
Note that the equality $s^i_\lambda=s^i_{\lambda+1}$ can take place if and only if $r'=1$.
As in Equation (\ref{sr}), we get for any $1\le\mu\le\lambda$
\begin{equation}\label{tele22}
r_\mu s^u_\mu=s^u_{\mu-1}-s^u_{\mu+1},
\end{equation}
where 
\[
s^u_0=x_{u,\gamma_\nu}+\cdots+x_{u,\tau_{\nu+1}^{\fp_u}}.
\]
Because $\gamma_\nu\le i\le u$, Proposition \ref{tg} yields $\tau_{\nu+1}^{\fp_u}=\min\{u,\tau_{\nu+1}\}$ so that
\[
s^u_0=x_{u,\gamma_\nu}+\cdots+x_{u,\tau_{\nu+1}}=x_{u,\gamma_\nu}+\rho_{u,\gamma_\nu}.
\]
It follows from Proposition \ref{1} that for $\lambda>0$ and $u\ge i$
\begin{equation}\label{tele33}
s^i_\lambda=r's^i_{\lambda+1},\text{ while }
s^u_\lambda\ge r's^u_{\lambda+1}
\end{equation}
This holds true also in the case $\lambda=0$. Indeed, if $\lambda=0$, then 
this is a direct consequence of the definition of $s^u_0$.

Grouping the terms in $[X_i\cdot X_j]_\nu$ by using Equation (\ref{tele11}) gives now
\begin{equation}\label{tele1}
\begin{array}{lll}
[X_i\cdot X_j]_\nu&=&\displaystyle\sum_{k=\gamma_\nu}^{i}x_{i,k}x_{j,k}-x_{i,\gamma_\nu}x_{j,\gamma_\nu}\\
&=&\displaystyle\sum_{\mu=1}^{\lambda}r_\mu s^i_\mu s^j_\mu+r's^i_{\lambda+1}s^j_{\lambda+1}
-s^i_1s^j_1.
\end{array}
\end{equation}
Assume that $i\in U$ in which case $\lambda$ is even. Using Equation 
(\ref{tele22}) yields 
\[
\begin{array}{lll}
\displaystyle\sum_{\mu=1}^{\lambda}r_\mu s^i_\mu s^j_\mu
&=&(s^i_0-s^i_2)s^j_1+s^i_2(s^j_1-s^j_3)+\cdots\vspace{-2pt}\\
&&\phantom{(s^i_0}+(s^i_{\lambda-2}-s^i_\lambda)s^j_{\lambda-1}+s^i_\lambda(s^j_{\lambda-1}-s^j_{\lambda+1})\bigskip\\
&=&s^i_0s^j_1-s^i_\lambda s^j_{\lambda+1}
\end{array}
\]
and
\[
\begin{array}{lll}
\displaystyle\sum_{\mu=1}^{\lambda}r_\mu s^i_\mu s^j_\mu
&=&(s^j_0-s^j_2)s^i_1+s^j_2(s^i_1-s^i_3)+\cdots\vspace{-2pt}\\
&&\phantom{(s^i_0}+(s^j_{\lambda-2}-s^j_\lambda)s^i_{\lambda-1}+s^j_\lambda(s^i_{\lambda-1}-s^i_{\lambda+1})\bigskip\\
&=&s^j_0s^i_1-s^j_\lambda s^i_{\lambda+1}.\smallskip
\end{array}
\]
By Equation (\ref{tele1}) we then obtain
\[
\begin{array}{lll}
[X_i\cdot X_j]_\nu&
=&s^i_0s^j_1-s^i_\lambda s^j_{\lambda+1}+r's^i_{\lambda+1}s^j_{\lambda+1}-s^i_1s^j_1\vspace{3pt}\\
&=&s^j_0s^i_1-s^j_\lambda s^i_{\lambda+1}+r's^i_{\lambda+1}s^j_{\lambda+1}-s^i_1s^j_1.
\end{array}
\]
Furthermore, by Equation (\ref{tele33})
\[
r's^i_{\lambda+1}s^j_{\lambda+1}-
s^i_{\lambda}s^j_{\lambda+1}=0\text{ and }r's^i_{\lambda+1}s^j_{\lambda+1}-s^i_{\lambda+1}s^j_{\lambda}\le 0.
\]
Therefore 
\begin{equation}\label{sisi}
[X_i\cdot X_j]_\nu=(s^i_0-s^i_1)s^j_1\le (s^j_0-s^j_1)s^i_1.
\end{equation}

Similarily, if $i\notin U$, then $\lambda$ is odd and 
\[
\begin{array}{lll}
\displaystyle\sum_{\mu=1}^{\lambda}r_\mu s^i_\mu s^j_\mu
&=&(s^i_0-s^i_2)s^j_1+s^i_2(s^j_1-s^j_3)+\cdots\vspace{-2pt}\\
&&\phantom{(s^i_0}+s^i_{\lambda-1}(s^j_{\lambda-2}-s^j_\lambda)+(s^i_{\lambda-1}-s^i_{\lambda+1})s^j_\lambda\bigskip\\
&=&s^i_0s^j_1-s^i_{\lambda+1}s^j_\lambda.
\end{array}
\]
and 
\[
\begin{array}{lll}
\displaystyle\sum_{\mu=1}^{\lambda}r_\mu s^i_\mu s^j_\mu
&=&(s^j_0-s^j_2)s^i_1+s^j_2(s^i_1-s^i_3)+\cdots\vspace{-2pt}\\
&&\phantom{(s^i_0}+s^j_{\lambda-1}(s^i_{\lambda-2}-s^i_\lambda)+(s^j_{\lambda-1}-s^j_{\lambda+1})s^i_\lambda\bigskip\\
&=&s^j_0s^i_1-s^j_{\lambda+1}s^i_\lambda.\smallskip
\end{array}
\]
As above by using Equations (\ref{tele33}) and (\ref{tele1}) this yields that 
\[
[X_i\cdot X_j]_\nu=(s^j_0-s^j_1)s^i_1\le s^i_1(s^i_0-s^j_1).
\]
Because $(s^u_0,s^u_1)=(x_{u,\gamma_\nu}+\rho_{u,\gamma_\nu},x_{u,\gamma_\nu})$ for $u=i,j$, this together 
with Equation (\ref{sisi}) gives the claim.
\end{proof}

\begin{cor}\label{key11}
Assume $1\le i\le j\le n$. Set $\eta:=\gamma_\nu$ and $\gamma:=\gamma_{\nu+1}$, where $\nu$ is such that 
$\eta<i\le\gamma$ whenever $i>1$, and $\nu=0$ if $i=1$. With the notation above,  
\[
X_i\cdot X_j=x_{i,\eta}x_{j,\eta}X_{\eta}^2+[X_i\cdot X_j]_\nu.
\] 
Especially, if $i=\gamma\le j$, then 
\[
X_\gamma\cdot  X_j=x_{\gamma,\eta}(x_{j,\eta}X_{\eta}^2+\rho_{j,\eta})
=x_{j,\eta}(x_{\gamma,\eta}X_{\eta}^2+\rho_\nu). 
\]
\end{cor}

\begin{proof}
Obviously, $X_i\cdot X_j=X_i\cdot X_j^{\le\eta}+X_i\cdot X_j^{>\eta}$. Because $i\le\gamma$, we get 
$X_i\cdot X_j^{>\eta}=[X_i\cdot X_j]_\nu$, and by Proposition \ref{tr} 
\[
X_i\cdot X_j^{\le\eta}=X_i^{\le\eta}\cdot X_j^{\le\eta}=x_{i,\eta}x_{j,\eta}X_{\eta}^2.
\] 
Thus the first claim is clear, and the second claim follows from Proposition \ref{key1}. 
\end{proof}

\begin{cor}\label{key22}
Suppose that $\eta:=\gamma_\nu\le j\le k\le\gamma_{\nu+1}=:\gamma$ for some $\nu\in\{0,\dots,g\}$. Set  
\[
\sigma_1(j,k):=\frac{x_{j,\eta}}{x_{k,\eta}}\text{ and }
\sigma_2(j,k):=\frac{x_{j,\eta}X_\eta^2+\rho_{j,\eta}}{x_{k,\eta} X_\eta^2+\rho_{k,\eta}}.
\]
For any $i\in\{1,\dots,n\}$ the following equalities hold:
\begin{itemize}
\item[i)] If $i\le j$, then
\[
\frac{X_i\cdot X_j}{X_i\cdot X_k}=
\left\lbrace
\begin{array}{ll}
\displaystyle
 \sigma_1(j,k) &\text{if }i\in U\vspace{3pt}\text{ or }i\le\eta;\\
\displaystyle
 \sigma_2(j,k)&\text{if }i\notin U\text{ and }i>\eta.
\end{array}
\right.
\]
\item[ii)]
If $j<i$, then
\[
\begin{split}
\frac{X_i\cdot X_j}{X_i\cdot X_k}
=\left\lbrace
\begin{array}{ll}
\displaystyle
 \sigma_1(j,k)
&\text{if }j\notin U,\text{ and }
k\ge i\in U\text{ or }i>k\notin U;\vspace{3pt}\\
\displaystyle
 \sigma_2(j,i)\sigma_1(i,k)
&\text{if }j\in U,\text{ and }k\ge i\in U\text{ or }i>k\notin U;\vspace{3pt}\\
\displaystyle
 \sigma_1(j,i)\sigma_2(i,k)
&\text{if }j\notin U,\text{ and }k\ge i\notin U\text{ or }i>k\in U;\vspace{3pt}\\
\displaystyle
 \sigma_2(j,k)
&\text{if }j\in U,\text{ and }k\ge i\notin U\text{ or }i>k\in U.\vspace{3pt}\\
\end{array}
\right.
\end{split}
\]
\item[iii)] If $k'+1>k$ is free, then $X_i\cdot X_j/X_i\cdot X_k$ is constant for $i\ge k'$
\item[iv)] For every $1\le i\le n$
\[
\sigma_{3-\upsilon}(j,k)
\ge\frac{X_i\cdot X_j}{X_i\cdot X_k}
\ge\sigma_\upsilon(j,k)
=\frac{X_k\cdot X_j}{X_k^2},
\]
where $\upsilon=1$ for $j\notin U$ and $\upsilon=2$ for $j\in U$.
\end{itemize}
\end{cor}

\begin{proof} 
Before embarking the proof we first make two observations.
Assume that $i\le\eta$. As $\eta\le j\le k$ we then get by using Proposition \ref{tr} that 
\begin{equation}\label{xxxa}
\frac{X_i\cdot X_j}{X_i\cdot X_k}=\frac{X_i\cdot X_j^{\le\eta}}{X_i\cdot X_k^{\le\eta}}
=\frac{x_{j,\eta}X_i\cdot X_\eta}{x_{k,\eta}X_i\cdot X_\eta}=\frac{x_{j,\eta}}{x_{k,\eta}}=\sigma_1(j,k).
\end{equation}
Suppose then that $i\ge\eta$. By using Corollary \ref{key11} we get 
\begin{equation}\label{xxxx}
\frac{X_i\cdot X_j}{X_i\cdot X_k}=
\frac{x_{i,\eta}x_{j,\eta} X_\eta^2+[X_i\cdot X_j]_\nu}{x_{i,\eta}x_{k,\eta} X_\eta^2+[X_i\cdot X_k]_\nu}.
\end{equation}

(i) Proposition \ref{key1} implies that if $j\ge i\in U$, then $[X_i\cdot X_j]_\nu=x_{j,\eta}\rho_{i,\eta}$ and 
$[X_i\cdot X_k]_\nu=x_{k,\eta}\rho_{i,\eta}$, Furthermore, if $j\ge i\notin U$, then 
$[X_i\cdot X_j]_\nu=x_{i,\eta}\rho_{j,\eta}$ and $[X_i\cdot X_k]_\nu=x_{i,\eta}\rho_{k,\eta}$. This together with 
Equations (\ref{xxxa}) and (\ref{xxxx}) proves (i). 

(ii) We have $i\ge\eta$, and again, we make use of Equation (\ref{xxxx}). Now $i>j$, and thus by Proposition 
\ref{key1} $[X_i\cdot X_j]_\nu$ is either $x_{j,\eta}\rho_{i,\eta}$ or $x_{i,\eta}\rho_{j,\eta}$, depending on whether 
$j\notin U$ or $j\in U$. Proposition \ref{key1} shows also that $[X_i\cdot X_k]_\nu=x_{k,\eta}\rho_{i,\eta}$ when 
$k\ge i\in U$ or $i>k\notin U$, and $[X_i\cdot X_k]_\nu=x_{i,\eta}\rho_{k,\eta}$ when 
$k\ge i\notin U$ or $i>k\in U$. Putting all this together we obtain (ii).

(iii) If $k'+1>k$ and $k'+1$ is free, then we get 
by Proposition \ref{tr} for every $i\ge k'$
\[
\frac{X_i\cdot X_j}{X_i\cdot X_k}=\frac{X_i\cdot X_j^{\le k'}}{X_i\cdot X_k^{\le k'}}
=\frac{X_i^{\le k'}\cdot X_j}{X_i^{\le k'}\cdot X_k}
=\frac{x_{i,k'}X_{k'}\cdot X_j}{x_{i,k'}X_{k'}\cdot X_k}
=\frac{X_{k'}\cdot X_j}{X_{k'}\cdot X_k}. 
\]

(iv). We may assume $j<k$, since the case $j=k$ is trivial. 
Choosing $i=k$ in Equation (\ref{xxxx}) we see that 
\[
\frac{X_k\cdot X_j}{X_k^2}=
\frac{x_{j,\eta}x_{k,\eta}X_\eta^2+[X_k\cdot X_j]_\nu}{x_{k,\eta}x_{k,\eta}X_\eta^2+x_{k,\eta}\rho_{k,\eta}}, 
\]
where $[X_k\cdot X_j]_\nu=\min\{x_{k,\eta}\rho_{j,\eta},x_{j,\eta}\rho_{k,\eta}\}$ by Proposition \ref{key1}. Thereby
\begin{equation}\label{www}
\frac{X_k\cdot X_j}{X_k^2}=\sigma_\upsilon(j,k)\le\sigma_{3-\upsilon}(j,k),
\end{equation}
where $\upsilon=1$ if $j\notin U$ and $\upsilon=2$ if $j\in U$. Clearly, the claim holds, if $i$ is such that 
\[
\frac{X_i\cdot X_j}{X_i\cdot X_k}\in\left\{\sigma_1(j,k),\sigma_2(j,k)\right\},
\]
Let us check the remaining cases. Suppose first that 
$i>j\in U$, and $k\ge i\in U$ or $i>k\notin U$. By Proposition 
\ref{key1} we know that  
$[X_i\cdot X_j]_\nu=x_{i,\eta}\rho_{j,\eta}\le x_{j,\eta}\rho_{i,\eta}$ and 
$[X_i\cdot X_k]_\nu=x_{k,\eta}\rho_{i,\eta}\le x_{i,\eta}\rho_{k,\eta}$.
Then
\[
\frac{x_{i,\eta}x_{j,\eta}X_\eta^2+x_{i,\eta}\rho_{j,\eta}}{x_{i,\eta}x_{k,\eta}X_\eta^2+x_{i,\eta}\rho_{k,\eta}}\le
\frac{x_{i,\eta}x_{j,\eta}X_\eta^2+x_{i,\eta}\rho_{j,\eta}}{x_{i,\eta}x_{k,\eta}X_\eta^2+x_{k,\eta}\rho_{i,\eta}}\le
\frac{x_{i,\eta}x_{j,\eta}X_\eta^2+x_{j,\eta}\rho_{i,\eta}}{x_{i,\eta}x_{k,\eta}X_\eta^2+x_{k,\eta}\rho_{i,\eta}},
\]
in other words,
\[
\sigma_2(j,k)\sigma_1(i,i)\le\sigma_2(j,i)\sigma_1(i,k)\le
\sigma_2(i,i)\sigma_1(j,k).
\]
Note that $\sigma_\upsilon(i,i)=1$ for $\upsilon\in\{1,2\}$ and $X_i\cdot X_j/X_i\cdot X_k=\sigma_2(j,i)\sigma_1(i,k)$ 
by (ii). This together with Equation (\ref{www}) gives the claim.

Suppose next that $i>j\notin U$, and $k\ge i\notin U$ or $i>k\in U$. According to Proposition \ref{key1} we have 
$[X_i\cdot X_j]_\nu=x_{j,\eta}\rho_{i,\eta}\le x_{i,\eta}\rho_{j,\eta}$ and 
$[X_i\cdot X_k]_\nu=x_{i,\eta}\rho_{k,\eta}\le x_{k,\eta}\rho_{i,\eta}$.
This gives
\[
\frac{x_{i,\eta}x_{j,\eta}X_\eta^2+x_{j,\eta}\rho_{i,\eta}}{x_{i,\eta}x_{k,\eta}X_\eta^2+x_{k,\eta}\rho_{i,\eta}}\le
\frac{x_{i,\eta}x_{j,\eta}X_\eta^2+x_{j,\eta}\rho_{i,\eta}}{x_{i,\eta}x_{k,\eta}X_\eta^2+x_{i,\eta}\rho_{k,\eta}}\le
\frac{x_{i,\eta}x_{j,\eta}X_\eta^2+x_{i,\eta}\rho_{j,\eta}}{x_{i,\eta}x_{k,\eta}X_\eta^2+x_{i,\eta}\rho_{k,\eta}},
\]
which says
\[
\sigma_1(j,k)\sigma_2(i,i)\le\sigma_1(j,i)\sigma_2(i,k)\le
\sigma_1(i,i)\sigma_2(j,k).
\]
Again, (ii) shows that $X_i\cdot X_j/X_i\cdot X_k=\sigma_1(j,i)\sigma_2(i,k)$, which together with Equation (\ref{www}) 
gives the claim. Thus the proof is complete.
\end{proof}

\begin{prop}\label{key2}
Let $0\neq Z=r_1X_1+\cdots+r_nX_n$, where $r_1,\dots,r_n\in\mathbb N$. 
We then have for any $j\le k$.  
\[
\frac{Z\cdot X_j}{Z\cdot X_k}\ge\frac{X_j\cdot X_k}{X_k^2}. 
\]
\end{prop}

\begin{proof}
Assume first that $Z=X_i$ is a row of $P^{-1}$. If $i>k$, then 
\[
\frac{Z\cdot X_j}{Z\cdot X_k}=\frac{X_i^{\le k}\cdot X_j}{X_i^{\le k}\cdot X_k}.
\] 
Proposition \ref{tr} yields $X_i^{\le k}=x_{i,k}X_k+\varrho X_h$, where 
$\varrho=x_{i,h}-(x_{i,h+1}+\cdots+x_{i,k})$ in the case $k+1$ is a satellite to $h$ and otherwise $\varrho=0$. 
Let us recall the following elementary fact: If $a,b,c,d,e,f \in \mathbb N$ such that $bdf \neq 0$, then
\begin{equation}\label{abc1}
\displaystyle\frac{a}{b}\sim\frac{e}{f},\frac{c}{d}\sim\frac{e}{f}\Rightarrow\frac{a+c}{b+d}\sim\frac{e}{f},
\end{equation}
where $\sim$ is one of the relations $=$, $<$ or $>$.
Applying this gives
\[
\frac{X_h\cdot X_j}{X_h\cdot X_k}\ge\frac{X_k\cdot X_j}{X_k^2}\Longrightarrow\frac{Z\cdot X_j}{Z\cdot X_k}=
\frac{x_{i,k}X_k\cdot X_j+\varrho X_h\cdot X_j}{x_{i,k}X_k^2+\varrho X_h\cdot X_k}\ge\frac{X_k\cdot X_j}{X_k^2}.
\] 
This shows that it is enough to consider the case $i\le k$. 
Moreover, we may suppose $i\le j$, because  
\[
\frac{X_i\cdot X_j}{X_i\cdot X_k}\ge\frac{X_k\cdot X_j}{X_k^2}
\Longleftrightarrow\frac{X_i\cdot X_j}{X_k\cdot X_j}=\frac{X_i\cdot X_k}{X_k^2},
\]
and so, if $j<i\le k$, then we may simply switch the roles of $X_i$ and $X_j$.

By the above, we may restrict ourselves to the case $i\le j\le k$. Let $\nu\in\{0,\dots,g\}$ be such that 
$\gamma_\nu\le j\le\gamma_{\nu+1}$ and write $\gamma:=\gamma_{\nu+1}$. If $k\le\gamma$, then the claim is clear by 
Corollary \ref{key22} (iv). Especially, 
\begin{equation}\label{zut}
\frac{X_i\cdot X_j}{X_i\cdot X_\gamma}\ge\frac{X_\gamma \cdot X_j}{X_\gamma^2}. 
\end{equation}
If $k>\gamma$, then according to Proposition \ref{tr} we have  
\[
X_i\cdot X_k=X_i\cdot X_k^{\le\gamma}=x_{k,\gamma}X_i\cdot X_\gamma\text{ and }X_k\cdot X_j=X_k^{\le\gamma}\cdot X_j
=x_{k,\gamma}X_\gamma\cdot X_j.
\] 
By Proposition \ref{tr} we know that 
\[
x_{k,\gamma}^2X_\gamma^2=(X_k^{\le\gamma})^2<X_k^2.
\] 
Then we get by using these and Equation (\ref{zut})
\[
\frac{Z\cdot X_j}{Z\cdot X_k}=\frac{X_i\cdot X_j}{x_{k,\gamma}X_i\cdot X_\gamma}
\ge\frac{X_\gamma\cdot X_j}{x_{k,\gamma}X_\gamma^2}=
\frac{x_{k,\gamma}X_\gamma\cdot X_j}{x_{k,\gamma}^2 X_\gamma^2}\ge\frac{X_k\cdot X_j}{X_k^2}.
\]
Thereby the claim holds for any row $Z$ of $P^{-1}$.

Suppose then that $Z=r_1X_1+\cdots+r_nX_n$ for some  $(r_1,\dots,r_n)\in\mathbb N^n\smallsetminus\{0\}$. By the above 
\[
\frac{X_i\cdot X_j}{X_i\cdot X_k}\ge\frac{X_k\cdot X_j}{X_k^2}
\]
for every $i=1,\dots,n$. Applying Equation (\ref{abc1}) we obtain
\[
\frac{Z\cdot X_j}{Z\cdot X_k}\ge\frac{X_k\cdot X_j}{X_k^2}
\]
as desired.
\end{proof}

\begin{rem}\label{val}
Note that by Equation (\ref{d}) we have $X_i\cdot X_j=-\hat E_i\cdot\hat E_j$. Moreover, by setting 
$v_i:=\ord_{\ga_i}$ we 
may write 
\[
\hat E_i=\sum_{k=1}^nv_k(\fp_i)E_k=\mathbf V(\fp_i)E,
\] 
where $\fp_i$ is the simple complete ideal 
in $\ga$ containing $\fa$ and having the point basis $X_i$ (cf. Remark \ref{VF}). Because 
$E_k\cdot\hat E_j=-\delta_{k,j}$, we obtain
\[
X_i\cdot X_j=-\hat E_i\cdot\hat E_j=v_j(\fp_i)=v_i(\fp_j),
\]
where the last equality (known as \textit{reciprocity}, see \cite[p. 247, Proposition 21.4]{L}) is now obvious.
Therefore Proposition \ref{key2} says especially for any $j<i$ and for any $k$ that
\[
\frac{v_k(\fp_j)}{v_k(\fp_i)}\ge\frac{v_i(\fp_j)}{v_i(\fp_i)}.
\]
\end{rem}

In the sequel we write for any vector $X=(x_1,\dots,x_n)$ 
\begin{equation}\label{SX}
\Sigma X:=x_1+\cdots+x_n.
\end{equation}

\begin{rem}\label{kexi}
Observe that if $X_i$ is a row of the inverse of the proximity matrix of $\fa$, then 
$\Sigma X_i=-K\cdot\hat E_i=-K\cdot X_i E^*$, where $K=E^*_1+\cdots+E^*_n$ is the canonical divisor 
(see page \pageref{can}).
\end{rem}

\begin{prop}\label{key0}
Let $X_i$ be a row of the inverse of the proximity matrix of $\fa$. Let 
$\mu$ be such that $\gamma_\mu<i\le\gamma_{\mu+1}$ unless $i=1$ in which case we set $\mu:=0$. 
Then  
\[
\Sigma X_i^{>\gamma_k}+1-x_{i,\gamma_{k}}=\rho_{i,\gamma_{k}}+\cdots+\rho_{i,\gamma_\mu}.
\] 
for every $k=0,\dots,\mu$ 
\end{prop}

\begin{proof} 
Clearly, 
\[
\Sigma X_i^{>\gamma_k}=\Sigma X_i^{(\gamma_k,\gamma_{k+1}]}+\cdots+\Sigma X_i^{(\gamma_g,\gamma_{g+1}]}.
\] 
If $\mu<\nu$, then $i\le\gamma_\nu$, and we see that $\Sigma X_i^{(\gamma_\nu,\gamma_{\nu+1}]}=0$. Thus 
\[
\Sigma X_i^{>\gamma_k}=\Sigma X_i^{(\gamma_k,\gamma_{k+1}]}+\cdots+\Sigma X_i^{(\gamma_\mu,\gamma_{\mu+1}]}. 
\] 
We observe that the claim holds if 
\begin{equation*}\label{keya}
\Sigma X_i^{(\gamma_\nu,\gamma_{\nu+1}]}=
\left\lbrace
\begin{array}{ll}
x_{i,\gamma_\nu}+\rho_{i,\gamma_\nu}-x_{i,\gamma_{\nu+1}}&\text{ for }\nu<\mu;\\
x_{i,\gamma_\mu}+\rho_{i,\gamma_\mu}-1&\text{ for }\nu=\mu. 
\end{array}
\right.
\end{equation*}
It follows from Proposition \ref{tr} that  
$\Sigma X_i^{(\gamma_\nu,\gamma_{\nu+1}]}=x_{i,\gamma_{\nu+1}}\Sigma X_{\gamma_{\nu+1}}^{(\gamma_\nu,\gamma_{\nu+1}]}$, 
while Corollary \ref{ax} yields 
$x_{i,\gamma_\nu}+\rho_{i,\gamma_\nu}-x_{i,\gamma_{\nu+1}}=
x_{i,\gamma_{\nu+1}}(x_{\gamma_{\nu+1},\gamma_\nu}+\rho_\nu-1)$.
Subsequently, it is enough to verify that
\[
\Sigma X_i^{(\gamma_\mu,\gamma_{\mu+1}]}=x_{i,\gamma_\mu}+\rho_{i,\gamma_\mu}-1
\]
for every $i$ (and then especially in the case $i=\gamma_{\mu+1}$).

Consider the transform $\fb:=\fp_i^{\ga_{\gamma_\mu}}$. Recall that the point basis of $\fb$ corresponds to 
$X_i^{[\gamma_\mu,i]}$. Moreover, it follows from Proposition \ref{tg} that there are no terminal satellites 
between $\gamma_\mu$ and $i\le\gamma_{\mu+1}$, and therefore $\bar\Gamma_\fb=\{i\}$. By Proposition \ref{rhow} we may 
reduce to the situation $i=n$ and $\bar\Gamma=\{n\}$.

As in Formula (\ref{sv}), the sequence of the multiplicities $a_1,\dots,a_n$ is represented by some positive integers 
$m=m_1$, $r_j=r_{1,j}$ and $s_j=s_{1,j}$ for $1\le j\le m$, where $a_1=s_1>\cdots>s_m=a_n$. By Equation (\ref{sr}) we 
have $r_js_j=s_{j-1}-s_{j+1}$ for every $j=1,\dots,m$. Therefore   
\[
a_1+\cdots+a_n=r_1s_1+\cdots+r_ms_m+s_m=s_0+s_1-s_m-s_{m+1}+s_m, 
\]
where $s_0=a_1+\cdots+a_{\tau_1}$ and $s_{m+1}=s_m=1$. Then 
\[
a_2+\cdots+a_n=s_0-s_m=a_1+\rho_{n,\gamma_0}-1.
\] 
\end{proof}

\section{Multiplier ideals and jumping numbers}

Let $\fa$ be an ideal in a two-dimensional regular local ring $\ga$. Let $\mathcal X\rightarrow\Spec(\ga)$ be a 
\textit{log-resolution} of $\fa$, i.e., a projective birational morphism $\pi:\mathcal X\rightarrow\Spec(\ga)$ such 
that $\mathcal X$ regular and $\fa\cO_\mathcal X=\cO_\mathcal X(-D)$ for an effective Cartier divisor $D$ on 
$\mathcal X$ with the property that $D+\textrm{Exc}(\pi)$ has simple normal crossing support. Here $\textrm{Exc}(\pi)$ 
denotes the sum of the exceptional divisors of $\pi$. Note that (\ref{res}) is always a log-resolution. Recall that the 
relative canonical sheaf $\omega_{\mathcal X}$ can be defined as the dual of the relative Jacobian sheaf 
$\cJ_{\mathcal X}$ (cf. \cite[p. 203, (2.3)]{L4}). The \textit{canonical divisor}\label{can} $K:=K_{\mathcal X}$ of 
$\mathcal X$ is the unique exceptional divisor on $\mathcal X$ for which $\cO_{\mathcal X}(K)=\omega_{\mathcal X}$.

\begin{defn} \label{multip}
For a non-negative rational number $c$, the multiplier ideal $\cJ(\fa^c)$ is defined to be the ideal 
\[
\cJ(\fa^c):=\Gamma\left(\mathcal X,\cO_{\mathcal X}\left(K-\left\lfloor cD\right\rfloor\right)\right)\subset\ga, 
\]
where $D\in\Lambda$ is the effective divisor satisfying $\fa\cO_\mathcal X=\cO_{\mathcal X}(-D)$ and 
$\left\lfloor cD\right\rfloor$ denotes the integer part of $cD$. 
\end{defn}

\begin{lem}\label{Lauf}
Let $\fa$ be a simple complete ideal of finite colength in a two-dimensional regular local ring $\ga$. Then the base 
points of the multiplier ideal $\cJ(\fa^c)$ are among the base points of $\fa$ for every non-negative rational exponent 
$c$.
\end{lem}

\begin{proof} 
Let $\mathcal X\rightarrow\Spec(\ga)$ be the resolution of $\fa$ as in (\ref{res}), and let $D$ be the antinef 
divisor satisfying $\fa=\Gamma(\mathcal X,\cO_{\mathcal X}(-D))$. Then, by (\ref{GL})
\[
\cJ(\fa^c)=\Gamma(\mathcal X,\cO_{\mathcal X}(-(\left\lfloor cD \right\rfloor -K)))
=\Gamma(\mathcal X,\cO_{\mathcal X}(-(\left\lfloor cD \right\rfloor -K)^{\sim})). 
\] 
Since the antinef closure is antinef, we observe by Proposition \ref{nef} that $\cJ(\fa^c)$ generates an invertible 
$\cO_{\mathcal X}$-ideal.
\end{proof}

\begin{defn}\label{Hyppy}
Let $\fa$ be an ideal in a two-dimensional regular local ring $\ga$. By 
\cite[Lemma 1.3]{ELSV} there is an increasing discrete sequence $0=\xi_0<\xi_1<\xi_2<\cdots$ of rational numbers 
$\xi_i$ characterized by the properties that $\cJ(\fa^{c})=\cJ(\fa^{\xi_i})\text{ for }c\in[\xi_i,\xi_{i+1})$, while 
$\cJ(\fa^{\xi_{i+1}})\subsetneq\cJ(\fa^{\xi_i})$ for every $i$. The numbers $\xi_1,\xi_2,\dots$, are called the jumping 
numbers of $\fa$. We set 
\[
\cH_{\fa}=\{\xi_i\mid i=1,2,\dots\}. 
\]
\end{defn}

\begin{rem}
For practical reasons we don't consider $0$ as a jumping number in contrary to \cite[Definition 1.4]{ELSV}. Clearly, 
this is no restriction. Note that if $\fa=\ga$, then $\cJ(\fa^c)=\ga$ for every $c$, which means that the set of the 
jumping numbers is empty.
\end{rem}

\begin{defn}\label{ND}
Let $\fa,\fb$ be ideals of finite colength in a two-dimensional regular local ring $\ga$.
We define the \textit{log-canonical threshold of $\fa$ with respect to $\fb$} to be 
\[
c_\fb=c_\fb^\fa:=\inf\{c\in\mathbb Q_{>0}\mid\cJ(\fa^{c})\nsupseteq\fb\}.
\]
Note that if $\fb=\ga$, then $c_\fb$ is the usual log-canonical threshold.
\end{defn}

\begin{rem}\label{cbr}
By \cite[Theorem 11.1.1]{Lz} $\cJ(\fa^c)=\fa\cJ(\fa^{c-1})$, when $c\ge2$. If $\fa$ is proper, then we may 
find $c\gg0$ such that $\cJ(\fa^{c})\nsupseteq\fb_R$, and so 
\[
\{c\in\mathbb Q_{>0}\mid\cJ(\fa^{c})\nsupseteq\fb\}\neq\emptyset.
\] 
It follows from Definition \ref{Hyppy} that 
$c_\fb=\xi_i\in\mathbb Q$ for some $i=1,2\dots$, provided that $\fa$ is 
proper. If $\fa=\ga$, then the above set is empty and $c_\fb=\infty$ for any $\fb$.
\end{rem}

Let $I$ denote the point basis of $\fa$, and let $P^{-1}$ be the inverse of the proximity matrix of $\fa$ with the 
rows $X_1,\dots, X_n$. For an arbitrary vector $R = (r_1,\dots,r_n) \in \mathbb N^n$, we set 
\[
\fb_{R}:=\prod_{i=1}^n\fp_i^{r_i}, 
\] 
where $\fp_i$ is the simple complete $v$-ideal containing $\fa$ and having the point basis $X_i$. We write
\[
c_R=c_R^\fa:=c_{\fb_R}^\fa.
\]

\begin{prop}\label{Ha}
Let $\ga$ be a two-dimensional regular local ring and let $\fa \subset \ga$ be a simple complete ideal of finite 
colength. Let $\cH_{\fa}$ denote the set of all jumping numbers of the ideal $\fa$. Then 
\[
\cH_{\fa}=\{c_R\in\mathbb Q_{>0}\mid R\in\mathbb N^n\}.
\]
\end{prop}

\begin{proof} 
We may assume $\fa\neq\ga$. Take $R\in\mathbb N^n$, and let $\fb_R$ be as above. As we observed in Remark \ref{cbr}, 
$c_R=\xi_i$ for some positive integer $i$. Hence
\[
\cH_{\fa}\supset\{c_R\in\mathbb Q_{\ge 0}\mid R\in\mathbb N^n\}.
\] 
To show the opposite inclusion, take a jumping number $\xi_i$ ($i > 0$). By Definition \ref{Hyppy} we get 
$\cJ(\fa^{c})=\cJ(\fa^{\xi_{i-1}})\nsubseteq\cJ(\fa^{\xi_i})$ for any $c\in[\xi_{i-1},\xi_{i})$. By Lemma \ref{Lauf} 
the base points of $\cJ(\fa^{c})$ are among the base points of $\fa$. Then $\cJ(\fa^{c})=\fb_{R}$ for some 
$R\in\mathbb N^n$, which means that $\xi_i=c_R$, i.e., 
\[
\cH_{\fa}\subset\{c_R\in\mathbb Q_{\ge 0}\mid R\in\mathbb N^n\}.
\]
\end{proof}

\begin{cor}
Let $\fa$ be a simple complete ideal of finite colength in a two-dimensional regular local ring $\ga$ and take 
$R\in\mathbb N^n$. Then the set $\{\fb_{\nu}\mid\nu\in\mathbb N^n,c_{\nu}=c_R\}$ has the largest element containing all 
the others. Furthermore, the set $\{\cJ(\fa^{c})\mid c<c_R\}$ has the least element, and these two coincide.
\end{cor}

\begin{proof}
By Proposition \ref{Ha} $c_R=\xi_i$ for some $i\in\mathbb Z_+$. Then by Definition \ref{Hyppy} we obtain 
$\cJ(\fa^{\xi_{i-1}})=\min\{\cJ(\fa^{c})\mid c<c_R\}$. On the other hand, as we observed above 
$\cJ(\fa^{\xi_{i-1}})=\fb_{\mu}$ for some $\mu\in\mathbb N^n$, and clearly $c_{\mu}=\xi_i=c_R$. If $\nu\in\mathbb N^n$ 
is such that $c_\nu=c_R$, then $\cJ(\fa^{c})\supset\fb_{\nu}$ for every $c<c_R$. Especially 
$\cJ(\fa^{\xi_{i-1}})=\fb_{\mu}\supset\fb_{\nu}$, and therefore $\fb_{\mu}=\max\{\fb_{\nu}\mid c_{\nu}=c_R\}$.
\end{proof}

\section{Key lemmas}

In order to determine the set of the jumping numbers, we make use of Proposition \ref{Ha}. For the main proofs we shall 
need a few technical results which are mostly gathered in this section. As above, $\fa$ is a simple complete ideal of 
finite colength in a two-dimensional regular local ring $\ga$ having the resolution (\ref{res}) and the base points 
$\ga=\ga_1\subset\cdots\subset\ga_n$. Let $P$ denote the proximity matrix and $I=(a_1\dots,a_n)$ the point basis of 
$\fa$ and let $\fb_R$ and $c_R$ be as in Definition \ref{ND}. 
Recall that $\cH_{\fa}=\{c_R \mid R \in \mathbb N^n\}$ according to Proposition \ref{Ha}.

\begin{ntt}
Let $X$ and $Y$ be row vectors of $P^{-1}$. For any $R \in \mathbb N^n$ set 
$\tilde{R}:=RP^{-1}$ and write
\begin{equation*}\label{RYX}
R_Y[X]:=\displaystyle\frac{\tilde{R}\cdot X+\Sigma X+1}{X\cdot Y},
\end{equation*}
where $\Sigma X$ is as defined in Equation (\ref{SX}).
In the following, we usually write $R[X]:=R_I[X]$.
\end{ntt}

\begin{prop}\label{ci}
Let $X_1,\dots,X_n$ be the rows of $P^{-1}$. Then for any $R \in \mathbb N^n$ 
\[
c_R=\min\big\{R[X_i]\mid i=1,\dots,n\big\}.
\]
\end{prop}

\begin{proof}
Set $D=\hat E_n$ so that $\fa\cO_\mathcal X=\cO_{\mathcal X}(-D)$. 
We have by Definition \ref{multip} and Equation (\ref{GL})
\[
\cJ(\fa^c)=\Gamma(\mathcal X,\cO_\mathcal X(-(\left\lfloor cD \right\rfloor-K)))=
           \Gamma(\mathcal X,\cO_\mathcal X(-(\left\lfloor cD \right\rfloor-K)^{\sim})), 
\]
where $K$ denotes the canonical divisor and $(\left\lfloor cD\right\rfloor-K)^{\sim}$ stands for the antinef closure of 
$\left\lfloor cD\right\rfloor-K$. By Proposition \ref{nef} 
$\cJ(\fa^c)\cO_\mathcal X=\cO_\mathcal X(-(\left\lfloor cD \right\rfloor-K)^{\sim})$ is invertible. Also 
$\fb_R\subset\ga$ is invertible. Since $\cJ(\fa^c)$ and $\fb_R$ are both complete, we have 
$\cJ(\fa^c)\supset\fb_R$ exactly when $\cJ(\fa^c)\cO_\mathcal X\supset\fb_R\cO_\mathcal X$, 
which is equivalent to
\[
\cO_\mathcal X(-(\left\lfloor cD \right\rfloor-K)^{\sim})\supset\cO_\mathcal X(-\mathbf V(\fb_R)E),
\]
where $E$ is as in Equation (\ref{E}). This means that
\[
(\left\lfloor cD\right\rfloor-K)^{\sim}\le \mathbf V(\fb_R)E. 
\]
Because $(\left\lfloor cD \right\rfloor-K)^{\sim}$ is by definition the minimal antinef divisor 
satisfying $(\left\lfloor cD \right\rfloor-K)\le(\left\lfloor cD \right\rfloor-K)^{\sim}$, we see that this holds 
if and only if
\[
\left\lfloor cD \right\rfloor-K \le \mathbf V(\fb_R)E.
\]
Recall that $\mathbf V(\fb_R)=(\tilde R\cdot X_1,\dots,\tilde R\cdot X_n)$ and 
$\mathbf V(\fa)=(I\cdot X_1,\dots,I\cdot X_n)$ by Equation (\ref{d*}). Similarly, 
$K=E_1^*+\cdots+E_n^*=(\Sigma X_1,\cdots,\Sigma X_n)E$. 
Therefore the inequality above is equivalent to
\[
\left\lfloor cX_i\cdot I \right\rfloor - \Sigma X_i \le \tilde R\cdot X_i
\]
for every $i=1,\dots,n$. So 
$\cJ(\fa^{c}) \nsupseteq \fb_R$ exactly if $\left\lfloor cX_i\cdot I\right\rfloor>\tilde R\cdot X_i+\Sigma X_i$
for some $i=1,\dots,n$, or equivalently,
\[
cX_i\cdot I \ge \tilde R\cdot X_i + \Sigma X_i + 1.
\]
for some $i=1,\dots,n$. This means that $c \ge R[X_i]$ for some $i = 1,\dots,n$. Now $c_R$ is by Definition \ref{ND} 
the smallest rational number $c$, for which $\cJ(\fa^{c}) \nsupseteq \fb_R$. Thus we get the claim. 
\end{proof}

The remaining problem is to tell for which $X_i$ we reach the minimum of the $R[X_i]$:s.
Let $\gamma_0,\dots,\gamma_{g+1}$ be as given in Notation \ref{gamma}. We shall now present a few 
useful equalities and equivalences which will be needed in calculating and comparing the $R[X_i]$:s.

\begin{lem}\label{cix}
Let $U$ be as in Equation (\ref{U}) and suppose that $u\in U$. 
When $u>1$ let $\nu$ satisfy $\gamma_\nu<u\le\gamma_{\nu+1}$, whereas $\nu=0$ if $u=1$.
Set $\eta:=\gamma_\nu$ and $\gamma:=\gamma_{\nu+1}$. Furthermore, take $R=(r_1,\dots,r_n)\in\mathbb N^n$ and write
\[
\xi:=\sum_{j\in J}r_j \rho_{j,\eta}\text{ and }\zeta:=\sum_{j\notin J}r_j x_{j,\eta}
\]
where we set $J:=\{1,\dots,\eta-1\}\cup\{j\mid\eta\le j<u\text{ and }j\in U\}$. Write also 
\[
\delta:=\tilde R\cdot X_\eta+\Sigma X_\eta+1-(\zeta+1)X_\eta^2.
\] 
Suppose that $\eta\le k\le\gamma$. Then
\vspace{5pt}
\begin{itemize}
\item[i)]  $\Sigma X_k+1=x_{k,\eta}(\Sigma X_\eta+1)+\rho_{k,\eta};$\vspace{6pt}
\item[ii)] $\tilde R\cdot X_k^{>\eta}\le\rho_{k,\eta}\zeta+x_{k,\eta}\xi;$\vspace{6pt}
\item[iii)]$\tilde R\cdot X_k+\Sigma X_k+1\le(\delta+\xi)x_{k,\eta}+(\zeta+1)(x_{k,\eta}X_\eta^2+\rho_{k,\eta});$
            \vspace{6pt}
\item[iv)] If $u=1$, then $\delta=1;$
\item[v)]  If $u>1$, then $R[X_\eta]\sim R[X_u]\Leftrightarrow\delta\sim x_{u,\eta}X_\eta^2\xi:\rho_{u,\eta}$,
           \vspace{6pt}
\end{itemize}
where $\sim$ denotes any of the relations $=,<$ or $>$. Moreover, the equality holds in ii) and iii) if $k=u$.
\end{lem}

\begin{proof}
(i) Clearly, 
\[
\Sigma X_k+1=\Sigma X_k^{\le\eta}+x_{k,\eta}+\Sigma X_k^{>\eta}+1-x_{k,\eta}. 
\]
By using Proposition \ref{tr} we get $\Sigma X_k^{\le\eta}+x_{k,\eta}=x_{k,\eta}(\Sigma X_\eta+1)$, and further, 
by Proposition \ref{key0} we obtain $\Sigma X_k^{>\eta}+1-x_{k,\eta}=\rho_{k,\eta}$.

(ii) We first observe that 
\[
\tilde R\cdot X_k^{>\eta}
=\sum_{j=1}^nr_jX_j\cdot X_k^{>\eta}.
\]
Because $k\le\gamma$, we have $X_j\cdot X_k^{>\eta}=[X_k\cdot X_j]_\nu$. Proposition \ref{key1} then yields that
\begin{equation*}\label{rixxxxx}
\begin{array}{rlllll}
\tilde R\cdot X_k^{>\eta}
&=&\displaystyle\sum_{j=1}^nr_j[X_k\cdot X_j]_\nu\vspace{3pt}\\
&=&\displaystyle\sum_{j=1}^nr_j\min\{x_{k,\eta}\rho_{j,\eta},x_{j,\eta}\rho_{k,\eta}\}\vspace{3pt}\\
&=&\displaystyle\sum_{j\in J}r_j\min\{x_{k,\eta}\rho_{j,\eta},x_{j,\eta}\rho_{k,\eta}\}+
\sum_{j\notin J}r_j\min\{x_{k,\eta}\rho_{j,\eta},x_{j,\eta}\rho_{k,\eta}\}\vspace{3pt}\\
&\le&\displaystyle\sum_{j\in J}r_jx_{j,\eta}\rho_{k,\eta}+\sum_{j\notin J}r_j\rho_{j,\eta}x_{k,\eta}\vspace{3pt}\\
&=&\rho_{k,\eta}\zeta+x_{k,\eta}\xi.
\end{array}
\end{equation*}
Let us then show that the equality holds here if $k=u$. By the above it is enough to prove that 
\[
[X_u\cdot X_j]_\nu=
\left\lbrace
\begin{array}{ll}
x_{u,\eta}\rho_{j,\eta},&\text{if } j\in J;\\
x_{j,\eta}\rho_{u,\eta},&\text{if } j\notin J. 
\end{array}
\right.
\]
Indeed, suppose first that $j\in J$. Then either $j<\eta$ or $j\in U$ with $\eta\le j<u$. In the first case we have 
$[X_u\cdot X_j]_\nu=0=x_{u,\eta}\rho_{j,\eta}$, while in the second case $[X_u\cdot X_j]_\nu=x_{u,\eta}\rho_{j,\eta}$ 
by Proposition \ref{key1} as wanted. Suppose then that $j\notin J$. If $\eta\le j<u$, then $j\notin U$, and Proposition 
\ref{key1} gives $[X_u\cdot X_j]_\nu=x_{j,\eta}\rho_{u,\eta}$. If $j\ge u$, then the same holds, as $u\in U$.

(iii) Using Proposition \ref{tr} we have 
\[
\tilde R\cdot X_k=\tilde R\cdot X_k^{\le\eta}+\tilde R\cdot X_k^{>\eta}
=x_{k,\eta}\tilde R\cdot X_\eta+\tilde R\cdot X_k^{>\eta}.
\] 
By i) $\Sigma X_k+1=x_{k,\eta}(\Sigma X_\eta+1)+\rho_{k,\eta}$. Putting these together yields
\[
\begin{array}{lll}
\tilde R\cdot X_k+\Sigma X_k+1\vspace{2pt}&=&
x_{k,\eta}(\tilde R\cdot X_\eta+\Sigma X_\eta+1)+\rho_{k,\eta}+\tilde R\cdot X_k^{>\eta}\vspace{2pt}
\\&=&x_{k,\eta}\delta+(\zeta+1)x_{k,\eta}X_\eta^2+\rho_{k,\eta}+\tilde R\cdot X_k^{>\eta}.
\end{array}
\]
By ii) $\tilde R\cdot X_k^{>\eta}\le\rho_{k,\eta}\zeta+x_{k,\eta}\xi$, where the equality holds if $k=u$. Thus the 
claim is clear.

(iv) If $u=1$, then $\eta=1$, which implies that $X_\eta^2=1=\Sigma X_\eta$. Furthermore, we have $J=\emptyset$, which 
yields
\[
\zeta=\sum_{j\notin J}r_jx_{j,\eta}=\sum_{j=1}^nr_jX_j\cdot X_\eta=\tilde R\cdot X_\eta.
\]
Thus we see that $\delta=\tilde R\cdot X_\eta+\Sigma X_\eta+1-\zeta-1=1$.

(v) Observe that $X_u\cdot I=a_\eta x_{u,\eta}X_\eta^2+[X_u\cdot I]_\nu$ by Corollary \ref{key11}, where 
$I=(a_1,\dots,a_n)=X_n$. Because $u\in U$, we have $[X_u\cdot I]_\nu=a_\eta\rho_{u,\eta}$ by Proposition \ref{key1}. 
Therefore we obtain by using iii) 
\[
\begin{array}{lll}
R[X_u]&=&\displaystyle\frac{\tilde R\cdot X_u+\Sigma X_u+1}{X_u\cdot I}\vspace{3pt}\\
&=&\displaystyle
\frac{x_{u,\eta}(\delta+\xi)+(\zeta+1)(x_{u,\eta}X_\eta^2+\rho_{u,\eta})}{X_u\cdot I}\vspace{3pt}\\
&=&\displaystyle
\frac{x_{u,\eta}(\tilde R\cdot X_\eta+\Sigma X_\eta+1)+x_{u,\eta}\xi+(\zeta+1)\rho_{u,\eta}}
{a_\eta x_{u,\eta}X_\eta^2+a_\eta\rho_{u,\eta}}.
\end{array}
\]
Noting that $X_\eta\cdot I=X_\eta\cdot I^{\le\eta}=a_\eta X_\eta^2$ by Proposition \ref{tr} we then see
that $R[X_\eta]\sim R[X_u]$ is equivalent to
\begin{equation}\label{abo}
\frac{\tilde R\cdot X_\eta+\Sigma X_\eta+1}{X_\eta^2}\sim
\frac{x_{u,\eta}(\tilde R\cdot X_\eta+\Sigma X_\eta+1)+x_{u,\eta}\xi+(\zeta+1)\rho_{u,\eta}}
{x_{u,\eta}X_\eta^2+\rho_{u,\eta}}.
\end{equation}
If $u>1$, in which case $u>\eta$, then $\rho_{u,\eta}>0$. An elementary fact similar to Equation (\ref{abc1}) says that 
if $a,b,c,d \in \mathbb N$ such that $bd \neq 0$, then
\begin{equation}\label{abc2}
\displaystyle\frac{a+c}{b+d}\sim\frac{a}{b}\Leftrightarrow\frac{c}{d}\sim\frac{a}{b}
\end{equation}
where $\sim$ is one of the relations $=$, $<$ or $>$. By using this we see that Equation (\ref{abo}) is further 
equivalent to 
\[
\frac{\tilde R\cdot X_\eta+\Sigma X_\eta+1}{X_\eta^2}\sim
\frac{x_{u,\eta}\xi+(\zeta+1)\rho_{u,\eta}}
{\rho_{u,\eta}}.
\]
which is the same as
\[
\delta=\tilde R\cdot X_\eta+\Sigma X_\eta+1-(\zeta+1)X_\eta^2\sim\frac{x_{u,\eta}X_\eta^2\xi}{\rho_{u,\eta}}.
\]
The claim has thus been proven.
\end{proof}

Next we will show that the relevant indices in searching the minimum of the $R[X_i]$:s are in the set 
$\bar\Gamma=\Gamma^*\cup\{n\}$, where $\Gamma^*$ denotes the set of the indices corresponding to the star vertices of 
the dual graph as before.

\begin{lem}\label{c=RHo}
Let $X_1,\dots,X_n$ be the rows of $P^{-1}:=(x_{i,j})_{n\times n}$. Let $i\in\{1,\dots,n\}$. For $1<i\le n$, let 
$\nu\in\{0,\dots,g\}$ satisfy $\gamma_\nu<i\le\gamma_{\nu+1}$, whereas $\nu=0$ if $i=1$. Write $\eta:=\gamma_\nu$ and 
$\gamma:=\gamma_{\nu+1}$. Then 
\[
R[X_\eta]\ge R[X_i]\Rightarrow R[X_i]\ge R[X_\gamma].
\]
Moreover, in the case $n>1$ we have $R[X_1]>R[X_{\gamma_1}]$.
\end{lem}

\begin{proof}
It follows from Proposition \ref{tr} that
\begin{equation*}\label{tito}
\frac{X_i\cdot X_n}{X_\gamma\cdot X_n}
=\frac{X_i^{\le\gamma}\cdot X_n}{X_\gamma^{\le\gamma}\cdot X_n}
=\frac{X_i\cdot X_n^{\le\gamma}}{X_\gamma\cdot X_n^{\le\gamma}}
=\frac{x_{n,\gamma}X_i\cdot X_\gamma}{x_{n,\gamma}X_\gamma^2}
=\frac{X_i\cdot X_\gamma}{X_\gamma^2},
\end{equation*}
which shows that 
\[
R[X_i]
=\frac{\tilde R\cdot X_i+\Sigma X_i+1}{X_i\cdot X_n}
\ge
\frac{\tilde R\cdot X_\gamma +\Sigma X_\gamma+1}{X_\gamma\cdot X_n}
=R[X_\gamma]
\] 
is equivalent to
\begin{equation}\label{taito}
\frac{\tilde R\cdot X_i+\Sigma X_i+1}{\tilde R\cdot X_\gamma +\Sigma X_\gamma+1}\ge\frac{X_i\cdot X_\gamma}{X_\gamma^2}. 
\end{equation}
Clearly, $R[X_i]>R[X_\gamma]$ if and only if the inequality here is strict.

Suppose first that $i\notin U$. According to Corollary \ref{key22} (iv) we have  
\[
\frac{x_{i,\eta}X_\eta^2+\rho_{i,\eta}}{x_{\gamma,\eta}X_\eta^2+\rho_{\gamma,\eta}}=\sigma_2(i,\gamma)
\ge\sigma_1(i,\gamma)=\frac{x_{i,\eta}}{x_{\gamma,\eta}}=\frac{X_i\cdot X_\gamma}{X_\gamma^2}.
\]
As $1\in U$ we must have $1\le\eta<i\le\gamma$. Then $\rho_{\gamma,\eta}>0$, and Equation (\ref{abc2}) implies
\[
\frac{\rho_{i,\eta}}{\rho_{\gamma,\eta}}
\ge 
\frac{x_{i,\eta}}{x_{\gamma,\eta}}.
\]
It follows from Lemma \ref{cix} (i) that
\[
\frac{\Sigma X_i+1}{\Sigma X_\gamma+1}=
\frac{x_{i,\eta}(\Sigma X_\eta+1)+\rho_{i,\eta}}{x_{\gamma,\eta}(\Sigma X_\eta+1)+\rho_{\gamma,\eta}}.
\]
An application of Equation (\ref{abc2}) then gives 
\begin{equation*}
\frac{\Sigma X_i+1}{\Sigma X_\gamma+1}\ge
\frac{x_{i,\eta}}{x_{\gamma,\eta}}
=\frac{X_i\cdot X_\gamma}{X_\gamma^2}.
\end{equation*}
In the case $R=0$ this is the same as Inequality (\ref{taito}). If $R\neq 0$, then
Proposition \ref{key2} yields
\[
\frac{\tilde R\cdot X_i}{\tilde R\cdot X_\gamma}
\ge
\frac{X_i\cdot X_\gamma}{X_\gamma^2}.
\]
Applying Equation (\ref{abc1}) to these two inequalities implies Inequality (\ref{taito}).

Suppose then that $i\in U$. 
According to Corollary \ref{key22} (iv) we have 
\begin{equation}\label{tatu}
\frac{x_{i,\eta}}{x_{\gamma,\eta}}=\sigma_1(i,\gamma)
\ge\sigma_2(i,\gamma)=
\frac{x_{i,\eta}X_\eta^2+\rho_{i,\eta}}{x_{\gamma,\eta}X_\eta^2+\rho_{\gamma,\eta}}
=\frac{X_i\cdot X_\gamma}{X_\gamma^2}. 
\end{equation}
Choosing $u=i$ and $k=i,\gamma$ in Lemma \ref{cix} (iii) gives
\begin{equation}\label{piitu}
\frac{\tilde R\cdot X_i+\Sigma X_i+1}{\tilde R\cdot X_\gamma +\Sigma X_\gamma+1}
\ge
\frac{x_{i,\eta}(\delta+\xi)+(\zeta+1)(x_{i,\eta}X_\eta^2+\rho_{i,\eta})}
     {x_{\gamma,\eta}(\delta+\xi)+(\zeta+1)(x_{\gamma,\eta}X_\eta^2+\rho_{\gamma,\eta})}.
\end{equation}
By Lemma \ref{cix} (iv) and (v) $R[X_\eta]\ge R[X_i]$ implies $\delta\ge0$. Then $\delta+\xi\ge0$.
If $\delta+\xi=0$, this already proves Inequality (\ref{taito}) by Inequality (\ref{tatu}). 
Thus it remains to consider the case $\delta+\xi>0$. Then we may apply Equation (\ref{abc2}) to (\ref{tatu}) to get
\begin{equation}\label{janne}
\frac{x_{i,\eta}(\delta+\xi)+(\zeta+1)(x_{i,\eta}X_\eta^2+\rho_{i,\eta})}
     {x_{\gamma,\eta}(\delta+\xi)+(\zeta+1)(x_{\gamma,\eta}X_\eta^2+\rho_{\gamma,\eta})}
\ge
\frac{x_{i,\eta}X_\eta^2+\rho_{i,\eta}}{x_{\gamma,\eta}X_\eta^2+\rho_{\gamma,\eta}}.
\end{equation}
Combining this to Inequalities (\ref{tatu}) and (\ref{piitu}) shows that Inequality (\ref{taito}) holds.
The first claim has thus been proven.

The second claim follows from putting $i=1$ in the inequalities above. Note that Inequality
(\ref{tatu}) is now strict, which implies that so are Inequalities (\ref{janne}) and (\ref{taito}), too. Indeed, 
$1\in U$, and $n>1$ yields that $\gamma_1=\gamma>\eta=1$ so that $\rho_{\gamma,\eta}>0$ while $\rho_{i,\eta}=0$. 
\end{proof}

\begin{prop}\label{c=RH}
Let $X_1,\dots,X_n$ be the rows of $P^{-1}:=(x_{i,j})_{n\times n}$. Then 
\[
c_R=\min\big\{R[X_\gamma]\mid\gamma\in\bar\Gamma\big\}
\]
for any $R=(r_1,\dots,r_n)\in\mathbb N^n$. Moreover, $c_R<R[X_1]$ if $n>1$.
\end{prop}

\begin{proof}
According Proposition \ref{ci} we have $c_R=\min\{R[X_j]\mid j=1,\dots,n\}$. Hence the claim is obvious in the case 
$n=1$. 

Suppose that $n>1$. 
By Lemma 
\ref{c=RHo} we have $c_R\le R[X_{\gamma_1}]<R[X_1]$. Thus the last claim is clear, and furthermore, if $i$ is such 
that $c_R=R[X_i]$ for some $i$, then $i>1$. It follows that 
$\gamma_\nu<i\le\gamma_{\nu+1}$ for some $\nu\in\{0,\dots,g^*\}$. Because $R[X_{\gamma_\nu}]\ge R[X_i]$, Lemma 
\ref{c=RHo} yields $R[X_i]\ge R[X_{\gamma_{\nu+1}}]$, which further implies that $c_R=R[X_{\gamma_{\nu+1}}]$.
\end{proof}

\begin{lem}\label{R=R'} 
Suppose that $\gamma\in\Gamma^*$, and take any $R \in \mathbb N^n$ such that $c_R = R[X_{\gamma}]$. 
Then $c_{R^{<\gamma}} = R^{<\gamma}[X_{\gamma}]\le a_{\gamma}^{-1}$, and for some $m \in \mathbb N$, 
\[
c_R = c_{R^{<\gamma}} + \frac{m}{a_{\gamma}}.
\]
\end{lem}

\begin{proof}
Write $R = (r_1,\dots,r_n)$ and set $R':=R^{<\gamma}$. In order to show that $c_{R'}=R'[X_\gamma]$, it is 
by Proposition \ref{c=RH} enough to verify, that $R'[X_\gamma]\le R'[X_\eta]$ for any $\eta\in\bar\Gamma$.

Write $R=R'+R''$, where $R'':=R^{\ge\gamma}$. Furthermore, set $\tilde R':=R'P^{-1}$ and $\tilde R'':=R''P^{-1}$. 
Then 
\[
\tilde R=RP^{-1}= R'P^{-1}+R''P^{-1}=\tilde R'+\tilde R'',
\] 
and for any row $X$ of $P^{-1}$ we may write
\begin{equation}\label{h1}
R[X]=\frac{(\tilde R'+\tilde R'')\cdot X+\Sigma X+1}{I\cdot X} 
= R'[X]+\frac{\tilde R''\cdot X}{I\cdot X}. 
\end{equation} 
Assume that $\eta\in\bar\Gamma$. Proposition \ref{ci} implies that $c_R=R[X_\gamma] \le R[X_\eta]$. So 
\begin{equation*}\label{h4}
R'[X_\gamma] - R'[X_\eta] \le 
\frac{\tilde R''\cdot X_\eta}{I\cdot X_\eta} - \frac{\tilde R''\cdot X_\gamma}{I\cdot X_\gamma}=:\Delta.
\end{equation*}
We will show that $\Delta\le0$. 
Note that the claim is trivial, if $R=R'$. Thus we may presume that $R''\neq0$.

Suppose first that $\eta\le\gamma$. Then $X_\eta=X_\eta^{\le\gamma}$. Since 
$\tilde R''=\sum_{i = \gamma}^{n}r_iX_i$, we see by using Proposition \ref{tr} that for any $\eta\le\gamma$
\begin{equation}\label{h3}
\frac{\tilde R''\cdot X_\eta}{I\cdot X_\eta} = \frac{\tilde R''\cdot X_\eta^{\le \gamma}}{I\cdot X_\eta^{\le \gamma}} 
=\frac{\sum_{i=\gamma}^{n}r_iX_i^{\le \gamma}\cdot X_\eta}{I^{\le \gamma}\cdot X_\eta}
=\frac{\sum_{i=\gamma}^{n}r_ix_{i,\gamma}X_\gamma\cdot X_\eta}{a_\gamma X_\gamma\cdot X_\eta}
=\frac{m}{a_\gamma}, 
\end{equation}
where $m:=\sum_{i=\gamma}^{n}r_ix_{i,\gamma}$. Hence $\Delta=(m-m)/a_\gamma=0$, if $\eta\le\gamma$.

Assume next that $\eta>\gamma$. Then Proposition \ref{tr} implies that
\[
\frac{I\cdot X_\gamma}{I\cdot X_\eta}=\frac{I\cdot X_\gamma^{\le\eta}}{I\cdot X_\eta^{\le\eta}}
=\frac{I^{\le\eta}\cdot X_\gamma}{I^{\le\eta}\cdot X_\eta}
=\frac{a_\eta X_\eta\cdot X_\gamma}{a_\eta X_\eta^2}=\frac{X_\eta\cdot X_\gamma}{X_\eta^2}. 
\]
By Proposition \ref{key2} we have for every $i=1,\dots,n$
\[
\frac{X_i\cdot X_\gamma}{X_i\cdot X_\eta}\ge\frac{X_\eta\cdot X_\gamma}{X_\eta^2}. 
\] 
By applying Equation (\ref{abc1}) we then get 
\[
\frac{\tilde R''\cdot X_\gamma}{\tilde R''\cdot X_\eta}
=\frac{\sum_{i=\gamma}^{n}r_iX_i\cdot X_\gamma}{\sum_{i=\gamma}^{n}r_iX_i\cdot X_\eta}
\ge\frac{X_\eta\cdot X_\gamma}{X_\eta^2}=\frac{I\cdot X_\gamma}{I\cdot X_\eta}.
\]
Thus 
\[
\Delta=\frac{\tilde R''\cdot X_\eta}{I\cdot X_\eta}-\frac{\tilde R''\cdot X_\gamma}{I\cdot X_\gamma}\le0,
\]
and therefore $c_{R'}=R'[X_\gamma]$. Subsequently, by Equations (\ref{h1}) and (\ref{h3})
\[
c_R = R[X_\gamma]=R'[X_\gamma]+\frac{\tilde R''\cdot X_\gamma}{I\cdot X_\gamma}=c_{R'}+\frac{m}{a_\gamma}.
\]

It remains to show $c_{R'} \le 1/a_\gamma$. Set $i:=\min\{\eta\in\bar\Gamma\mid\gamma<\eta\}$. Observe that such an 
index exists, because by assumption $\gamma<n$. Recall that $i\in U$ (see Remark \ref{keyU}). Since 
$c_{R'}=R'[X_\gamma]$, we have  $R'[X_\gamma]\le R'[X_i]$. By Lemma \ref{cix} (v) this implies
\begin{equation}\label{1a}
\tilde R'\cdot X_\gamma+\Sigma X_\gamma+1-(\zeta+1)X_\gamma^2\le\frac{x_{i,\gamma}X_\gamma^2\xi}{\rho_{i,\gamma}}.  
\end{equation}
We have $\rho_{i,\gamma}\zeta+x_{i,\gamma}\xi=\tilde R'\cdot X_i^{>\gamma}$ by Lemma \ref{cix} (ii). On the other hand,
\[
\tilde R'\cdot X_i^{>\gamma}=\sum_{j=1}^{\gamma-1}r_jX_j\cdot X_i^{>\gamma}=0,
\]
Therefore Equation (\ref{1a}) gives
\[
\tilde R'\cdot X_\gamma+\Sigma X_\gamma+1\le 
\left(\frac{\rho_{i,\gamma}\zeta+x_{i,\gamma}\xi}{\rho_{i,\gamma}}+1\right)X_\gamma^2=X_\gamma^2. 
\]
Since $X_\gamma\cdot I=X_\gamma\cdot I^{\le\gamma}=a_\gamma X_\gamma^2$ by Proposition \ref{tr}, we obtain from this 
\[
R'[X_\gamma]=\frac{\tilde R'\cdot X_\gamma+\Sigma X_\gamma+1}{X_\gamma\cdot I}\le \frac{1}{a_\gamma}.
\] 
Because $c_{R'}=R'[X_\gamma]$, we get the claim. 
\end{proof}

\section{Jumping numbers of a simple ideal}

In this section we will give a formula for the jumping numbers of a simple complete ideal $\fa$ in a two-dimensional 
regular local ring $\ga$ in terms of the multiplicities of the point basis $I=(a_1,\dots,a_n)$ of $\fa$. 
Let $\gamma_0,\dots,\gamma_{g+1}$ and $\Gamma^*$ be as in Notation \ref{gamma}. Recall from Proposition \ref{star} that 
$\Gamma^*$ the set of indices corresponding to the stars of the associated dual graph. 
Before we proceed to the main theorem, we prove the following lemma:

\begin{lem}\label{1ain} 
Let $\ga$ be a two-dimensional regular local ring, and let $\fa$ be a simple complete $\fm_\ga$-primary ideal in 
$\ga$ having the point basis $I=(a_1,\dots,a_n)$. 
For $\nu=0,\dots,g$, write $\gamma_\nu=\eta$ and 
$\gamma_{\nu+1}=\gamma$, and set 
\[
b_\nu:=\frac{I\cdot I^{\le\gamma}}{a_{\eta}}.
\] 
Then 
\[
b_\nu=\frac{a_\gamma I\cdot X_\gamma}{a_\eta}=a_\gamma(x_{\gamma,\eta}X_{\eta}^2+\rho_\nu).
\]
Especially, $b_\nu$ is an integer, for which 
\[
\gcd\{a_\eta,b_\nu\}=a_\gamma.
\]
Moreover, $a_\eta\le b_\nu$, where the equality holds if and only if $n=1$.
\end{lem}

\begin{proof} Proposition \ref{tr} gives $I^{\le\gamma}=a_\gamma X_\gamma$. Therefore 
\[
b_\nu=\frac{I\cdot I^{\le\gamma}}{a_{\eta}}=\frac{a_\gamma I\cdot X_\gamma}{a_\eta}.
\] 
Using Corollary \ref{key11} we obtain 
$X_\gamma\cdot I=a_\eta (x_{\gamma,\eta}X_\eta^2+\rho_\nu)$, 
and so we get the first claim. Clearly, $b_\nu$ is an integer, and since 
$a_\eta^2\le I\cdot I^{\le\gamma}=a_\eta b_\nu$, we see that $a_\eta\le b_\nu$. 
Here the equality holds if only if $n=1$.
By Corollary \ref{ax} $a_\gamma x_{\gamma,\eta}=a_\eta$ and $a_\gamma\rho_\nu=\rho_{n,\eta}$. Then
\[
\begin{array}{lll}
\gcd\{a_\eta,b_\nu\}
&=&\gcd\{a_\eta,a_\gamma(x_{\gamma,\eta}X_\eta^2+\rho_\nu)\}\\
&=&\gcd\{a_\eta,a_\eta X_\eta^2+\rho_{n,\eta}\}\\
&=&\gcd\{a_\eta,\rho_{n,\eta}\}.
\end{array}
\]
By definition $\rho_{n,\eta}=a_{\eta+1}+\cdots+a_{\tau_{\nu+1}}$, and $a_\eta=\cdots=a_{\tau_{\nu+1}-1}$ by Proposition 
\ref{taudit}. Thus $\gcd\{a_\eta,\rho_{n,\eta}\}=\gcd\{a_\eta,a_{\tau_{\nu+1}}\}$, and the claim follows from 
Proposition \ref{aiaj1}. 
\end{proof}

\begin{rem}\label{Za}
The integers $a_1,b_0,\dots,b_g$ in fact coincide with the so called \textit{Zariski exponents} 
$\bar\gb_0,\dots,\bar\gb_{g+1}$. The Zariski exponents can be defined recursively as follows. Let 
$\gb'_1,\dots,\gb'_{g+1}$ be the Puiseux exponents (see Notation \ref{dpuis} and Remark \ref{puis}). Proposition 
\ref{aiaj1} gives that $\gcd\{a_{\gamma_{\nu-1}}+\cdots+a_{\tau_\nu},a_{\gamma_{\nu-1}}\}=a_{\gamma_\nu}$ for every 
$\nu=1,\dots,g+1$, and then by using Corollary \ref{ax} we see that
\[
\gb'_\nu=\frac{x_{\gamma_\nu,\gamma_{\nu-1}}+\cdots+x_{\gamma_\nu,\tau_\nu}}{x_{\gamma_\nu,\gamma_{\nu-1}}},
\]
where $\gcd\{x_{\gamma_\nu,\gamma_{\nu-1}}+\cdots+x_{\gamma_\nu,\tau_\nu},x_{\gamma_\nu,\gamma_{\nu-1}}\}=1$. Recall 
that by Corollary \ref{ax} $a_{\gamma_\nu}=x_{\gamma_{\nu+1},\gamma_\nu}\cdots x_{\gamma_{g+1},\gamma_g}$ for every 
$\nu=0,\dots,g$. Note that the integers $a_{\gamma_0},\dots,a_{\gamma_{g+1}}$ and 
$x_{\gamma_1,\gamma_0},\dots,x_{\gamma_{g+1},\gamma_{g}}$ are usually denoted by $e_0,\dots,e_{g+1}$ and 
$n_1,\dots,n_{g+1}$, respectively (cf. \cite[p. 130]{S1}). Write $n_0:=1$. Following \cite[Equation 6.1]{S1} 
(clearly, $e_{i+1}$ is a misprint in the cited equation), set 
\[
\bar\gb_0:=e_0,\text{ and }\bar\gb_\nu:=(\gb'_\nu-1)e_{\nu-1}+\bar\gb_{\nu-1}n_{\nu-1}\text{ for }1\le\nu\le g+1.
\] 

Let us prove that these are the integers $a_{\gamma_0},b_0\dots,b_g$. 
From the definition above we see that $\bar\gb_0=a_{\gamma_0}=I\cdot X_{\gamma_0}$. In order to verify that 
$\bar\gb_\nu=b_{\nu-1}$ for $1\le\nu\le g+1$, let us first show that $\bar\gb_\nu n_\nu=I\cdot X_{\gamma_\nu}$ also 
for $\nu>0$. Suppose that $\bar\gb_{\nu-1}n_{\nu-1}=I\cdot X_{\gamma_{\nu-1}}$ holds for some $1\le\nu\le g+1$. Remark 
\ref{rhob} shows that $(\gb'_\nu-1)e_{\nu-1}=\rho_{n,\gamma_{\nu-1}}$, and 
$I\cdot X_{\gamma_{\nu-1}}=a_{\gamma_{\nu-1}}X_{\gamma_{\nu-1}}^2$ by Proposition \ref{tr}. Subsequently, we obtain by 
using Corollary \ref{key11}
\begin{equation}\label{Zar}
\bar\gb_\nu n_\nu=n_\nu(I\cdot X_{\gamma_{\nu-1}}+\rho_{n,\gamma_{\nu-1}}) 
=x_{\gamma_{\nu},\gamma_{\nu-1}}(a_{\gamma_{\nu-1}}X_{\gamma_{\nu-1}}^2+\rho_{n,\gamma_{\nu-1}}) 
=I\cdot X_{\gamma_\nu}. 
\end{equation}
Thus an induction on $\nu$ shows that 
\[
\bar\gb_\nu=\frac{I\cdot X_{\gamma_\nu}}{x_{\gamma_{\nu},\gamma_{\nu-1}}}
=\frac{a_{\gamma_\nu}I\cdot X_{\gamma_\nu}}{a_{\gamma_{\nu-1}}}
\] 
for $1\le\nu\le g+1$, where the last equality follows from Corollary 
\ref{ax}. Together with Lemma \ref{1ain} this completes the proof.

We have also an alternative characterization for the Zariski exponents:
\[
\bar\gb_\nu=v_{\tau_\nu}(\fa)=v(\fp_{\tau_\nu})\hspace{25pt}(\nu=0,\dots,g+1).
\]
Clearly, this holds for $\nu=0$. Let us then verify this in the case $\nu>0$. Equation (\ref{Zar}) yields 
$\bar\gb_\nu=I\cdot X_{\gamma_{\nu-1}}+\rho_{n,\gamma_{\nu-1}}$ for $1\le\nu\le g+1$. 
On the other hand, since $x_{\tau_{\nu},i}=1$ for every $\gamma_{\nu-1}\le i\le\tau_{\nu}$ by Proposition \ref{tau}, we 
get $I\cdot X_{\tau_\nu}^{>\gamma_{\nu-1}}=\rho_{n,\gamma_{\nu-1}}$. An application of Proposition \ref{tr} then 
gives $I\cdot X_{\tau_\nu}^{\le\gamma_{\nu-1}}=I\cdot X_{\gamma_{\nu-1}}$. Subsequently, for $1\le\nu\le g+1$, 
\[
\bar\gb_\nu=I\cdot X_{\gamma_{\nu-1}}+\rho_{n,\gamma_{\nu-1}}=
I\cdot X_{\tau_\nu}^{\le\gamma_{\nu-1}}+I\cdot X_{\tau_\nu}^{>\gamma_{\nu-1}}=I\cdot X_{\tau_\nu}.
\]
Because $I\cdot X_{\tau_\nu}=v_{\tau_\nu}(\fa)=v(\fp_{\tau_\nu})$, as observed in Remark \ref{val}, we thus obtain the 
desired characterization.
\end{rem}

\begin{thm}\label{main} 
Let $\ga$ be a two-dimensional regular local ring, and let $\fa$ be a simple complete $\fm_\ga$-primary ideal in $\ga$ 
having the point basis $I=(a_1,\dots,a_n)$. Let $\Gamma^*=\{\gamma_1,\dots,\gamma_{g^*}\}$ the set of indices 
corresponding to the stars of the dual graph associated to $\fa$, and write $\gamma_0:=1$ and $\gamma_{g^*+1}:=n$. For 
$\nu=0,\dots,g^*$, set $$b_\nu:=\frac{I\cdot I^{\le \gamma_{\nu+1}}}{a_{\gamma_\nu}}$$ and then define for 
$s,t,m\in\mathbb N$ 
\[
H_\nu:=\left\lbrace \frac{s+1}{a_{\gamma_\nu}}+\frac{t+1}{b_\nu}+\frac{m}{a_{\gamma_{\nu+1}}}
\hspace{2pt}\middle |\hspace{2pt} s,t,m\in\mathbb N, 
\frac{s+1}{a_{\gamma_\nu}}+\frac{t+1}{b_\nu}\le\frac{1}{a_{\gamma_{\nu+1}}}\right\rbrace
\]
for $\nu=0,\dots,g^*-1$ and 
\[
H_{g^*}:=
\left\lbrace\frac{s+1}{a_{\gamma_{g^*}}}+\frac{t+1}{b_{g^*}}
\hspace{2pt}\middle|\hspace{2pt} s,t\in\mathbb N\right\rbrace.
\]
The set of the jumping numbers of the ideal $\fa$ is then  
\[
\cH_\fa=H_0\cup\cdots\cup H_{g^*}.
\]
\end{thm}

\begin{proof} 
If $n=1$, then $\fa=\fm_\ga$. Subsequently, $\cH_\fa=\{m\in\mathbb N\mid m>1\}$, and the case is clear. Thus we may 
suppose throughout the proof that $n>1$.

To see that every jumping number is in $H_0\cup\cdots\cup H_{g^*}$, recall first that by Proposition \ref{Ha} any 
element in $\cH_\fa$ is of the form $c_R$ for some $R\in\mathbb N^n$. Take an arbitrary 
$R=(r_1,\dots,r_n)\in\mathbb N^n$. 
By Proposition \ref{c=RH} there is 
$\nu\in\{0,\dots,g^*\}$ (recall that $\bar\Gamma=\{\gamma_1,\dots,\gamma_{g^*+1}\}$) such that  
\begin{equation}\label{dgg}
c_R=R[X_\gamma]<R[X_\eta],
\end{equation}
where $\gamma:=\gamma_{\nu+1}$ and $\eta:=\gamma_\nu$. As observed in Remark \ref{keyU}, $\gamma\in U$. 
Thus we may apply Lemma \ref{cix} (iii) to get
\begin{equation*}
R[X_\gamma]
=\frac{\tilde{R}\cdot X_\gamma+\Sigma X_\gamma+1}{I\cdot X_\gamma}
 =\frac{(x_{\gamma,\eta}X_\eta^2+\rho_\nu)(\zeta+1)+x_{\gamma,\eta} (\xi+\delta)}{I\cdot X_\gamma}. 
\end{equation*}
By Corollary \ref{ax} we have $a_\eta/a_\gamma=x_{\gamma,\eta}$. It then follows from Lemma \ref{1ain} that 
\begin{equation}\label{dg}
{I\cdot X_\gamma}=a_\eta(x_{\gamma,\eta}X_\eta^2+\rho_\nu)=x_{\gamma,\eta}b_\nu.
\end{equation}
Putting all together we get 
\begin{equation}\label{dg1}
c_R=\frac{\zeta+1}{a_\eta}+\frac{\xi+\delta}{b_\nu}. 
\end{equation}
Clearly, $\zeta+1$ is a positive integer, and by Lemma \ref{cix} (v) we see from Equation (\ref{dgg}) that $\delta$ 
and thereby also $\xi+\delta$ are positive integers. 

If $\nu=g^*$, then Equation (\ref{dg1}) proves that 
$c_R\in H_{g^*}$. If $\nu<g^*$, then it follows from Lemma \ref{R=R'} that $c_R$ is 
in $H_\nu$ exactly, when $c_{R^{<\gamma}}$ is. Thus we may assume $R=R^{<\gamma}$. This case is now clear by Equation 
(\ref{dg1}), as Lemma \ref{R=R'} guarantees that $c_{R^{<\gamma}}\le a_\gamma^{-1}$. Subsequently, 
\[
\cH_\fa\subset H_0\cup\cdots\cup H_{g^*}
\]
as wanted.

In order to prove the opposite inclusion, we first need two more lemmas.

\begin{lem} \label{syt} 
Let $\nu\in\{0,\dots,g\}$, where $g$ is the number of the terminal satellites of $\fa$, and let 
$t_{\nu+1},m_{\nu+1}\in\mathbb N$. Then there exists two sequences of pairs $(s_1,t_1),\dots,(s_\nu,t_\nu)$ and 
$(s_{\nu+2},m_{\nu+2}),\dots,(s_g,m_g)$ of non negative integers satisfying the following conditions: 
\begin{itemize}
\item[i)] for every $1\le i\le\nu$ 
\begin{equation*}
t_{i+1}+1+X_{\gamma_i}^2=(s_i+1)(x_{\gamma_i,\gamma_{i-1}}X_{\gamma_{i-1}}^2 + \rho_{i-1}) 
+ (t_{i}+1)x_{\gamma_i,\gamma_{i-1}}, 
\end{equation*} 
\item[ii)] for every $\nu<i<g$ we have $m_i=m_{i+1}x_{\gamma_{i+1},\gamma_i}+s_{i+1}$ and 
$s_{i+1}<x_{\gamma_{i+1},\gamma_i}$. Moreover, writing $s_{g+1}:=m_g$ we get for every $\nu<\mu\le g$ 
\[
m_\mu=\sum_{i=\mu}^g s_{i+1}x_{\gamma_i,\gamma_\mu}.
\]
\item[iii)] 
For any $k=\nu+2,\dots,g+1$, we have
\[
\Phi_k:=\frac{(m_{\nu+1}+1)\rho_{\gamma_k,\gamma_{\nu+1}}+\cdots+(m_{k-1}+1)\rho_{\gamma_k,\gamma_{k-1}}}
{a_{\gamma_{\nu+1}}\rho_{\gamma_k,\gamma_{\nu+1}}
+\cdots+a_{\gamma_{k-1}}\rho_{\gamma_k,\gamma_{k-1}}}
\ge 
\frac{m_{\nu+1}+1}{a_{\gamma_{\nu+1}}}.
\]
\end{itemize} 
\end{lem}

\begin{proof}
(i)
Using descending induction, suppose that $t_{i+1}\in\mathbb N$ is given for some $1<i\le\nu$. By Corollary \ref{key11} 
\[
x_{\gamma_i,\gamma_{i-1}}(x_{\gamma_i,\gamma_{i-1}}X_{\gamma_{i-1}}^2+\rho_{i-1})=X_{\gamma_i}^2. 
\]
Note that $x_{\gamma_i,\gamma_{i-1}}=a_{\gamma_{i-1}}/a_{\gamma_i}$ by Corollary \ref{ax}. Moreover, we also have 
$x_{\gamma_i,\gamma_{i-1}}X_{\gamma_{i-1}}^2+\rho_{i-1}=b_{i-1}/a_{\gamma_i}$ by Lemma \ref{1ain}, which then further 
yields 
\[
\gcd\{x_{\gamma_i,\gamma_{i-1}},x_{\gamma_i,\gamma_{i-1}}X_{\gamma_{i-1}}^2+\rho_{i-1}\}=
\frac{\gcd\{a_{\gamma_{i-1}},b_{i-1}\}}{a_{\gamma_i}}=1.
\] 
The existence of a pair $(s_i,t_i)$ now follows from Lemma \ref{syt2} below, and thereby we obtain a sequence 
$(s_1,t_1),\dots,(s_\nu,t_\nu)$.

(ii)
Given a non negative integer $m_i$ for $\nu<i<g$, we have $m_{i+1},s_{i+1} \in \mathbb N$ with 
$m_i=m_{i+1}x_{\gamma_{i+1},\gamma_i}+s_{i+1}$ and $s_{i+1}<x_{\gamma_{i+1},\gamma_i}$. Arguing inductively, 
the existence of the pairs $(s_{\nu+2},m_{\nu+2}),\dots,(s_{g},m_{g})$ is then clear. 
Observe that $m_g=s_{g+1}x_{\gamma_g,\gamma_g}$. Moreover, assuming that
\[
m_{\mu+1}=\sum_{i=\mu+1}^g s_{i+1}x_{\gamma_i,\gamma_{\mu+1}}
\]
holds for any $\nu<\mu<g$ we see by using Corollary \ref{ax} that 
\[
m_\mu=\sum_{i=\mu+1}^g s_{i+1}x_{\gamma_i,\gamma_{\mu+1}}x_{\gamma_{\mu+1},\gamma_\mu}+s_{\mu+1}
=\sum_{i=\mu}^g s_{i+1}x_{\gamma_i,\gamma_\mu}.
\]
Hence this holds for every $\nu<\mu\le g$.

(iii) 
To prove the last claim we first note that for every $\nu<i<g$
\[
 \frac{m_{i+1}+1}{a_{\gamma_{i+1}}}
=\frac{m_{i+1}x_{\gamma_{i+1},\gamma_i}+x_{\gamma_{i+1},\gamma_i}}{a_{\gamma_{i+1}}x_{\gamma_{i+1},\gamma_i}}.
\]
As $x_{\gamma_{i+1},\gamma_i}\ge s_{i+1}+1$ and $a_{\gamma_{i+1}}x_{\gamma_{i+1},\gamma_i}=a_{\gamma_i}$ by Corollary 
\ref{ax}, we see that
\[
 \frac{m_{i+1}+1}{a_{\gamma_{i+1}}}
\ge\frac{m_{i+1}x_{\gamma_{i+1},\gamma_i}+s_{i+1}+1}{a_{\gamma_{i+1}}x_{\gamma_{i+1},\gamma_i}}
=\frac{m_i+1}{a_{\gamma_i}}.
\]
Because $\rho_{\gamma_k,\gamma_i}>0$ for every $\nu<i<k$, the above implies 
\[
\frac{(m_i+1)\rho_{\gamma_k,\gamma_i}}{a_{\gamma_i}\rho_{\gamma_k,\gamma_i}}\ge
\frac{m_{\nu+1}+1}{a_{\gamma_{\nu+1}}}.
\]
The claim follows now from Equation (\ref{abc1}).
\end{proof}

\begin{lem} \label{syt2} 
If $a\le b$ are positive integers and $\gcd\{a,b\} = 1$, then for any positive integer $t$ there exist positive 
integers $u$ and $v$ such that $ua + vb = ab + t$.
\end{lem}

\begin{proof}
We can assume $a < b$, as the case $a=b=1$ is trivial, and clearly, we may reduce to the case $t<a$. 
Since $\gcd\{a,b\}=1$, we can find a positive integer $p$ such that 
$pa = 1 \mod b$. Then $tpa = t+qb$ for some integer $q$. 
We see that $b \nmid tp$, as otherwise $b\mid tpa-qb=t$, which is impossibile as $0<t<a<b$. 
Thus we may find integers $r$ and $u$ such that $tp = rb + u$ and 
$0 < u < b$.  We get $tpa = rba + ua = t + qb$, i.e., $ua = t + (q-ra)b$. As $0<ua$ and $t<b$, we see that 
$0 \le q-ra$. Since $u<b$ we get $ua = t +(q-ra)b < ab$. Especially, this gives $q-ra < a$. Set $v:= (r+1)a-q$. Then
we observe that $v>0$ and 
\[\begin{array}{lll}
ua+vb&=&ua+((r+1)a-q)b\\
&=&t+(q-ra)b+(r+1)ab-qb\\&=&t+ab.
\end{array}
\] 
\end{proof}

Choose any $\nu\in\{0,\dots,g^*\}$, and take an arbitrary element $c\in H_\nu$. Then there exist 
$s,t,m\in\mathbb N$ such that 
\[
c=\frac{s+1}{a_\eta}+\frac{t+1}{b_\nu}+\frac{m}{a_\gamma},\text{ moreover, }c\le\frac{m+1}{a_\gamma}\text{ if }\nu<g^*. 
\]
As above, write $\eta=\gamma_\nu$ and $\gamma=\gamma_{\nu+1}$. In order to prove that $c=c_R$ for some $R\in\mathbb N^n$, we 
shall proceed in three steps. 
\medskip

A) 
To begin with, we shall construct a suitable $R \in \mathbb N^n$.
Let $g$ be the number of the terminal satellites of $\fa$. Set $t_{\nu+1}:=t$ and $m_{\nu+1}:=m$, and let 
$(s_1,t_1),\dots,(s_\nu,t_\nu)$ and $(s_{\nu+2},m_{\nu+2}),\dots,(s_g,m_g)$ be sequences constructed as in Lemma 
\ref{syt}. From these we obtain a sequence $s_0,\dots,s_{g+1}\in\mathbb N$ by setting $s_{\nu+1}:=s$, $s_0:=t_1$ and 
$s_{g+1}:=m_g$. Define $R=(r_1,\dots,r_n)$ so that 
\[
r_i=
\left\lbrace
\begin{array}{ll}
s_\nu&\text{if }i=\tau_\nu\text{ for some }\nu=0,\dots,g+1;\\ 
0  &\text{otherwise}.
\end{array}
\right.
\]
We observe that 
\[
\tilde R = s_0X_{\tau_0}+\cdots+s_{g+1}X_{\tau_{g+1}}
=\tilde R_k + \tilde S_k
\]
for any $k=0,\dots,g+1$, where
\[
\tilde R_k := \sum_{i = 0}^k s_iX_{\tau_i}\text{ and }
\tilde S_k := \sum_{i = k+1}^{g+1} s_iX_{\tau_i}. 
\]

\medskip
B) 
We shall show $c=R[X_\gamma]$. For $k=0,\dots,g+1$, set 
\[
\psi(k):=\tilde R_k\cdot X_{\gamma_k}+\Sigma X_{\gamma_k}+1.
\] 
Note that 
\begin{equation}\label{pssw}
R[X_\gamma]=\frac{\tilde R\cdot X_\gamma+\Sigma X_\gamma+1}{I\cdot X_\gamma}
=\frac{\psi(\nu+1)}{I\cdot X_\gamma}+\frac{\tilde S_{\nu+1}\cdot X_\gamma}{I\cdot X_\gamma}.
\end{equation}
According to Proposition \ref{tr} we have $X_{\tau_{i+1}}^{\le\gamma_i}=x_{\tau_{i+1},\gamma_i}X_{\gamma_i}$, where 
$x_{\tau_{i+1},\gamma_i}=1$ by Proposition \ref{tau}. If $i\ge k$, then $\gamma_k\le\gamma_i$ so that
\[
X_{\tau_{i+1}}\cdot X_{\gamma_k}=X_{\tau_{i+1}}\cdot X_{\gamma_k}^{\le\gamma_i}
=X_{\tau_{i+1}}^{\le\gamma_i}\cdot X_{\gamma_k}
=X_{\gamma_i}\cdot X_{\gamma_k}.
\] 
Moreover, an application of Proposition \ref{tr} gives 
\[
I\cdot X_{\gamma_k}=I^{\le\gamma_i}\cdot X_{\gamma_k}=a_{\gamma_i}X_{\gamma_i}\cdot X_{\gamma_k}.
\] 
Therefore
\begin{equation*}
\frac{\tilde S_k\cdot X_{\gamma_k}}{I\cdot X_{\gamma_k}}
=\sum_{i=k}^{g}\frac{s_{i+1}X_{\tau_{i+1}}\cdot X_{\gamma_k}}{I\cdot X_{\gamma_k}} 
=\sum_{i=k}^{g}\frac{s_{i+1}X_{\gamma_{i}}\cdot X_{\gamma_k}}{a_{\gamma_i}X_{\gamma_{i}}\cdot X_{\gamma_k}}
=\sum_{i=k}^{g}\frac{s_{i+1}}{a_{\gamma_i}}.  
\end{equation*}
As $a_\gamma/a_{\gamma_i}=x_{\gamma_i,\gamma}$ by Corollary \ref{ax}, it follows from Lemma \ref{syt} (ii) that 
\begin{equation}\label{si}
\frac{\tilde S_{\nu+1}\cdot X_\gamma}{I\cdot X_\gamma}
=\sum_{i=\nu+1}^{g}\frac{s_{i+1}}{a_{\gamma_i}}
=\sum_{i=\nu+1}^{g}\frac{s_{i+1}x_{\gamma_i,\gamma}}{a_\gamma}
=\frac{m}{a_\gamma}.
\end{equation}
We aim to show that 
\begin{equation}\label{psi4}
\psi({\nu+1})=(t+1)x_{\gamma,\eta} + (s+1)(x_{\gamma,\eta}X_{\eta}^2+ \rho_\nu).
\end{equation}
Suppose for a moment that this holds.
As we already saw in Equation (\ref{dg}), 
$I\cdot X_\gamma=a_\eta(x_{\gamma,\eta}X_\eta^2+\rho_\nu)=x_{\gamma,\eta}b_\nu$. 
By Equations (\ref{pssw}) and (\ref{si}) we then have 
\[
\begin{array}{lll}
R[X_\gamma]&=&\displaystyle
\frac{(t+1)x_{\gamma,\eta}}{I\cdot X_\gamma}+\frac{(s+1)(x_{\gamma,\eta}X_{\eta}^2+ \rho_\nu)}
{I\cdot X_\gamma}+\frac{m}{a_\gamma}\vspace{3pt}\\
&=&\displaystyle\frac{(t+1)}{b_\nu}+\frac{(s+1)}{a_\eta}+\frac{m}{a_\gamma}=c
\end{array}
\]
as desired.

In order to verify Equation (\ref{psi4}), recall first that by Proposition \ref{key0} 
\[
\Sigma X_{\gamma_k}^{>\gamma_0}+1-x_{\gamma_k,\gamma_0}=\rho_{\gamma_k,\gamma_0}+\cdots+\rho_{\gamma_k,\gamma_{k-1}}.
\]
By Corollary \ref{ax} we get  
$\rho_{\gamma_k,\gamma_i}=x_{\gamma_k,\gamma_{i+1}}\rho_i$ for $0\le i<k$, and so  
\[
\Sigma X_{\gamma_k}+1=2x_{\gamma_k,\gamma_0}+x_{\gamma_k,\gamma_1}\rho_{0}+\cdots+x_{\gamma_k,\gamma_k}\rho_{k-1}.
\] 
Setting $\rho_{(-1)}:=2$ allows us to write
\begin{equation}\label{psi1}
\psi(k)=\tilde R_k\cdot X_{\gamma_k}+\Sigma X_{\gamma_k}+1 
        =\sum_{i=0}^k\big(s_iX_{\tau_i}\cdot X_{\gamma_k}+x_{\gamma_k,\gamma_i}\rho_{i-1}\big). 
\end{equation}
If $i<k$, then $\gamma_{k-1}\ge\gamma_i\ge\tau_i$, and Proposition \ref{tr} implies that  
\[
X_{\tau_i}\cdot X_{\gamma_k}=X_{\tau_i}\cdot X_{\gamma_k}^{\le\gamma_{k-1}}
=x_{\gamma_k,\gamma_{k-1}}X_{\tau_i}\cdot X_{\gamma_{k-1}}.
\] 
Furthermore, by Corollary \ref{ax} $x_{\gamma_k,\gamma_i}=
x_{\gamma_k,\gamma_{k-1}}x_{\gamma_{k-1},\gamma_i}$. Hence for $i<k$
\[
s_iX_{\tau_i}\cdot X_{\gamma_k}+x_{\gamma_k,\gamma_i}\rho_{i-1}=x_{\gamma_k,\gamma_{k-1}}
\big(s_iX_{\tau_i}\cdot X_{\gamma_{k-1}}+x_{\gamma_{k-1},\gamma_i}\rho_{i-1}\big). 
\] 
Therefore we get by Equation (\ref{psi1})
\[
\begin{array}{lll}
\psi(k)&=&\displaystyle
\sum_{i=0}^{k-1}\big(s_iX_{\tau_i}\cdot X_{\gamma_k}+x_{\gamma_k,\gamma_i}\rho_{i-1}\big)
+s_kX_{\tau_k}\cdot X_{\gamma_k}+x_{\gamma_k,\gamma_k}\rho_{k-1}\\
&=&\displaystyle
x_{\gamma_k,\gamma_{k-1}}
\sum_{i=0}^{k-1}\big(s_iX_{\tau_i}\cdot X_{\gamma_{k-1}}+x_{\gamma_{k-1},\gamma_i}\rho_{i-1}\big)
+s_kX_{\tau_k}\cdot X_{\gamma_k}+\rho_{k-1}.
\end{array}
\]
But a look at Equation (\ref{psi1}) again shows that this is the same as
\[
\psi(k)=x_{\gamma_k,\gamma_{k-1}}\psi(k-1)+s_kX_{\tau_k}\cdot X_{\gamma_k}+\rho_{k-1}. 
\] 
Proposition \ref{tau} says that $x_{\tau_k,\gamma_{k-1}}=\cdots=x_{\tau_k,\tau_k}=1$, which implies that 
\[
\begin{array}{lll}
X_{\tau_k}\cdot X_{\gamma_k}
&=&X_{\tau_k}^{\le\gamma_{k-1}}\cdot X_{\gamma_k}^{\le\gamma_{k-1}}+X_{\tau_k}\cdot X_{\gamma_k}^{>\gamma_{k-1}}\\
&=&x_{\tau_{k},\gamma_{k-1}}x_{\gamma_{k},\gamma_{k-1}}X_{\gamma_{k-1}}^2+\rho_{k-1}\\
&=&x_{\gamma_{k},\gamma_{k-1}}X_{\gamma_{k-1}}^2+\rho_{k-1},
\end{array}
\] 
as $X_{\gamma_k}^{\le\gamma_{k-1}}=x_{\gamma_{k},\gamma_{k-1}}X_{\gamma_{k-1}}$ and
$X_{\tau_k}^{\le\gamma_{k-1}}=x_{\tau_{k},\gamma_{k-1}}X_{\gamma_{k-1}}$ by Proposition \ref{tr}. 
This yields a recursion formula 
\begin{equation}\label{psi2}
\psi(k)=x_{\gamma_k,\gamma_{k-1}}\psi(k-1)+s_k (x_{\gamma_{k},\gamma_{k-1}}X_{\gamma_{k-1}}^2+\rho_{k-1})+\rho_{k-1}.
\end{equation}
We claim that for $k=0,\dots,\nu$
\begin{equation*}\label{psi3}
\psi(k)=t_{k+1}+1+X_{\gamma_k}^2.
\end{equation*}
Taking $k=\nu$, Equation (\ref{psi4}) results from the recursion formula above.

We use induction on $k$. 
By definition $s_0=t_1$ and $\gamma_0=\tau_0=1$, so that Equation (\ref{psi1}) gives 
\[
\psi(0)=s_0+2=t_1+1+X_{\gamma_0}^2.
\]  
Assume next that $\psi(k-1) = t_k+1 + X_{\gamma_{k-1}}^2$ for some $1\le k\le\nu$. Then the recursion in Equation 
(\ref{psi2}) yields 
\begin{equation*}
\begin{array}{crl}
\psi(k)&=&(t_k+1 + X_{\gamma_{k-1}}^2)x_{\gamma_k,\gamma_{k-1}} 
           + s_k(x_{\gamma_k,\gamma_{k-1}}X_{\gamma_{k-1}}^2+ \rho_{k-1})+\rho_{k-1}\\
       &=&(t_k+1)x_{\gamma_k,\gamma_{k-1}} + (s_k+1)(x_{\gamma_k,\gamma_{k-1}}X_{\gamma_{k-1}}^2+ \rho_{k-1}).
\end{array}
\end{equation*}
Subsequently, by Lemma \ref{syt} (i) we get $\psi(k)=t_{k+1}+1+X_{\gamma_k}^2$ as needed.
This completes the step B).

\medskip
C) 
It remains to show that $R[X_\gamma]=c_R$.
If this is not the case, then by Proposition \ref{c=RH} 
there is an integer $k\neq\nu+1$ satisfying
\begin{equation}\label{count}
c_R=R[X_{\gamma_k}]<R[X_\gamma]
\end{equation}
($\gamma=\gamma_{\nu+1}$).
Assume first that $k<\nu+1$. Then $\tau_k<\gamma_k<\tau_{k+1}$ by Equation (\ref{tgt}), which yields 
$\tilde R_k=(R^{<\gamma_k})^{\sim}$, and subsequently, 
\[
R^{<\gamma_k}[X_{\gamma_k}]=\frac{\tilde R_k\cdot X_{\gamma_k}+\Sigma X_{\gamma_k}+1}{I\cdot X_{\gamma_k}}
=\frac{\psi(k)}{I\cdot X_{\gamma_k}}. 
\] 
By Lemma \ref{R=R'} this gives
\[
\frac{\psi(k)}{I\cdot X_{\gamma_k}}
\le\frac{1}{a_{\gamma_k}}.
\]
On the other hand, using Proposition \ref{tr} we get 
\[
I\cdot X_{\gamma_k}=I\cdot X_{\gamma_k}^{\le\gamma_k}=I^{\le\gamma_k}\cdot X_{\gamma_k}=
a_{\gamma_k}X_{\gamma_k}^2,
\] 
and then by Equation (\ref{psi3}) 
\[
\frac{\psi(k)}{I\cdot X_{\gamma_k}}=\frac{t_{k+1}+1+X_{\gamma_k}^2}{a_{\gamma_k}X_{\gamma_k}^2} 
>\frac{1}{a_{\gamma_k}},
\]
which is a contradiction. Therefore we must have $k>\nu+1$. Then $\nu<g^*$, and by assumption and the step B) 
\begin{equation}\label{count1}
R[X_\gamma]=
c\le\frac{m+1}{a_\gamma}.
\end{equation}
Clearly, we may rewrite inequality (\ref{count}) 
\[
\begin{array}{rl}
&\displaystyle\frac{\tilde R\cdot X_{\gamma_k}+\Sigma X_{\gamma_k}+1}{I\cdot X_{\gamma_k}}\vspace{5pt}\\
=&\displaystyle
\frac{\tilde R\cdot X_{\gamma_k}^{\le\gamma}+\Sigma X_{\gamma_k}^{\le\gamma}+x_{\gamma_k,\gamma}+
\tilde R\cdot X_{\gamma_k}^{>\gamma}+\Sigma X_{\gamma_k}^{>\gamma}+1-x_{\gamma_k,\gamma}}
{I\cdot X_{\gamma_k}^{\le\gamma}+I\cdot X_{\gamma_k}^{>\gamma}}\vspace{5pt}\\
=&\displaystyle
\frac{x_{\gamma_k,\gamma}(\tilde R\cdot X_\gamma+\Sigma X_\gamma+1)+
\tilde R\cdot X_{\gamma_k}^{>\gamma}+\Sigma X_{\gamma_k}^{>\gamma}+1-x_{\gamma_k,\gamma}}
{x_{\gamma_k,\gamma}I\cdot X_\gamma+I\cdot X_{\gamma_k}^{>\gamma}}\vspace{5pt}\\
<&\displaystyle\frac{\tilde R\cdot X_\gamma+\Sigma X_\gamma+1}{I\cdot X_\gamma},
\vspace{5pt}
\end{array}
\]
where the second equality follows from Proposition \ref{tr}. Applying Equation (\ref{abc2}) to this, we get
\begin{equation}\label{cdc2}
\frac{\tilde R\cdot X_{\gamma_k}^{>\gamma}+\Sigma X_{\gamma_k}^{>\gamma}+1-x_{\gamma_k,\gamma}}
{I\cdot X_{\gamma_k}^{>\gamma}} <R[X_\gamma]. 
\end{equation}
We aim to prove that 
\[
\frac{\tilde R\cdot X_{\gamma_k}^{>\gamma}+\Sigma X_{\gamma_k}^{>\gamma}+1-x_{\gamma_k,\gamma}}
{I\cdot X_{\gamma_k}^{>\gamma}}=\Phi_k, 
\]
where $\Phi_k$ is as given in Lemma \ref{syt}. This will lead to a contradiction proving the claim, because 
$\Phi_k\ge(m+1)/a_\gamma$ by Lemma \ref{syt} (iii), and then
\begin{equation}\label{cdc4}
R[X_\gamma]\le\frac{m+1}{a_\gamma}\le\Phi_k<R[X_\gamma].
\end{equation}

Because
\[
\tilde R=\sum_{i=0}^{g+1}s_iX_{\tau_i}
\text{ and }
X_{\gamma_k}^{>\gamma}=\sum_{j=\nu+1}^{k-1}X_{\gamma_k}^{(\gamma_j,\gamma_{j+1}]},
\]
we get 
\[
\tilde R\cdot X_{\gamma_k}^{>\gamma}
=\sum_{i=0}^{g+1}\sum_{j=\nu+1}^{k-1}s_iX_{\tau_i}\cdot X_{\gamma_k}^{(\gamma_j,\gamma_{j+1}]}
=\sum_{i=0}^{g+1}\sum_{j=\nu+1}^{k-1}s_i[X_{\tau_i}\cdot X_{\gamma_k}]_j.
\]
Here $[X_{\tau_i}\cdot X_{\gamma_k}]_j=0$ whenever $i\le j<k$, because then $\tau_i<\gamma_j$. Therefore 
\[
\tilde R\cdot X_{\gamma_k}^{>\gamma}
=\sum_{j=\nu+1}^{k-1}\sum_{i=j+1}^{g+1}s_i[X_{\tau_i}\cdot X_{\gamma_k}]_j.
\]
As observed in Remark \ref{keyU}, $\gamma_k\in U$ and $\tau_i\notin U$ when $0<i<g+1$. Note also that 
$\gamma_k\le\tau_{g+1}=n$. Therefore $[X_{\tau_i}\cdot X_{\gamma_k}]_j=x_{\tau_i,\gamma_j}\rho_{\gamma_k,\gamma_j}$ 
for $0<i\le g+1$ by Proposition \ref{key1}. 
If $j<i$, then $\gamma_j\le\gamma_{i-1}\le\tau_i$, and by Corollary \ref{ax} we get 
$x_{\tau_i,\gamma_j}=x_{\tau_i,\gamma_{i-1}}x_{\gamma_{i-1},\gamma_j}$, where $x_{\tau_i,\gamma_{i-1}}=1$ according to 
Proposition \ref{tau}. Hence
\begin{equation}\label{cdc3}
[X_{\tau_i}\cdot X_{\gamma_k}]_j=x_{\tau_i,\gamma_j}\rho_{\gamma_k,\gamma_j}
=x_{\gamma_{i-1},\gamma_j}\rho_{\gamma_k,\gamma_j}
\end{equation}
 for $\nu<j<k$ and $j<i<g$. Subsequently, by Lemma \ref{syt} (ii)
\[
\tilde R\cdot X_{\gamma_k}^{>\gamma}
=\sum_{j=\nu+1}^{k-1}\left(\sum_{i=j+1}^{g+1}s_ix_{\gamma_{i-1},\gamma_j}\right)\rho_{\gamma_k,\gamma_j}
=\sum_{j=\nu+1}^{k-1}m_j\rho_{\gamma_k,\gamma_j}.
\]
By Proposition \ref{key0} we know that 
\[
\Sigma X_{\gamma_k}^{>\gamma}+1-x_{\gamma_k,\gamma}=
\sum_{j=\nu+1}^{k-1}\rho_{\gamma_k,\gamma_j}.
\]
Thus we get
\[
\tilde R\cdot X_{\gamma_k}^{>\gamma}+\Sigma X_{\gamma_k}^{>\gamma}+1-x_{\gamma_k,\gamma}=
\sum_{j=\nu+1}^{k-1}(m_j+1)\rho_{\gamma_k,\gamma_j}.
\]
Moreover, $[I\cdot X_{\gamma_k}]_j=a_{\gamma_j}\rho_{\gamma_k,\gamma_j}$ for every 
$j<k$ by Proposition \ref{key1}, and so
\begin{equation*}
I\cdot X_{\gamma_k}^{>\gamma}=
\sum_{j=\nu+1}^{k-1}I\cdot X_{\gamma_k}^{(\gamma_j,\gamma_{j+1}]}=
\sum_{j=\nu+1}^{k-1}[I\cdot X_{\gamma_k}]_j=\sum_{j=\nu+1}^{k-1}a_{\gamma_j}\rho_{\gamma_k,\gamma_j}.
\end{equation*}
Hence we finally obtain 
\[
\frac{\tilde R\cdot X_{\gamma_k}^{>\gamma}+\Sigma X_{\gamma_k}^{>\gamma}+1-x_{\gamma_k,\gamma}}
{I\cdot X_{\gamma_k}^{>\gamma}}=
\frac{\sum_{j=\nu+1}^{k-1}(m_j+1)\rho_{\gamma_k,\gamma_j}}
     {\sum_{j=\nu+1}^{k-1}a_{\gamma_j}\rho_{\gamma_k,\gamma_j}}=\Phi_k,
\]
which leads to the contradiction (\ref{cdc4}). 
The proof is thus complete.
\end{proof}

As a corollary, we give the result for the monomial case.

\begin{cor} \label{A}
Let $\fa$ be a simple complete $\fm_\ga$-primary ideal in a two-dimensional regular local ring $\ga$. Let 
$I$ be a point basis of $\fa$. Suppose that $\Gamma^*=\emptyset$. Write 
$a:=ord(\fa)$ and $b:=I^2/a$. Then the set $\cH_{\fa}$ of the jumping numbers of the ideal $\fa$ is
\[
H:=\left\{\frac{s+1}{a}+\frac{t+1}{b}\hspace{2pt}\middle|\hspace{2pt}s,t\in\mathbb N\right\}.
\]
\end{cor}

\begin{lem}\label{hö}
In the setting of Theorem \ref{main}, we have for every $\nu\in\{0,\dots,g^*\}$
\[
\frac{1}{a_{\gamma_\nu}}+\frac{1}{b_\nu}\in H_\nu.
\]
In particular, the subsets $H_0,\dots,H_{g^*}$ are non empty.
\end{lem}

\begin{proof}
Obviously, the claim holds for $\nu=g^*$. Suppose that $0\le\nu<g^*$. Note that then $\nu<g$, which implies that $n>1$, and 
further, $b_\nu>a_{\gamma_\nu}$ by Lemma \ref{1ain}. Recall also that $a_{\gamma_\nu}>a_{\gamma_{\nu+1}}$ by 
Proposition \ref{taudit}. Because $a_{\gamma_{\nu+1}}\mid a_{\gamma_\nu}$ by Proposition \ref{aiaj1}, it follows that 
$b_\nu>a_{\gamma_\nu}\ge2a_{\gamma_{\nu+1}}$. Therefore
\[
\frac{1}{a_{\gamma_\nu}}+\frac{1}{b_\nu}\le\frac{1}{a_{\gamma_{\nu+1}}}, 
\]
and then the claim follows from Theorem \ref{main}. Hence $H_\nu$ is a non empty, in fact an infinite set for every 
$\nu\in\{0,\dots,g^*\}$. 
\end{proof}

\begin{ntt}\label{höntt}
For every $\nu\in\{0,\dots,g^*\}$, write
\[
\xi'_\nu:=\frac{1}{a_{\gamma_\nu}}+\frac{1}{b_\nu}\,\;\;\big(=\min H_\nu\,\big).
\]
\end{ntt}

\begin{prop}\label{1ei}
Let $\fa$ be a simple complete ideal in a two-dimensional regular local ring $\ga$ with the point basis 
$(a_1,\dots,a_n)$. Then $a_{\gamma_i}^{-1}\notin\cH_{\fa}$ for any $i\in\{0,\dots,g^*+1\}$. In particular, 
$1\notin\cH_\fa$.
\end{prop}

\begin{proof}
Let $\nu\in\{0,\dots,g^*\}$ and let us write $\gamma:=\gamma_{\nu+1}$ and $\eta:=\gamma_\nu$. According to Corollary 
\ref{ax} we have $a_\eta=a_\gamma x_{\gamma,\eta}$. Furthermore, by Lemma \ref{1ain} we have 
$b_\nu=a_\gamma(x_{\gamma,\eta}X_\eta^2+\rho_\nu)$ and $\gcd\{a_\eta,b_\nu\}=a_\gamma$. 
Let us write $a:=a_\eta/a_\gamma$ and $b:=b_\nu/a_\gamma$. Then $a$ and $b$ are positive integers with $\gcd\{a,b\}=1$.

Assume that $1/a_\gamma\in H_\nu$. Theorem \ref{main} then shows that there are positive 
integers $u$ and $v$ with 
\[
\frac{v}{a_\eta}+\frac{u}{b_\nu}=\frac{1}{a_\gamma}. 
\]
So $vb=(b-u)a$. This is impossible, because $\gcd\{a,b\}=1$. Therefore $1/a_\gamma\notin H_\nu$, which further implies 
that $m/a_\gamma\notin H_\nu$ for any $m\in\mathbb N$ and $\nu<g^*$. Especially, if $i>\nu$, then we may choose 
$m=x_{\gamma_i,\gamma}$, and we see by Corollary \ref{ax} that $1/a_{\gamma_i}\notin H_\nu$. 
If $i\le\nu$, then $a_{\gamma_i}\ge a_\eta$, which shows that
\[
\frac{1}{a_{\gamma_i}}<\frac{1}{a_\eta}+\frac{1}{b_\nu}=\min H_\nu.
\]
Therefore $a_{\gamma_i}^{-1}\notin H_\nu$ for any $i$ and for any $\nu$, i.e., $a_{\gamma_i}^{-1}\notin\cH_\fa$. 
\end{proof}

\begin{rem}\label{1ei2}
By Proposition \ref{1ei} we observe that we could define the sets $H_\nu$ for $\nu=0,\dots,g^*-1$ in Theorem \ref{main} 
as follows:
\[
H_\nu:=\left\lbrace \frac{s+1}{a_{\gamma_\nu}}+\frac{t+1}{b_\nu}+\frac{m}{a_{\gamma_{\nu+1}}}
\hspace{2pt}\middle |\hspace{2pt} s,t,m\in\mathbb N, 
\frac{s+1}{a_{\gamma_\nu}}+\frac{t+1}{b_\nu}<\frac{1}{a_{\gamma_{\nu+1}}}\right\rbrace.
\]
\end{rem}

\begin{prop}\label{höm}
Let $\fa$ be a simple complete ideal in a two-dimensional regular local ring $\ga$ with the point basis 
$I=(a_1,\dots,a_n)$. Then 
\[
\left\lbrace c\in H_{g^*}\hspace{2pt}\middle|\hspace{2pt} c>1\right\rbrace=
\left\lbrace1+\frac{k+1}{I^2}\hspace{2pt}\middle|\hspace{2pt} k\in\mathbb N\right\rbrace.
\]
Especially, every integer greater than one is a jumping number. Moreover, $1+1/I^2$ is the smallest jumping number 
at least one, while $1-1/(I^{\le\gamma_g})^2$ is the greatest jumping number at most one, whenever $g>0$. 
\end{prop}

\begin{proof}
Let the sets $H_0,\dots,H_{g^*}$ be as in Theorem \ref{main}. For $0\le\nu\le g^*$, let us write 
$\eta:=\gamma_\nu$ and $\gamma:=\gamma_{\nu+1}$. Set $a:=a_\eta/a_\gamma$ and $b:=b_\nu/a_\gamma$.
Proposition \ref{tr} yields $I^{\le\gamma}=a_\gamma X_\gamma$, and thus $ab=X_\gamma^2$.
By Lemma \ref{1ain} $\gcd\{a_\eta,b_\nu\}=a_\gamma$, i.e., $\gcd\{a,b\}=1$. It 
follows now from Lemma \ref{syt2} that we can find $s,t\in\mathbb N$ for any $k\in\mathbb N$ satisfying
\begin{equation}\label{iik}
ab+k+1=(t+1)a+(s+1)b.
\end{equation}
Recall also that by  Theorem \ref{main} every element in $H_\nu$ is of the form
\begin{equation}\label{iikkis}
\frac{(t+1)a+(s+1)b+mab}{a_\gamma ab} 
\end{equation}
for some $s,t,m\in\mathbb N$. Thereby we observe that 
\begin{equation}\label{iikki}
H_\nu\subset\left\{\frac{k+1}{(I^{\le\gamma})^2}\hspace{2pt}\middle|\hspace{2pt}k\in\mathbb N\right\}
\end{equation}
for every $\nu\in\{0,\dots,g^*\}$.

Choosing $\nu=g^*$ gives $X_\gamma=X_{\gamma_{g^*+1}}=I$ and $a_\gamma=a_{\gamma_{g^*+1}}=1$. Moreover, then
$a=a_{\gamma_{g^*}}$ and $b=b_{\gamma_{g^*}}$.
Equations (\ref{iik}) and (\ref{iikki}) now yield that 
\[
\left\lbrace1+\frac{k+1}{I^2}\hspace{2pt}\middle|\hspace{2pt} k\in\mathbb N\right\rbrace=
\left\lbrace c\in H_{g^*}\hspace{2pt}\middle|\hspace{2pt} c>1\right\rbrace.
\] 
Subsequently, $1+1/I^2$ is the least element in $H_{g^*}$ greater than one, and every integer greater than one is in 
$H_{g^*}$. Furthermore, Equation (\ref{iikki}) implies that $1+1/I^2$ must be the smallest jumping number at least one, 
as $1\notin\cH_\fa$ by Proposition \ref{1ei},

For the last claim, note first that if $g=0$, then $g^*=0$. In that case Theorem \ref{main} implies that 
$\cH_\fa=H_{g^*}$, and $a_{\gamma_{g^*+1}}=a_{\gamma_{g^*}}=1$ by Proposition \ref{taudit}. Especially, Equation 
(\ref{iikki}) then shows that every jumping number of $\fa$ is greater than one.

Suppose next that $g>0$. Observe that if $c\in H_\nu$ and $c\le1$, then it follows from Equation (\ref{iikkis}) that 
$\nu<g$, because $a_{\gamma_g}=a_{\gamma_{g+1}}=1$ by Proposition \ref{taudit}. Set $\nu=g-1$.Then we get 
$a>a_\gamma=1$ by Proposition \ref{taudit}, and By Lemma \ref{1ain} $a\le b$. Let $k=a-2$ in Equation (\ref{iik}), and 
take $s,t\in\mathbb N$ accordingly so that $ab+a-1=(s+1)b+(t+1)a$. Note that $t$ must be positive as otherwise 
$b\mid ab-1$ implying $b=1$. But then $1<a\le b=1$, which is impossible. Thus
\[
(I^{\le\gamma_g})^2-1=ab-1=(s+1)b+(t'+1)a,
\] 
where $t'=t-1\in\mathbb N$. In particular, Theorem \ref{main} now shows that
\[
1-\frac{1}{(I^{\le\gamma_g})^2}=\frac{s+1}{a}+\frac{t'+1}{b}\in H_{g-1}. 
\]
Moreover, it follows now from Equation (\ref{iikki}) that this must be the maximal jumping number at most one.
\end{proof}

\begin{cor} \label{HS} 
Let $\fa$ be a simple complete ideal in a two-dimensional regular local ring $\ga$ and let $\cH_\fa$ denote the set 
of the jumping numbers of $\fa$. The Hilbert-Samuel multiplicity of $\fa$ is 
\[
e(\fa)=(\xi'-1)^{-1},
\] 
where $\xi':=\min\{\xi\in\cH_\fa\mid\xi>1\}$. 
\end{cor}

\begin{proof}
It follows from the Hoskin-Deligne formula that the Hilbert-Samuel multiplicity of $\fa$ is $I^2$, where $I$ denotes 
the point basis vector of $\fa$ (see \cite[Corollary 3.8]{L1}). By Proposition \ref{höm} we now know that 
$\xi'=(1+I^2)^{-1}$.
\end{proof}

\begin{lem}\label{höllä}
For $\mu,\nu\in\{0,\dots,g^*\}$ set 
\[
\theta_{\mu,\nu}:=\min\left\{\xi\in H_\mu\mid\xi\ge\frac{1}{a_{\gamma_\nu}}\right\}.
\]
Then 
\[
\theta_{\mu,\nu}=
\left\{
\begin{array}{ll}
\displaystyle
\xi'_\mu+\frac{1}{a_{\gamma_\nu}}&\text{ if }\mu<\nu;\\
\xi'_\mu&\text{ if }\mu\ge\nu.
\end{array}
\right.
\]
In particular, $\theta_{\mu,\nu}\ge\xi'_\nu$, where the equality holds if and only if $\mu=\nu$. 
\end{lem}

\begin{proof}
Since $1/a_{\gamma_\nu}<\xi'_\nu=\min H_\nu$, we observe that $\theta_{\nu,\nu}=\xi'_\nu$.

Suppose that $\mu>\nu$. Then $a_{\gamma_{\nu+1}}\ge a_{\gamma_\mu}$ by Proposition \ref{taudit}. As $\nu<g^*$ and 
$\xi'_\nu\in H_\nu$ by Lemma \ref{hö}, 
Theorem \ref{main} shows that $\xi'_\nu\le1/a_{\gamma_{\nu+1}}$. Hence 
\[
\xi'_\nu=\frac{1}{a_{\gamma_\nu}}+\frac{1}{b_\nu}\le\frac{1}{a_{\gamma_{\nu+1}}}\le\frac{1}{a_{\gamma_\mu}}
<\frac{1}{a_{\gamma_\mu}}+\frac{1}{b_\mu}=\xi'_\mu,
\]
which yields $\xi'_\nu<\xi'_\mu=\theta_{\mu,\nu}$.

Suppose next that $\mu<\nu$. In particular, this gives $\mu<g^*$. As $\theta_{\mu,\nu}\in H_\mu$, we then know by 
Theorem \ref{main} that $\theta_{\mu,\nu}=c+m/a_{\gamma_{\mu+1}}$ for some $m\in\mathbb N$ and $c\in H_\mu$ 
with $c\le1/a_{\gamma_{\mu+1}}$. By Proposition \ref{1ei} we know that $\theta_{\mu,\nu}>1/a_{\gamma_\nu}$, and since 
$1/a_{\gamma_\nu}=x_{\gamma_\nu,\gamma_{\mu+1}}/a_{\gamma_{\mu+1}}$ by Corollary \ref{ax}, 
we observe that $m=x_{\gamma_\nu,\gamma_{\mu+1}}$ and $c=\xi'_\mu$. Thereby 
$\theta_{\mu,\nu}=\xi'_\mu+1/a_{\gamma_\nu}$.

Because $a_{\gamma_\mu}b_\mu=(I^{\le\gamma_{\mu+1}})^2\le(I^{\le\gamma_{\nu+1}})^2=a_{\gamma_\nu}b_\nu$, we have 
$b_\mu\le b_\nu$. Thus
\[
\xi'_\mu=\frac{1}{a_{\gamma_\mu}}+\frac{1}{b_\mu}>\frac{1}{b_\nu}=\xi'_\nu-\frac{1}{a_{\gamma_\nu}},
\]
which shows that $\theta_{\mu,\nu}>\xi'_\nu$, as wanted. 
\end{proof}

\begin{prop}\label{höl}
In the setting of Theorem \ref{main} we get for $\nu\in\{0,\dots,g^*\}$
\[
\xi'_\nu=\min\left\{\xi\in\cH_\fa\hspace{2pt}\middle|\hspace{2pt}\xi\ge\frac{1}{a_{\gamma_\nu}}\right\}.
\]
Moreover, $\xi'_\nu\in H_\mu$ if and only if $\mu=\nu$, and 
\[
\frac{1}{a_{\gamma_0}}<\xi'_0<\cdots<\frac{1}{a_{\gamma_{g^*}}}<\xi'_{g^*}.
\] 
\end{prop}

\begin{proof}
Note first that $1/a_{\gamma_\nu}<\xi'_\nu$ for every $\nu=0,\dots,g^*$. Furthermore, $\xi'_\nu\in H_\nu$ by Lemma 
\ref{hö}. By Remark \ref{1ei2} we then see that $\xi'_\nu<1/a_{\gamma_{\nu+1}}$ for $\nu<g^*$. 
If $\xi'_\nu\in H_\mu$ for some $\mu,\nu\in\{0,\dots,g^*\}$, then $\theta_{\mu,\nu}\le\xi'_\nu$ by definition as 
$\xi'_\nu>1/a_{\gamma_\nu}$. It follows from Lemma \ref{höllä} that $\theta_{\mu,\nu}=\xi'_\nu$, and further, 
$\mu=\nu$. Moreover, $\xi'_\nu=\min\left\{\theta_{\mu,\nu}\mid\mu=0,\dots,g^*\right\}$ for any $\nu=0,\dots,g^*$. Hence
\[
\xi'_\nu=\min\left\{\xi\in\cH_\fa\hspace{2pt}\middle|\hspace{2pt}\xi\ge\frac{1}{a_{\gamma_\nu}}\right\}.
\]
\end{proof}

\begin{cor} \label{maincor} 
Let $\fa$ be a simple complete ideal of finite colength in a two-dimensional regular local ring $\ga$. The sequence of 
pairs $(a_{\gamma_0},b_0),\dots,(a_{\gamma_{g^*}}b_{g^*})$ and thereby the set of the jumping numbers of $\fa$, is 
totally determined by the numbers $\xi'_0,\dots,\xi'_{g^*}$.
\end{cor}

\begin{proof} For  $\nu\in\{0,\dots,g^*\}$, write $\gamma:=\gamma_{\nu+1}$, $\eta:=\gamma_\nu$, 
$u_\nu:=a_{\eta}/a_{\gamma}$ and $v_\nu:=b_\nu/a_{\gamma}$. Corollary \ref{ax} gives $u_\nu=x_{\gamma,\eta}$ while 
$v_\nu=x_{\gamma,\eta}X^2_{\eta}+\rho_\nu$  by Lemma \ref{1ain}. Then
\[
\xi'_\nu
=\frac{1}{a_\eta}+\frac{1}{b_\nu}
=\frac{a_{\eta}+b_\nu}{a_{\eta}b_\nu}
=\frac{1}{a_\gamma}\cdot
\frac{u_\nu+v_\nu}{u_\nu v_\nu}.
\]
Moreover, by Lemma \ref{1ain} we have
\[
\gcd\{u_\nu+v_\nu,u_\nu v_\nu\}=\gcd\{u_\nu,v_\nu\}=\frac{\gcd\{a_\eta,b_\nu\}}{a_\gamma}=1.
\]
Thus $a_\gamma$ and $\xi'_\nu$ determine $u_\nu+v_\nu$ and $u_\nu v_\nu$. 
Note that $u_\nu$ and $v_\nu$ are the 
roots of the quadratic equation $\omega^2-(u_\nu+v_\nu)\omega+u_\nu v_\nu=0$. So $u_\nu$ and $v_\nu$ are uniquely 
determined by $a_\gamma$ and $\xi'_\nu$.

Given all the numbers $\xi'_0,\dots,\xi'_{g^*}$, suppose that we would know 
the integer $a_\gamma$ as well as the pair $(u_\nu,v_\nu)$ for some $\nu>0$. Then we obtain $a_\eta=a_\gamma u_\nu$, 
and from the product $a_\eta\xi'_{\nu-1}$ we get the pair $(u_{\nu-1},v_{\nu-1})$ as described above. Because 
$a_{\gamma_{g^*+1}}=1$, we see that $\xi_{g^*}$ yields $(u_{g^*},v_{g^*})=(a_{\gamma_{g^*}},b_{g^*})$ so that we 
eventually get all the pairs $(u_0,v_0),\dots,(u_{g^*},v_{g^*})$ as well as the integers 
$a_{\gamma_0},\dots,a_{\gamma_{g^*}}$. Subsequently, we get the sequence 
$(a_{\gamma_0},b_0),\dots,(a_{\gamma_{g^*}},b_{g^*})$, and the claim now follows from Theorem \ref{main}.
\end{proof}

\begin{lem}\label{211}
In the setting of Theorem \ref{main}, consider the set 
\[
H':=\left\{\xi\in\cH_\fa\mid\xi\le\frac{1}{a_{\gamma_1}}\right\}.
\]
Write $a:=a_{\gamma_0}/a_{\gamma_1}$ and $b:=b_0/a_{\gamma_1}$. Then $H'=\{\xi\in H_0\mid\xi<1/a_{\gamma_1}\}$, and
\begin{itemize}
\item[(i)]  $H'$ is empty, if and only if $g=0;$ 
\item[(ii)] $H'$ has exactly one element, if and only if $(a,b)=(2,3);$
\item[(ii)] $H'$ has exactly two elements, if and only if $(a,b)=(2,5)$.
\end{itemize}
\end{lem}

\begin{proof} 
Observe that Lemma \ref{1ain} yields $\gcd\{a_{\gamma_0},b_0\}=a_{\gamma_1}$ and $a_{\gamma_0}\le b_0$. Thus $a\le b$ 
are positive integers. Note also that $\xi'_\nu>1/a_{\gamma_1}$ whenever $\nu\ge1$ by Proposition \ref{höl}. Thereby 
$H'=\{\xi\in H_0\mid\xi<1/a_{\gamma_1}\}$ where the inequality is strict by Remark \ref{1ei2}. 

(i) $H'=\emptyset$ exactly when $a_{\gamma_1}\xi'_0=1/a+1/b>1$. This is the case if and only if $a=1$, i.e., 
$a_{\gamma_0}=a_{\gamma_1}$. By Proposition \ref{taudit} this is equivalent to $g=0$.

(ii) Theorem \ref{main} implies that $H'$ has exactly one element, if and only if $\xi'_0\in H'$, while 
$1/a_{\gamma_0}+2/b_0\notin H'$. This happens exactly when
\[
\frac{1}{a}+\frac{1}{b}<1\text{ and }\frac{1}{a}+\frac{2}{b}>1.
\]
Clearly, this takes place if and only if $(a,b)=(2,3)$.

(iii) It follows from Theorem \ref{main} that $H'$ has exactly two elements, if and only if 
both $\xi'_0$ and $1/a_{\gamma_0}+2/b_0$ are in $H'$ while $1/a_{\gamma_0}+3/b_0$ and $2/a_{\gamma_0}+1/b_0$ are not. 
This is equivalent to the condition
\[
\frac{1}{a}+\frac{2}{b}<1<\min\left\{\frac{1}{a}+\frac{3}{b},\frac{2}{a}+\frac{1}{b}\right\}.
\]
Obviously $a>1$. If $a\ge3$ then $2/a+1/b\le1$ for any $b\ge a$. Therefore $a=2$, and the inequality on the left implies 
that $b>4$. If $b\ge6$ then $1/2+2/b<1$. Thus the only possible solution to this is $(a,b)=(2,5)$.  
\end{proof}

\begin{thm}\label{1maincor} 
The point basis $I=(a_1,\dots,a_n)$ of a simple complete ideal $\fa$ of finite colength in a two-dimensional regular 
local ring $\ga$ can be read off from the set $\cH_\fa$ of its jumping numbers. In particular, the jumping numbers 
belonging to the subset $\cH':=\{c\in\cH_\fa\mid 0<c<1\}$ determine the multiplicities $a_1,\dots,a_{\gamma_g}$.
\end{thm}

\begin{proof} 
As $a_{\gamma_1}\ge1$, Lemma \ref{211} (i) implies that $\cH'=\emptyset$, if and only if $g=0$. Then 
$a_1=\cdots=a_n=1$, i.e., $I^2=n$. Subsequently, Proposition \ref{höm} yields $\xi'_0=1+1/n$ so that $n=1/(\xi'_0-1)$, 
and we are done in this case.
   
Consider then the case $\cH'\neq\emptyset$, i.e., $g\ge1$. Suppose first that besides the jumping numbers of $\fa$, we 
would already know the multiplicities $a_1,\dots,a_{\gamma_g}$. According to Proposition \ref{taudit} $a_i=1$ for every 
$i=\gamma_g,\dots,n$. Thus it remains to determine the number $n$. Let $\xi'$ be as above, and let 
$\xi'':=\max\{c\in\cH_\fa\mid c<1\}$. By Proposition \ref{höm} we know that $\xi'=1+1/I^2$, while 
$\xi''=1-1/(I^{\le\gamma_g})^2$.  Therefore we see that 
\[
n-\gamma_g=I^2-(I^{\le\gamma_g})^2=\frac{1}{\xi'-1}+\frac{1}{\xi''-1}.
\]

Let us then prove that the set $\cH'$ determines the multiplicities $a_1,\dots,a_{\gamma_g}$. According to Proposition 
\ref{aiaj1} it is sufficient to find the rational numbers $\gb'_1,\dots,\gb'_g$. Recall that 
$\gb'_{\nu+1}:=1+\rho_{n,\gamma_\nu}/a_{\gamma_\nu}$ by Remark \ref{rhob}. Since 
\[
a_{\gamma_\nu} b_\nu=(I^{\le \gamma_{\nu+1}})^2=a_{\gamma_0}^2+[I^2]_0+\cdots+[I^2]_\nu
\]
for every $\nu=0,\dots,g-1$, we get by Proposition \ref{key1}
\[
a_{\gamma_\nu} b_\nu=a_{\gamma_0}^2+a_{\gamma_0}\rho_{n,\gamma_0}+\cdots+a_{\gamma_\nu}\rho_{n,\gamma_\nu}.
\]
As $a_{\gamma_\nu} b_\nu-a_{\gamma_{\nu-1}} b_{\nu-1}=a_{\gamma_\nu}\rho_{n,\gamma_\nu}$, it is enough to find out all 
the pairs 
\[
(a_{\gamma_0},b_0),\dots,(a_{\gamma_{g-1}},b_{g-1}).
\]

Arguing inductively, suppose that for some $1\le k\le g-1$  we know the pairs 
$(a_{\gamma_0},b_0),\dots,(a_{\gamma_{k-1}},b_{k-1})$. Then we obtain by Lemma \ref{1ain} 
\[
a_{\gamma_k}=\gcd\{a_{\gamma_{k-1}},b_{k-1}\},
\] 
and by Proposition \ref{taudit} we know that $a_{\gamma_k}>a_{\gamma_g}=1$. The least element in $\cH_\fa$ greater than
$1/a_{\gamma_k}$ is by Proposition \ref{höl} 
\[
\xi'_k=\frac{1}{a_{\gamma_k}}+\frac{1}{b_k},
\]
and because $1<a_{\gamma_k}\le b_k$ by Lemma \ref{1ain} we see that $\xi'_k$ lies in the set $\cH'$.
Now $b_k=(\xi'_k-1/a_{\gamma_k})^{-1}$. Thereby we get also the pair $(a_{\gamma_k},b_k)$.

The problem is to find the first pair of integers $(a_{\gamma_0},b_0)$. By Theorem \ref{ord} below the three 
smallest jumping numbers determine the order of the ideal $\fa$, which is precisely the integer $a_{\gamma_0}$. Then 
$b_0=(\xi'_0-1/a_{\gamma_0})^{-1}$. Thus everything is clear, if $\cH'$ has at least three elements.

It remains to consider the two special cases where $\cH'$ has either two or only one element. Let us show that then 
$g=1$. We already saw that $\cH'\neq\emptyset$ implies $g\ge1$. Suppose that we would have $g>1$. Then it follows from 
Proposition \ref{taudit} that $a_{\gamma_1}>a_{\gamma_g}=1$, i.e., $a_{\gamma_1}\ge2$. Proposition \ref{höl} implies 
that $\xi'_0<1/a_{\gamma_1}$, and further, $\xi'_0<\xi'_1<\xi'_0+1/a_{\gamma_1}<2/a_{\gamma_1}\le 1$. This means that 
$\xi'_0$, $\xi'_1$ and $\xi'_0+1/a_{\gamma_1}$ are all in $\cH'$, which is impossible. Therefore we have $g=1$ in both 
cases so that $a_{\gamma_1}=1$ by Proposition \ref{taudit}. Write $a:=a_{\gamma_0}$ and $b:=b_0$.

Suppose first that $\xi'_0$ is the only element in $\cH'$. By Lemma \ref{211} (ii) we must then have $(a,b)=(2,3)$. 
In this case we get
\[
(a_1,\dots,a_{\gamma_g})=(2,1,1).
\]
Suppose next that $\cH'$ has just two elements. Then by Lemma \ref{211} (iii) we have $(a,b)=(2,5)$, in which case 
\[
(a_1,\dots,a_{\gamma_g})=(2,2,1,1).
\]
Thereby we are done as soon as we have proven the next theorem.
\end{proof}

\begin{thm} \label{ord} 
Let $\fa$ be a simple complete ideal in a two-dimensional regular local ring $\ga$. Let $\xi<\psi<\zeta$ be 
the three smallest jumping numbers of $\fa$, and let $H_0$ be as in Theorem \ref{main}. Then the order of the ideal 
$\fa$ is
\[
\ord(\fa)=
\left\{
\begin{array}{cl}
\displaystyle\frac{5}{3\xi}&\text{ if } 6\xi=10\psi-5\zeta;
\vspace{10pt}\\
\displaystyle\frac{1}{2\xi-\psi}&\text{ if } 6\xi\neq10\psi-5\zeta.
\end{array}
\right.
\]
Moreover, $\ord(\fa)=1$ if and only if $\xi>1$, and if $\xi<1<\zeta$, then $\ord(\fa)=2$. 
\end{thm}

\begin{proof}
Let $a_{\gamma_0},\dots,a_{\gamma_g}$ and $b_0,\dots,b_{g^*}$ be as in Theorem \ref{main}. 
Note that $\ord(\fa)=a_{\gamma_0}$. By Proposition \ref{höl} we know that 
\begin{equation*}\label{q0}
\xi=\xi'_0=\frac{1}{a_{\gamma_0}}+\frac{1}{b_0}\in H_0. 
\end{equation*}
We shall first show that $\psi\in H_0$ implies $\ord(\fa)=1/(2\xi-\psi)$, while $\psi\notin H_0$ gives 
$\ord(\fa)=5/3\xi$. We shall then verify that $6\xi=10\psi-5\zeta$ is equivalent to $\psi\notin H_0$. This will prove 
the first claim.

Suppose first that $\psi\in H_0$. Because $a_{\gamma_0}\le b_0$ by Lemma \ref{1ain}, it follows from 
Theorem \ref{main} that necessarily
\[
\psi=\frac{1}{a_{\gamma_0}}+\frac{2}{b_0},
\]
provided that in the case $g^*\ge1$ we can prove that $\psi\le1/a_{\gamma_1}$. Indeed, suppose that we would have 
$\psi>1/a_{\gamma_1}$. As $g^*\ge1$, we know by Proposition \ref{höl} that $\xi'_1\notin H_0$. This means that 
$\psi<\xi'_1$. Lemma \ref{höllä} yields $\xi'_1<\theta_{0,1}$ so that $1/a_{\gamma_1}<\psi<\theta_{0,1}$ 
contradicting the definition of $\theta_{0,1}$. Thus we observe that if $\psi\in H_0$, then 
\[
a_{\gamma_0}=\frac{1}{2\xi-\psi}.
\]

Suppose next that $\psi\notin H_0$. Then $g^*\ge1$. By Proposition \ref{höl} we get $\psi=\xi'_1>1/a_{\gamma_1}$. So 
$\xi$ is the only jumping number at most $1/a_{\gamma_1}$, in which case Lemma \ref{211} (ii) gives 
$(a_{\gamma_0},b_0)=(2a_{\gamma_1},3a_{\gamma_1})$. Hence 
\[
\xi=\frac{5}{6a_{\gamma_1}},\text{ or equivalently }a_{\gamma_0}=\frac{5}{3\xi},
\]
as wanted.

Let us now verify that in the case $\psi\notin H_0$ we have $6\xi=10\psi-5\zeta$. It is enough to prove that 
\begin{equation}\label{jetta}
\zeta=\frac{1}{a_{\gamma_1}}+\frac{2}{b_1},
\end{equation}
because then $2\psi-\zeta=1/a_{\gamma_1}=6\xi/5$, i.e., $6\xi=10\psi-5\zeta$. Observe that if $g^*=1$ or if $g^*>1$ and 
$1/a_{\gamma_1}+2/b_1<1/a_{\gamma_2}$, then $1/a_{\gamma_1}+2/b_1\in H_1$ by Theorem \ref{main}. Subsequently it is the 
smallest jumpling number in $H_1$ greater than $\psi$, as $a_{\gamma_1}<b_1$.

We aim to show that $1/a_{\gamma_1}+2/b_1<\theta_{0,1}$ and that 
$\theta_{0,1}<1/a_{\gamma_2}$ whenever $g^*>1$. These will then yield that $1/a_{\gamma_1}+2/b_1\in H_1$ and that 
$\zeta\notin H_0$. The condition $\theta_{0,1}<1/a_{\gamma_2}$ will also guarantee by Proposition \ref{höl} that 
$\zeta\notin H_\nu$ for any $\nu>1$ in the case $g^*>1$. But then the only possibility is $\zeta\in H_1$, which will then prove 
Equation (\ref{jetta}).

Since $[I^2]_1=a_{\gamma_1}\rho_{n,\gamma_1}$ by Proposition \ref{key1}, we see that 
\[
b_1=\frac{(I^{\le\gamma_2})^2}{a_{\gamma_1}}=\frac{(I^{\le\gamma_1})^2+[I^2]_1}{a_{\gamma_1}}
=\frac{a_{\gamma_0}b_0+a_{\gamma_1}\rho_{n,\gamma_1}}{a_{\gamma_1}}=6a_{\gamma_1}+\rho_{n,\gamma_1}. 
\]
This, together with Lemma \ref{höllä}, implies that
\begin{equation*}
\frac{1}{a_{\gamma_1}}+\frac{2}{b_1}<\frac{1}{a_{\gamma_1}}+\frac{2}{6a_{\gamma_1}}
<\xi+\frac{1}{a_{\gamma_1}}=\theta_{0,1}. 
\end{equation*}
If $g^*>1$, then it follows from Proposition \ref{taudit} that $a_{\gamma_1}>a_{\gamma_2}$ so that 
\[
\theta_{0,1}=\xi+\frac{1}{a_{\gamma_1}}<\xi+\frac{1}{a_{\gamma_2}}=\theta_{0,2} 
\]
by Lemma \ref{höllä}. This means that $\theta_{0,1}<1/a_{\gamma_2}$. The proof of Equation (\ref{jetta}) is thus 
complete.

Assume then that $6\xi=10\psi-5\zeta$, and that we would have $\psi\in H_0$. As we saw above $2\xi-\psi=1/a_{\gamma_0}$, 
which further yields $\psi-\xi=1/b_0$.

Suppose first that $\zeta<1/a_{\gamma_1}$. It follows from Proposition \ref{höl} that $\zeta\in H_0$.
By applying Theorem \ref{main} we observe that
\[
\zeta=\min\left\{\frac{s+1}{a_{\gamma_0}}+\frac{t+1}{b_0}\hspace{2pt}\middle|\hspace{2pt}
(s,t)\in\mathbb N^2\smallsetminus\{(0,0);(0,1)\}\right\}.
\]
so that 
\[
\zeta=
\frac{1}{a_{\gamma_0}}+\frac{3}{b_0}\hspace{4pt}\text{ or }\hspace{4pt}\zeta=\frac{2}{a_{\gamma_0}}+\frac{1}{b_0}.
\]
If $\zeta=1/a_{\gamma_0}+3/b_0$, then we get
\[
10\psi-6\xi=5\zeta=5(2\xi-\psi)+15(\psi-\xi)=9\xi+10\psi,
\]
which is impossible. If $\zeta=2/a_{\gamma_0}+1/b_0$, then
\[
10\psi-6\xi=5\zeta=10(2\xi-\psi)+5(\psi-\xi)=15\xi-5\psi,
\]
which yields $9/b_0=9\psi-9\xi=12\xi-6\psi=6/a_{\gamma_0}$, i.e., $a_{\gamma_0}/b_0=2/3$. It follows 
from Lemma \ref{211} (ii) that $\psi\notin H_0$, which is a contradiction.

We must thus have $\zeta>1/a_{\gamma_1}$. Note that $1/a_{\gamma_1}$ is not a jumping number by Proposition \ref{1ei}. 
Now $\xi$ and $\psi$ are the only jumping numbers at most 
$1/a_{\gamma_1}$. By Lemma \ref{211} (iii) we then have $a_{\gamma_0}=2a_{\gamma_1}$ and $b_0=5a_{\gamma_1}$. 
This gives
\[
\frac{1}{a_{\gamma_1}}\le\zeta=\frac{10\psi-6\xi}{5}=\frac{4}{5a_{\gamma_0}}+\frac{14}{5b_0}=\frac{24}{25a_{\gamma_1}}<
\frac{1}{a_{\gamma_1}},
\]
which is a contradiction. The first claim has thus been proven.

It remains to prove the last two claims. Proposition \ref{taudit} implies that $a_{\gamma_0}=1$ if and only if 
$g=0$. Moreover, it follows from Lemma \ref{211} (i) that this happens exactly when $\xi>1/a_{\gamma_1}=1$.

Assume then that $\xi<1<\zeta$, in which case $g>0$. This implies that the set $H'$ in Lemma \ref{211} has either one 
or two elements. Indeed, otherwise we would have $\zeta\in H'$ in which case $\zeta\le1/a_{\gamma_1}\le1$. So 
$a_{\gamma_0}=2a_{\gamma_1}$ by Lemma \ref{211} (ii) and (iii). We need to verify that $a_{\gamma_1}=1$. If 
$a_{\gamma_1}$ were greater than one, then we would have $g^*\ge1$ because of Proposition \ref{taudit}. Lemma 
\ref{höllä}  would then yield $\xi'_1<\theta_{0,1}=\xi+1/a_{\gamma_1}$. Since $\xi<1/a_{\gamma_1}<\xi'_1$ by 
Proposition \ref{höl}, we would then have $\xi<\xi'_1<\theta_{0,1}$ so that
\[
1<\zeta\le\theta_{0,1}<\frac{2}{a_{\gamma_1}}\le1,
\] 
which is a contradiction. Therefore $a_{\gamma_1}=1$, i.e., 
$a_{\gamma_0}=2$.
\end{proof}

\begin{exmp} \label{56} 
Suppose that $I=(2,1,1)$. It follows from Theorem \ref{main} that in 
Theorem \ref{ord} $\xi=5/6$, $\psi=7/6$ and $\zeta=8/6$. Then $6\xi=10\psi-5\zeta$.
On the other hand, in the case $I=(3,1,1,1)$ $\xi=7/12$, $\psi=10/12$ and $\zeta=11/12$ so that 
$6\xi\neq10\psi-5\zeta$.
\end{exmp}

\section{Jumping numbers of an analytically irreducible plane curve}

In this section, we aim to utilize Theorem \ref{main} in determining the jumping numbers of an analytically irreducible 
plane curve with an isolated singularity at the origin. As jumping numbers are compatible with localization, it is 
enough to consider the local situation. Therefore in the following we mean by a plane curve 
the subscheme $\cC_f$ of $\Spec \ga$ determined by an element $f$ in the maximal ideal $\fm_\ga$. We will next recall 
some basic facts about plane curves. For more details, we refer to \cite{C}, \cite{KG} and \cite{CA}.

As before, we assume that the residue field $\Bbbk$ of $\ga$ is algebraically closed. If $\gb\supset\ga$, then the 
\textit{strict transform} of $\cC_f$ is $\cC_f^{(\gb)}:=\Spec\gb/(f^{(\gb)})$, where $f^{(\gb)}$ denotes any generator 
of the transform $(f)^{(\gb)}$. The \textit{multiplicity} of $\cC_f$ at $\gb$ is $m_\gb(\cC_f):=\ord_\gb(f^{(\gb)})$. 
Following the terminology of \cite{KG}, we call the set of those $\gb\supset\ga$ for which $f^{(\gb)}\neq\gb$ the 
\textit{point locus} of $f$ (or $\cC_f$).

Suppose that $f$ is analytically irreducible. Then by \cite[Corollary 4.8]{KG} there is a unique quadratic transform 
$\gb\supset\ga$ belonging to the point locus of $f$. By \cite[Corollary 4.8]{KG} $f^{(\gb)}$ is analytically 
irreducible. It follows that the point locus of $f$ consists of a quadratic sequence 
\[
\ga=\ga_1\subset\ga_2\subset\cdots.
\]
This corresponds to a sequence of point blow ups 
\begin{equation}\label{rloc}
\Spec\ga=\mathcal X_1\xleftarrow{\pi_1}\mathcal X_2\xleftarrow{\pi_2}\cdots.
\end{equation}
There exists the smallest $\nu$ such that the total transform $(\pi_\nu\circ\cdots\circ\pi_1)^*\cC_f$ has normal crossing 
support. The morphism $\bar\pi:=\pi_\nu\circ\cdots\circ\pi_1:\mathcal X_{\nu+1}\rightarrow\mathcal X_1$ is called the 
\textit{standard resolution} of $\cC_f$. The \textit{multiplicity sequence} of $\cC_f$ is now 
\[
(m_1,\dots,m_\nu),
\] 
where $m_i=m_{\ga_i}(\cC_f)$.

Let us recall the notion of a \textit{general element} of an ideal. Fix a minimal system $f_1,\dots,f_\mu$ of 
generators for an ideal $\fa$ in $\ga$. Set $\bar\lambda:=\lambda+\fm_\ga\in\Bbbk$ for $\lambda\in\ga$. One says that a 
general element of $\fa$ has some property $\mathcal P$, if there is a non empty open subset $V\subset\Bbbk^\mu$ such 
that $f=\lambda_1f_1+\cdots+\lambda_\mu f_\mu$ has $\mathcal P$ whenever $(\bar\lambda_1,\dots,\bar\lambda_\mu)\in V$. 
Note that the ideal can always be generated by general elements.

Suppose that $\fa$ is a simple complete ideal. It follows from \cite[Corollary 4.10]{KG} that a general element 
$f\in\fa$ is analytically irreducible. As $f$ is general, it easily follows that if $I=(a_1,\dots,a_n)$ is the point 
basis of $\fa$, then $(a_1,\dots,a_\nu)$ is the multiplicity sequence of $\cC_f$. 
It is clear that the resolution (\ref{res}) of $\fa$ contains the 
standard resolution of $\cC_f$, i.e., if $\pi:\mathcal X_{n+1}\rightarrow\mathcal X_1$ is the resolution of $\fa$, then
\[
\pi=\pi_n\circ\cdots\circ\pi_{\nu+1}\circ\bar\pi.
\]

Let $\cC_f^{(i)}$ denote the strict transform of $\cC_f$ on $\mathcal X_i$.
Since $\nu$ is the least integer such that $\bar\pi^*\mathcal C_f$ has only normal crossings, we observe that either 
$\nu=1,2$ or $\cC_f^{\nu-1}$ intersects transversely the exceptional divisor of 
$\pi_{\nu-1}$ and the strict transform of some other exceptional divisor going through the center 
$\varsigma_\nu\in\mathcal X_\nu$ of $\pi_\nu$. Therefore 
$\ga_\nu$ must be a satellite point or $\nu=1,2$. Furthermore, since $\cC_f^{(i+1)}$ 
intersects for every $i\in\{\nu,\dots,n\}$ only one of the exceptional divisors and that transversely, we see that the 
points $\ga_{\nu+1},\dots,\ga_n$ are free. It follows that $\nu=\gamma_g=\max\Gamma_\fa$.

In a lack of a suitable reference we state the following lemma:

\begin{lem}\label{cf}
Let $\fa$ be a simple complete ideal of finite colength in a two-dimensional regular local ring $\ga$ having the 
resolution (\ref{res}) and the point basis $(a_1,\dots,a_n)$. Let $f$ be an analytically irreducible element 
in $\fm_\ga$. The following conditions are then equivalent:
\begin{itemize}
\item[(i)]   $m_{\ga_i}(\cC_f)=a_i$ for $i=1,\dots,n;$
\item[(ii)]  $\pi^*\cC_f=\cC_f^{(n+1)}+\hat E_n;$
\item[(iii)] $\cC_f^{(n+1)}\cdot E_i=\delta_{i,n}$ for $i=1,\dots,n;$
\item[(iv)]  if $E_i^{(n)}$ passes through $\varsigma_n\in \mathcal X_n$ for some $i<n$, then
             $\cC_f^{(n)}$ intersects $E_i^{(n)}$
             transversely at $\varsigma_n$.

\end{itemize}
\end{lem}

\begin{proof}
We know that 
\[
\pi^*\cC_f-\cC_f^{(n+1)}=\sum_{j=1}^nm_{\ga_j}(\cC_f)E_j^*.
\]
Since $\hat E_n=\sum_{j=1}^n a_jE_j^*$, this proves the equivalence of (i) and (ii). Write 
\[
\pi^*\cC_f-\cC_f^{(n+1)}=\sum_{j=1}^n\hat d_j\hat E_j. 
\]
By the projection formula $\pi^*\cC_f\cdot E_i=0$ for all $i=1,\dots,n$. As $\hat E_j\cdot E_i=\delta_{i,j}$ for all 
$j=1,\dots,n$, we now obtain $\hat d_i=\cC_f^{(n+1)}\cdot E_i$. It thus follows that (ii) and (iii) are equivalent. 

In order to prove the equivalence of (iii) and (iv) we first observe that in any case $\varsigma_n\in\cC_f^{(n)}$. 
Because $\cC_f$ is analytically irreducible, we see that if $E_i^{(n)}$ passes through $\varsigma_n$, 
then $\varsigma_n$ is the only point of intersection of $\cC_f^{(n)}$ and $E_i^{(n)}$. This implies that (iv) is 
equivalent to $\cC_f^{(n)}\cdot E_i^{(n)}=1$. As $\pi_n^*E_i^{(n)}=E_i+E_n$, we get by the projection formula 
\[
\cC_f^{(n)}\cdot E_i^{(n)}=\cC_f^{(n+1)}\cdot\pi_n^*E_i^{(n)}
=\cC_f^{(n+1)}\cdot E_i+\cC_f^{(n+1)}\cdot E_n.
\] 
This immediately shows that (iii) implies (iv). Conversely, assuming (iv) and observing 
that the both summands on the right hand side are non negative and $\cC_f^{(n+1)}\cdot E_n\neq 0$, we get 
$\cC_f^{(n+1)}\cdot E_i^{(n+1)}=0$ and $\cC_f^{(n+1)}\cdot E_n=1$. Thereby we see that (iii) holds.
\end{proof}

\begin{rem}\label{SP}
Following \cite[Definition 7.1]{S1} and \cite[Definition 1]{CPR} 
we could define an element $f\in\fa$ to be general, if the corresponding curve $\cC_f$ is analytically irreducible 
and $\cC_f^{(n)}$ intersects transversely at $\varsigma_n\in\mathcal X_n$ every $E_i^{(n)}$ ($i<n$) passing through 
$\varsigma_n$.
\end{rem}

We will now describe the correspondence between simple complete ideals and classes of equisingular plane curves 
following the exposition given in \cite[II.5, p. 433]{S}. The class of equisingular curves $\cL$ 
corresponding to the ideal $\fa$ is defined to be the set of the analytically irreducible plane curves whose 
strict transform on $\mathcal X_n$ intersects transversely at the point $\varsigma_n$ the strict transform of any 
exceptional divisor passing through $\varsigma_n$. Note that $\cL$ is specified by a pair $(\cC,t)$, where 
$\cC\in\cL$ is a curve and $t:=n-\nu$, so that $\cL$ is the collection of all the curves equisingular to $\cC$ and 
sharing the $\nu+t$ first points of its point locus with $\cC$.

Conversely, let $\cL$ be the class of equisingular plane curves specified by the pair $(\cC,t)$. Then the corresponding 
simple ideal $\fa$ is generated by the defining equations of the elements of $\cL$. If the standard resolution 
of $\cC$ is $\bar\pi:\mathcal X_{\nu+1}\rightarrow\mathcal X_1=\Spec\ga$, then $\fa$ is 
the ideal, whose resolution is 
$\pi=\pi_n\circ\cdots\circ\pi_{\nu+1}\circ\bar\pi:\mathcal X_{n+1}\rightarrow\mathcal X_1$, where $n=\nu+t$ and 
$\pi_i:\mathcal X_{i+1}\rightarrow\mathcal X_i$ is the blow up emerging in the sequence (\ref{rloc}) corresponding to 
the point locus of $\cC$. 

For the convenience of the reader we state the following variant of \cite[Proposition 9.2.28]{Lz} adjusted to our case:

\begin{prop}\label{curveJ} 
Let $\ga$ be a two-dimensional regular local ring and let $\fa$ be a simple complete ideal in $\ga$ having the 
resolution (\ref{res}). 
Suppose that $\mathcal C\subset\Spec\ga$ is an analytically irreducible plane curve, whose strict transform intersects 
transversely at the point $\varsigma_n$ the strict transform of any exceptional divisor passing through $\varsigma_n$.
Then the multiplier ideals of the curve $\cC$ and the multiplier ideals of the ideal $\fa$ coincide in the interval 
$[0,1[$. In particular, the jumping numbers of the curve $\cC$ and the ideal $\fa$ coincide in the interval $[0,1[$.
\end{prop}

\begin{proof}
Take a non-negative rational number $c$. The multiplier ideal $\cJ(c\cdot\mathcal C)\subset\ga$ is 
defined by 
\[
\cJ(c\cdot\mathcal C):=
\Gamma\big(\mathcal X,\cO_{\mathcal X}\big(K_{\mathcal X}-\lfloor c\cdot\pi^*\mathcal C\rfloor\big)\big),
\]
where $\pi^*\mathcal C$ is the total transform of $\cC$ on $\mathcal X$. 
Since $\cC^{(n)}$ intersects transversely 
at the point $\varsigma_n$ the strict transform of any exceptional divisor passing through $\varsigma_n$, 
Lemma \ref{cf} (ii) implies $\pi^*\cC=\cC^{(n)}+\hat E_n$. Moreover, 
$\fa\cO_{\mathcal X}=O_{\mathcal X}(-\hat E_n)$, and thus
\[
\cJ( c\cdot\mathcal C)= 
\Gamma\big(\mathcal X,\cO_{\mathcal X}\big(K_{\mathcal X}
-\lfloor c\cdot \hat E_n\rfloor-\lfloor c\cdot\cC^{(n)}\rfloor\big)\big), 
\]
Suppose that $0\le c<1$. Then $\lfloor c\cdot\cC^{(n)}\rfloor$ vanish, and we obtain
\[
\cJ( c\cdot\mathcal C)= 
\Gamma\big(\mathcal X,\cO_{\mathcal X}\big(K_{\mathcal X}-\lfloor c\cdot \hat E_n\rfloor\big)\big)=\cJ(\fa^c). 
\]
\end{proof}

\begin{thm}\label{maincurve} 
Let $\ga$ be a two-dimensional regular local ring, and let $\fa$ be a simple complete ideal in $\ga$ having the 
resolution (\ref{res}).
Let $\mathcal C\subset\Spec\ga$ be an analytically irreducible plane curve, whose strict transform intersects 
transversely at the point $\varsigma_n$ the strict transform of any exceptional divisor passing through $\varsigma_n$. 
Then the set of the jumping numbers of $\mathcal C$ is 
\[
\cH_{\mathcal C}=\{c+m\mid c\in\cH_\fa\cup\{1\},0<c\le1\text{ and }m\in\mathbb N\}, 
\]
where $\cH_\fa$ denotes the set of the jumping numbers of $\fa$. 
\end{thm}

\begin{proof} 
By the periodicity of the jumping numbers for integral divisors (see e.g. \cite[Remark 1.15]{ELSV}) we know that $c>0$ 
is a jumping number of $\mathcal C$ if and only if $c+1$ is. Thus it is enough to find out the jumping numbers of 
$\cC$ in the interval $]0,1]$.

If $c<1$, then by Proposition \ref{curveJ} $\cJ(c\cdot\mathcal C)=\cJ(\fa^c)$. 
Because $\left\lfloor c\cdot\cC^{(n)}\right\rfloor=0$ for $c<1$,
we must have $\cJ(c\cdot\cC)\supsetneq\cJ(\cC)$ for $c<1$ and thus $1$ must 
be also a jumping number of $\cC$.
Hence $c\in]0,1]$ is a jumping number of $\mathcal C$, if and only if $c \in \cH_\fa$ or $c=1$.
\end{proof}

This result can be utilized in determining the jumping numbers of an arbitrary analytically irreducible plane 
curve.

\begin{cor}\label{maincurve2} 
Let $\ga$ be a two-dimensional regular local ring and let 
$\mathcal C\subset\Spec\ga$ be an analytically irreducible plane curve. Let $\fa$ be the simple complete ideal in $\ga$ 
corresponding to the class of equisingular plane curves specified by the pair $(\cC,0)$. Then the set of the 
jumping numbers of $\mathcal C$ is
\[
\cH_{\mathcal C}=\{c+m\mid c\in\cH_\fa\cup\{1\},0<c\le1\text{ and }m\in\mathbb N\}. 
\]
\end{cor}

\begin{proof} 
As we observed above, the standard resolution of $\cC$ coincides with the resolution of $\fa$. Then the claim is 
a direct consequence of Theorem \ref{maincurve}.
\end{proof}

\begin{rem}\label{pax} 
By Corollary \ref{maincurve2} we observe that the jumping numbers of $\cC$ depend only on the 
equisingularity class of $\cC$, because by Theorem \ref{main} the jumping numbers of $\fa$ less than one are totally 
determined by the multiplicities $a_1,\dots,a_{\gamma_g}$ of the point basis $(a_1,\dots,a_n)$ of $\fa$. 
\end{rem}

\begin{rem}\label{pux}
Let $\mathcal C$ be an analytically irreducible plane curve, and let 
$(a_1,\dots,a_\nu)$ be the multiplicity sequence of $\mathcal C$. Let $X_i$ denote the $i$:th row of the inverse of 
the corresponding proximity matrix, and let $\{\gamma_1,\dots,\gamma_g\}$ be the indices corresponding to the terminal 
satellites. Note that $g$ is the genus of the curve $\mathcal C$, (see e.g. \cite[Definition 3.2.1]{C1}). The 
\textit{characteristic exponents} of $\mathcal C$ are 
\[
\beta_k:=a_1+\rho_{n,\gamma_0}+\cdots+\rho_{n,\gamma_{k-1}},\text{ where }k=0,\dots,g.
\] 
It is easy to see by induction on $k$ that $\gcd\{\beta_0,\dots,\beta_{k}\}=a_{\gamma_{k}}$. Indeed, we first observe 
that $\beta_0=a_1$ and $\beta_k=\beta_{k-1}+\rho_{n,\gamma_{k-1}}$ for every $0<k\le g$. 
Assume that $\gcd\{\beta_0,\dots,\beta_{k-1}\}=a_{\gamma_{k-1}}$. Then by Proposition \ref{aiaj1}
\[
\gcd\{\beta_0,\dots,\beta_{k}\}=\gcd\{a_{\gamma_{k-1}},\rho_{n,k-1}\}=a_{\gamma_{k}}.
\]

The \textit{characteristic pairs} or \textit{Puiseux pairs} (in the case $\Bbbk=\mathbb C$) 
are the pairs of integers $(m_k,n_k)$ for $k=1,\dots,g$, where 
\[
m_k:=\frac{\gb_k}{a_{\gamma_k}}\text{ and }n_k:=x_{\gamma_k,\gamma_{k-1}}
\] 
(see, e.g., \cite[Remark 3.1.6]{C1}). We can obtain the pairs $(a_{\gamma_0},b_0),\dots,(a_{\gamma_g},b_g)$ from these 
as follows: Corollary \ref{ax} yields  
\[
m_i=x_{\gamma_i,\gamma_0}+\rho_{\gamma_i,\gamma_0}+\cdots+\rho_{\gamma_i,\gamma_{i-1}}
\] 
and 
$a_{\gamma_{i-1}}=n_i\cdots n_g$ for $i=1,\dots,g$. 
Moreover, $\rho_{\gamma_i,\gamma_{i-1}}=m_i-n_im_{i-1}$ for $i=2,\dots,g$, which further yields 
$\rho_{n,\gamma_{i-1}}=(n_{i+1}\cdots n_g)(m_i-n_im_{i-1})$. Then also 
$a_{\gamma_{i-1}}/a_{\gamma_k}=n_i\cdots n_k$ for  $i=1,\dots,k+1$. Proposition \ref{key1} implies that 
\[
a_{\gamma_k}b_k=I\cdot I^{\le\gamma_{k+1}}=a_{\gamma_0}^2+a_{\gamma_0}\rho_{n,\gamma_0}+\cdots+a_{\gamma_k}\rho_{n,\gamma_k}. 
\]
Writing $\varphi_1=m_1$ and $\varphi_i=m_i-n_im_{i-1}$, we then get 
\[
a_{\gamma_k}=n_{k+1}\cdots n_g\text{ and }b_k=\sum_{i=1}^{k+1}(n_{i+1}\cdots n_g)(n_i\cdots n_k)\varphi_i
\]
for every $k=1,\dots,g$. 
\end{rem}

\begin{thm}\label{curveJh} 
Let $\ga$ be a two-dimensional regular local ring. The jumping numbers of an analytically irreducible plane curve 
$\mathcal C\subset\Spec\ga$ less than one determine the equisingularity class of $\cC$.
\end{thm}

\begin{proof}
Take the ideal $\fa$ corresponding to the class of equisingular curves specified by the pair $(\cC,0)$. Then the point 
basis of $\fa$ is the same as the multiplicity sequence $(a_1,\dots,a_{\gamma_g})$ of $\cC$. By Theorem 
\ref{maincurve} we know that the the jumping numbers of $\fa$, which are less than one, coincide with those of the 
curve $\cC$. By Theorem \ref{1maincor} they determine the sequence of multiplicities $(a_1,\dots,a_{\gamma_g})$. 
\end{proof}

In Theorem \ref{maincurve} we saw how to obtain the jumping numbers of the analytically irreducible plane curve 
determined by a general element of a simple complete ideal from the jumping numbers of the ideal. In the next theorem 
we consider the converse situation.

\begin{thm}\label{paxi}
Let $\ga$ be a two-dimensional regular local ring, and let $\fa$ be a simple complete ideal in $\ga$ having the 
resolution (\ref{res}).
Let $\mathcal C\subset\Spec\ga$ be an analytically irreducible plane curve, whose strict transform intersects 
transversely at the point $\varsigma_n$ the strict transform of any exceptional divisor passing through $\varsigma_n$. 
Then the set $\cH_\fa$ of the jumping numbers of $\fa$ is
\[
\begin{array}{l}
\cH_\fa=
\left(\cH_\cC\smallsetminus\{1\}\right)
\displaystyle\cup
\left\{1+\frac{k+1}{v(\fa)}\hspace{2pt}\middle|\hspace{2pt} k\in\mathbb N\right\},
\end{array}
\]
where $v=v_n$ denotes the divisorial valuation associated to $\fa$.
\end{thm}

\begin{proof}
As before, let $I=(a_1,\dots,a_n)$ denote the point basis of $\fa$ and let $\cH_\fa=H_0\cup\cdots\cup H_{g^*}$ be the 
set of the jumping numbers of $\fa$ where the subsets $H_\nu$ are as in Theorem \ref{main}.
Recall that $v(\fa)=I^2$ (see Remark \ref{val}). By Proposition \ref{höm}  
\[
\left\{1+\frac{k+1}{I^2}\hspace{2pt}\middle|\hspace{2pt} k\in\mathbb N\right\}=\left\{c\in H_{g^*}\mid c>1\right\}.
\]
Proposition \ref{curveJ} shows that $0<c<1$ is a jumping number of $\fa$, if and only if $c\in\cH_\cC$. Moreover, 
Proposition \ref{1ei} says that $1\notin\cH_\fa$. In order to prove the claim, it is thereby enough to show that every 
element of $H_\nu$ greater than one is in $\cH_\cC$ for any $\nu<g^*$, and that every element of $\cH_\cC$ gerater than 
one is in $\cH_\fa$.

Suppose that $\xi>1$ and $\xi\in H_\nu$ for some $\nu<g^*$. Theorem \ref{main} implies that 
$\xi=c+m/a_{\gamma_{\nu+1}}$ for some $c\in H_\nu$ and $m\in\mathbb N$, where $0<c\le a_{\gamma_{\nu+1}}^{-1}$. By 
Proposition \ref{1ei} we even know that $c<a_{\gamma_{\nu+1}}^{-1}$. Now $m>(1-c)a_{\gamma_{\nu+1}}$ yields 
$m\ge a_{\gamma_{\nu+1}}$. Write $m=m'a_{\gamma_{\nu+1}}+m''$, where $m'\in\mathbb N$ and $m''<a_{\gamma_{\nu+1}}$. Set 
$c':=c+m''/a_{\gamma_{\nu+1}}$. Then $c'\in H_\nu$ and $0<c'<1$. By Proposition \ref{curveJ} $c'\in\cH_\cC$. As 
$\xi=c'+m'$, it follows that $\xi\in\cH_\cC$.

Take $\xi\in\cH_\cC$ so that $\xi>1$. By Proposition \ref{höm} the case is clear, if $\xi$ is an integer. Let us then 
assume that $\xi$ is not an integer. Proposition \ref{maincurve} implies that $\xi=c+m$, where $m$ is a positive 
integer and $c\in H_\nu$ with $0<c<1$ for some $\nu\in\{0,\dots,g^*\}$. It follows from Theorem \ref{main} that 
$c=c'+m'/a_{\gamma_{\nu+1}}$ for some $m'\in\mathbb N$ and for some $c'\in H_\nu$ satisfying 
$c'<a_{\gamma_{\nu+1}}^{-1}$. Thus $\xi=c'+(a_{\gamma_{\nu+1}}m+m')/a_{\gamma_{\nu+1}}$, and $\xi\in\cH_\fa$ according 
to Theorem \ref{main}. This completes the proof.
\end{proof}

\begin{rem}\label{pox}
In the setting of Theorem \ref{paxi}, $\fa$ is the ideal corresponding to the pair $(\cC,t)$, where $t=n-\gamma_g$.
The jumping numbers of the curve $\cC$ together with the integer $t$, or equivalently the integer $n$, then determine 
the jumping numbers of $\fa$. Indeed, according to Theorem \ref{curveJh} the integer $\gamma_g$ as well as the 
entire multiplicity sequence $(a_1,\dots,a_{\gamma_g})$ of $\cC$ can be obtained from the jumping numbers of the curve 
$\cC$. Because $a_{\gamma_g}=\cdots=a_n=1$ by Proposition \ref{taudit}, we then get 
$v(\fa)=a_1^2+\cdots+a_{\gamma_g}^2+n-\gamma_g$. 
\end{rem}

\section*{Acknowledgements}
I am greatly indebted to Eero Hyry who suggested and introduced this topic to me and who patiently helped me with 
his advice during the preparation of this manuscript. I want to express my gratitude to Ana J. Reguera L\'opez and Karen E. Smith for 
their critical comments and valuable suggestions for improvements. I am also grateful for the financial support from 
the Academy of Finland, project 48556.

\end{document}